\documentclass[12pt,twoside,a4paper]{article}
\usepackage{fullpage}
\usepackage{amssymb}
\usepackage{amsmath}
\usepackage{epsfig}
\usepackage{psfrag}

\def \bR{\mathbb R}
\def \R{\mathbb R}
\def \P{\mathbb P}

\def \Z{\mathbb Z}
\def \Id{{\rm Id}}

\def \cM{\mathcal M}
\def \cO{\mathcal O}
\def \TS{{\rm Terseg}}

\def \tz0{\tau^Z_0}

\def \tzp0{\tau^{Z'}_0}

\newcommand{\bbE}{{\mathbb E}}
\newcommand{\bbP}{{\mathbb P}}

\def \ud{\frac{1}{2}}
\def \l{\lambda}
\def \8{\infty}
\def \1{\bf 1}

\def \a{\alpha}

\def \I0T{I_{0,T}}
\def \J0T{J_{0,T}}

\def\qed{\hfill\rule{.2cm}{.2cm}\par\medskip\par\relax}
\def \Lr{L^{\rm ref}}

\newcommand{\mysection}{\setcounter{equation}{0} \section}

{\markboth{\centerline{COMETS, DELARUE AND SCHOTT}}{\centerline{DISTRIBUTED ALGORITHMS}}}

\newtheorem{thm}{Theorem}[section]
\newtheorem{prop}[thm]{Proposition}
\newtheorem{cor}[thm]{Corollary}
\newtheorem{lem}[thm]{Lemma}

\begin{document}

\parindent=10pt

\thispagestyle{empty}

\title{Large Deviations Analysis for
Distributed Algorithms in an Ergodic  Markovian Environment}

\author{Francis Comets${}^{(a),}$\footnote{comets@math.jussieu.fr; http://www.proba.jussieu.fr/pageperso/comets} \hspace*{.05cm},
Fran{\c c}ois Delarue${}^{(a),}$\footnote{delarue@math.jussieu.fr; http://www.math.jussieu.fr/$\sim$delarue} \hspace*{.05cm} and
Ren{\'e} Schott${}^{(b),}$\footnote{schott@loria.fr; http://www.loria.fr/$\sim$schott}
 \\ \\
(a) Laboratoire de Probabilit{\'e}s et Mod{\`e}les
Al{\'e}atoires,\\ Universit{\'e} Paris 7, UFR de Math{\'e}matiques, Case 7012, \\ 2, Place Jussieu, 75251 Paris Cedex 05 - France.
\\
\\
(b) IECN and LORIA, Universit{\'e} Henri Poincar{\'e}-Nancy 1, \\ 54506 Vandoeuvre-l{\`e}s-Nancy - France.}
\maketitle

\begin{abstract}
We provide a large deviations analysis of deadlock phenomena occurring in distributed systems sharing common resources. 
In our model transition probabilities of resource allocation and deallocation are time and space dependent. The process is driven by an ergodic 
Markov chain and is reflected on the boundary of the $d$-dimensional cube. In the large resource limit, we prove Freidlin-Wentzell estimates, 
we study the asymptotic of the deadlock time and we show that the quasi-potential is a viscosity solution of a Hamilton-Jacobi equation
with a Neumann boundary condition.
We give a complete analysis of the colliding 2-stacks problem and show an example where the system has a stable attractor which is a limit cycle.

\end{abstract}

\footnotesize
\noindent{\bf Short Title:} Distributed Algorithms in an Ergodic Environment
 
\noindent{\bf Key words and phrases:} Large deviations, 
distributed algorithm, averaging principle,
Hamilton-Jacobi equation, viscosity solution

\noindent{\bf AMS subject classifications:} Primary 60K37; secondary
60F10, 60J10

\tableofcontents

\normalsize

\section{Introduction}
Distributed algorithms are related to resource sharing problems. Colliding stacks problems and the banker algorithm are among the examples which have attracted large 
interest over the last decades in the context of deadlock prevention on multiprocessor systems. Knuth \cite{knuth}, Yao \cite{yao}, Flajolet \cite{flajolet}, Louchard, 
Schott et al. \cite{louchard1,louchard2,louchard3} have provided combinatorial or probabilistic analysis of these algorithms in the $2$-dimensional case under the 
assumption that transition probabilities (of allocation or deallocation) are constant. Maier \cite{Maier} proposed a large deviations analysis of colliding stacks for
the more difficult case where the transition probabilities are non-trivially state-dependent. More recently Guillotin-Plantard and Schott \cite{GS1,GS2} analyzed a model 
of exhaustion of shared resources where allocation and deallocation requests are modeled by time-dependent dynamic random walks. In \cite{CDS1}, the present authors provided 
a probabilistic analysis of the $d$-dimensional banker algorithm when transition probabilities evolve, as time goes by, along the trajectory of an ergodic Markovian environment, whereas the spatial parameter just acts on long runs. The analysis in \cite{CDS1} relies on techniques from stochastic homogenization theory. In this paper, 
we consider a similar dynamics, but in a stable regime instead of a neutral regime as in our previous paper, and we provide an original large deviations analysis in the framework of Freidlin-Wentzell theory. 
Given the environment, the process of interest is a Markov process depending
on the number $m$ of available resource, with smaller and more frequent jumps 
as $m \to \8$, see (\ref{eq:last}).
A number of  monographs and papers have been written on this theory: \cite{FK}, \cite{FV} and \cite{GV} for random environment, \cite{dupuis 87}, \cite{Dupuis-Ellis}, \cite{ignatiouk01} and
 \cite{ignatiouk} for reflected processes,
\cite{AZ}, \cite{DZ} and \cite{OV} for homogeneous Markov processes. However, our framework, including both reflections on the boundary 
and averaging on the Markovian environment, is not covered by the current literature, and we establish here the large deviations principle.
We prove that the time of resource exhaustion then grows exponentially with the size of the system -- instead of polynomially in the neutral regime of \cite{CDS1}-- 
and has exponential law as limit distribution. Then, we study  the quasi-potential, which solves, according to general wisdom, some Hamilton-Jacobi equation: in view of the reflection on the hypercube, which boundary is non-regular, 
we prove this fact in the framework of  viscosity solution, and study the optimal paths (so-called instantons).


We investigate in details a particular situation introduced in a beautiful paper of Maier \cite{Maier}, where the motion in each direction depends on the corresponding coordinate only, with the additional dependence in the Markovian environment. In fact, we discover the quasi-potential by observing that the discrete process has an invariant measure, for which we study the large deviations properties. We can then use the characterization in terms of Hamilton-Jacobi equation to bypass  the Hamiltonian mechanics approach of \cite{Maier}.
For the deadlock phenomenon, we finally obtain a (even more) complete
picture (after an even shorter work). To the best of our knowledge, 
this is first such analysis developed for space-time inhomogeneous distributed
algorithms.
\\

The organization of this paper is as follows: we discuss our probabilistic model in Section 2. In Section 3 we prove a Large Deviations Principle. Deadlock phenomenon analysis is done rigorously with much details in Section 4. In Section 5 we illustrate with the two-stacks model. In Section 6 we work out an example where the system has (in the large scale resource limit $m\rightarrow\infty$) a stable attractor which is a limit cycle. Some technical proofs of results stated in Section 4 are deferred to Appendix (Section 7).

\section{The Model}

The environment is given by  a Markov chain  $(\xi_n)_{n \geq 0}$
defined on $(\Omega, {\cal A},  {\mathbb P})$
with values in a finite space $E$, $N= |E|$.
We denote by  $P$ its transition matrix,
$P(k,\ell) = {\mathbb P}(\xi_{n+1}=\ell|\xi_{n}=k)$ for $ k,\ell
\in E$.

The steps of the walker take place in the set
$${\mathcal V} = \{e_1,-e_1,\dots,e_d,-e_d\}\;,$$
where $(e_i)_{1 \leq i \leq d}$ denotes the canonical basis of ${\mathbb Z}^d$, and are reflected along the boundary of
the hypercube $[0,m]^d$, for a large integer $m$.

Following \cite{CDS1}, we first discuss the dynamics of the walk in the non-reflected setting.
The displacement of the walker has then law $(p(s/m,i,v); v
\in  {\mathcal V})$, when located at $s$ and when the environment is
$i$. To obtain a stochastic representation -- which is, in contrast to
\cite{CDS1}, needed here -- , we are also given, on $(\Omega,{\mathcal A},\P)$,
a sequence $(U_n)_{n \geq 0}$ of independent and uniformly distributed
random variables on $(0,1)$, independent of the family $(\xi_n)_{n \geq 0}$.
Denoting by $f : (s,i,v) \in \R^d \times E \times (0,1) \rightarrow {\mathcal V}$
the inverse of the cumulative distribution function of $(p(s,i,u))_{u \in {\mathcal V}}$
(for an arbitrary order on ${\mathcal V}$), we have
$$
{\mathbb P}\{ f(s, i, U_{n})=v\}=p(s,i,v) \;,
$$
so that the position of the walker can be defined recursively by
\begin{equation}
\label{walk}
S_{n+1}=S_n + f(S_n/m, \xi_n, U_{n}).
\end{equation}
The reflected walk is obtained by symmetry with respect to the faces of the hypercube. Denoting by $\Id$ the identity mapping on $\R^d$ and
by $\Pi$ (resp. $\Pi^{(m)}$) the projection on the hypercube $[0,1]^d$ (resp. $[0,m]^d$), we define recursively the position of the walker by
\begin{equation}
\label{reflected_walk}
\begin{split}
X_{n+1} &= (2 \Pi^{(m)} - \Id)\bigl( X_n + f(X_n/m,\xi_n,U_n) \bigr)
\\
&=
m (2\Pi-\Id) \bigl( X_n/m + (1/m)f(X_n/m,\xi_n,U_n) \bigr).
\end{split}
\end{equation}
When $X_n$ is on the boundary and $X_n + f(X_n/m,\xi_n,U_n)$ is outside the hypercube, $X_{n+1}$ is the symmetric point of $X_{n} + f(X_n,\xi_n,U_n)$ with respect to the face
containing $X_n$ and orthogonal to $f(X_n/m,\xi_n,U_n)$, i.e. $X_{n+1} = X_n - f(X_n/m,\xi_n,U_n)$.
The kernel $q$ of the walk $(X_n)_{n \geq 1}$ has the following form. When located at $x \in (0,m)^d \cap \Z^d$ and when the environment is $i$, the jump of the walker has
law $(q(x/m,i,v) =
p(x/m,i,v);v \in {\mathcal V})$. On the boundary, the reflection rules may be expressed as follows: if $x_{\ell}/m=1$
($x_{\ell}$ is the $\ell$th coordinate of
$x$), $q(x/m,i,e_{\ell})=0$
and $q(x/m,i,-e_{\ell})=p(x/m,i,e_{\ell})+p(x/m,i,-e_{\ell})$
; if $x_{\ell}=0$, $q(x/m,i,-e_{\ell})=0$
and $q(x/m,i,-e_{\ell})=p(x/m,i,e_{\ell})+p(x/m,i,-e_{\ell})$.

We could choose another reflection rule by setting $X_{n+1} = \Pi^{(m)} ( X_n + f(X_n/m,\xi_n,U_n))$. Such a choice wouldn't change anything to the proofs given in the paper,
except the proof of Theorem \ref{th:tloi} which uses the fact that the steps of $(X_n)_{n \geq 0}$ are always non-zero.

Following \eqref{walk}, we can write
\begin{equation}
\label{eq:form_g}
X_{n+1}=X_n + g(X_n/m,\xi_n,U_n) \; ,
\end{equation}
where $g : (x,i,v) \in \R^d \times E \times (0,1) \rightarrow {\mathcal V}$ is
the inverse of the cumulative distribution function of $(q(x,i,u))_{u \in {\mathcal V}}$. Of course, $g(x,i,v) = f(x,i,v)$ for $x \in (0,1)^d$. On the boundary,
$x_{\ell} =1 \Rightarrow g_{\ell}(x,i,v) \leq 0$ and $x_{\ell}=0
\Rightarrow g_{\ell}(x,i,v) \geq 0$.

The process $(\xi_n,X_n)_{n \geq 0}$ is a Markov chain with
transition probabilities
\begin{equation*}
 \ {\mathbb P}
\bigl\{\xi_{n+1}=k,X_{n+1}=v+X_n |{\mathcal F}^{\xi,X}_n\bigr\}
= P(\xi_{n},k) q(X_n/m,\xi_{n},v),
\end{equation*}
where ${\mathcal F}^{\xi,X}_n = \sigma
\{\xi_0,\dots,\xi_n,X_0,\dots,X_n\}$.
In particular, for $x \in {\mathbb Z}^d$ and for a probability measure
$\nu$ on $E$, we can write ${\mathbb P}^{\nu}_{x/m}$ to indicate that the
chain starts under the measure $\delta_x \otimes \nu$. In many cases,
we just write ${\mathbb P}_{x/m}$ (resp. ${\mathbb P}^{\nu}$): this means that the law of the
environment (resp. of the walker) is arbitrary. And, of course, the
notation $\P$ means that both the initial conditions of the walker and
of the environment are arbitrary.

\subsection{Main Assumptions}
\textit{In the whole paper,
$\langle \cdot,\cdot  \rangle$ and
 $|\cdot|$ stand for the Euclidean scalar product and the Euclidean norm in $\R^d$. The symbols
$|\cdot|_1$ and $|\cdot|_{\infty}$ denote the standard $\ell^1$ and $\ell^{\infty}$ norms in $\R^d$.}

\textit{From a purely practical point of view, the values of $p(x,i,v)$ for $x$ outside the hypercube $[0,1]^d$ are totally useless. In the sequel, we refer, for pedagogical
reasons, to the non-reflected walk: in such cases, we need $p(x,i,v)$ to be defined for all $x \in \R^d$. This is the reason why the variable $x$ lies in $\R^d$ in the
following assumptions.}

In formulas \eqref{walk} and \eqref{reflected_walk}, the
division by $m$ indicates that the dependence of the transition kernel
on the position of the walker takes place at scale $m$. For large  $m$,
the space dependence is mild, since we will assume all through the paper
 the following {\it smoothness} property:
\vspace{5pt}
\\
{\bf Assumption (A.1).} There exists a finite constant $K$ such that
$|p(x,i,v)-p(y,i,v)| \leq K |x-y|, \ x,y \in  {\mathbb R}^d,
i \in E, v \in   {\mathcal V}$.

\medskip

For technical reasons, which are explained in the paper, we impose the following {\it ellipticity} condition:

\medskip

\noindent
{\bf Assumption (A.2).} For all $x \in \R^d$, $i \in E$ and $v \in {\mathcal V}$, $p(x,i,v)>0$. By continuity,
$\inf \{p(x,i,v); \, x \in [0,1]^d, \, i \in E, \, v \in {\mathcal V}\}>0$.

\medskip

We also assume the environment to be {\it ergodic} and to obey the large
deviations principle for Markov
chains. We thus impose the following sufficient conditions:
\vspace{5pt}
\\
{\bf Assumption (A.3).} The matrix $P$ is irreducible on $E$. Its unique invariant probability
measure is denoted by $\mu$.
\vspace{5pt}
\\
In particular, the following vector-valued function is smooth:
\begin{equation}
\label{def:barf}
\bar f (x) = {\mathbb E}^\mu f(x,\xi,U)
 = \sum_{v \in  {\mathcal V}} v \; {\mathbb E}^\mu p(x,\xi,v)
\;,\quad x \in  {\mathbb R}^d\;;
\end{equation}
the above expectations are taken over independent variables $\xi,U$,
where $\xi$ has the distribution $\mu$ and $U$
is uniformly distributed on [0,1].

For the deadlock time analysis, another assumption will be necessary (see {\bf (A.4)} in Section \ref{sec-deadlock}).

\subsection{Continuous Counterpart and Skorohod Problem}
\label{subsection:skorohod}

Because of the reflection phenomenon, we briefly recall what the Skorohod problem is (we
refer to \cite{lions sznitman 1984} for a complete overview of the subject). For each continuous mapping $w : t \in [0,+\infty) \mapsto w_t \in \R^d$,
with $w_0 \in [0,1]^d$, there exists a unique continuous mapping $t \in [0,+\infty) \mapsto (x_t,k_t) \in [0,1]^d \times \R^d$, with $k$ of bounded variation
on any bounded sets, such that:
\begin{equation}
\label{Skorokhod}
\forall t \geq 0 \;, \quad w_t = x_t + k_t \;, \quad  k_t = \int_0^t n_s d|k|_s \;, \quad |k|_t = \int_0^t {\mathbf 1}_{\{x_s \in \partial [0,1]^d\}} d|k|_s \; ,
\end{equation}
where $n_s \in {\mathcal N}(x_s)$, ${\mathcal N}(x)$ denoting for $x \in \partial [0,1]^d$ the set of unit outward normals to $\partial [0,1]^d$ at $x$,
 that is
$$
{\mathcal N}(x)=\{ v \in \R^d: \, |v|=1, \, v_{\ell}=0 {\rm \ if\ } x_\ell \in (0,1), \,
  v_{\ell}\leq 0 {\rm \ if\ } x_\ell=0,  \, v_{\ell}\geq 0 {\rm \ if\ } x_\ell=1
\}.
$$
When $x$ is in the relative interior of a face of
the hypercube, ${\mathcal N}(x)$ is obviously empty.

It can be proved (see again \cite{lions sznitman 1984}) that, for every $T >0$, the mapping $\Psi : (w_t)_{0 \leq t \leq T} \mapsto
(x_t)_{0 \leq t \leq T}$ is continuous from ${\mathcal C}_{[0,1]^d}([0,T];\R^d)$ into itself with respect to the supremum norm (it is even $1/2$-H\"older continuous on compact
subsets of ${\mathcal C}_{[0,1]^d}([0,T];\R^d)$); here and below,
${\mathcal C}_A([0,T];\R^d)$ denotes the set of continuous functions from $[0,T]$
to $\R^d$ with an initial datum in $A$. Moreover, if $w$ is absolutely continuous, then $x$ and $k$ are also absolutely continuous
(see \cite[Theorem 2.2]{lions sznitman 1984}).

Equation \eqref{reflected_walk} corresponds to a Euler scheme for a Reflected Differential Equation (RDE in short).
An RDE is an ordinary differential equation, but driven by a pushing process $k$ as in \eqref{Skorokhod}.
For a given initial condition $x_0 \in [0,1]^d$ and a given jointly measurable and
$x$-Lipschitz continuous mapping $b : \R_+ \times [0,1]^d \rightarrow \R^d$, the RDE
\begin{equation}
\label{RDE}
\forall T > 0, \ (x_t)_{0 \leq t \leq T} = \Psi \biggl[ \biggl( x_0+
\int_0^t b(s,x_s) ds \biggr)_{0 \leq t \leq T} \biggr] \; ,
\end{equation}
admits a unique solution (see again \cite{lions sznitman 1984}). This solution satisfies the equation
\begin{equation*}
\forall t \geq 0, \ x_t = x_0 + \int_0^t b(s,x_s)ds - k_t \; ,
\end{equation*}
with $k$ as in \eqref{Skorokhod}. In this case, $x$ and $k$ are absolutely continuous.

Reflected equations driven by Lipschitz continuous coefficients are stable. By \cite[Lemma 3.1]{lions sznitman 1984}, we can prove that for every $T>0$, there exists a constant
$C_T \geq 0$, such that, for any $x_0,y_0 \in [0,1]^d$, the solutions $(x_t)_{0 \leq t \leq T}$ and $(y_t)_{0 \leq t \leq T}$ to \eqref{RDE} with $x_0$ and $y_0$ as initial
conditions satisfy $\sup_{0 \leq t \leq T}|x_t-y_t| \leq C_T |x_0-y_0|$.

When $b(s,x)=\bar{f}(x)$, we denote by $(\chi_t^{x_0})_{t \geq 0}$ the unique solution to the averaged reflected differential equation
\begin{equation}
\label{eq_lim_ref}
\forall t \geq 0, \ x_t = x_0 + \int_0^t \bar{f}(x_s)ds - k_t.
\end{equation}

\mysection{Large Deviations Principle}

We now denote the process $X$ by $X^{(m)}$ to indicate the dependence on the parameter $m$. In what follows, we
investigate an interpolated version of the rescaled process $(m^{-1} X_{\lfloor m t \rfloor}^{(m)})_{t \geq 0}$, namely
\begin{eqnarray} \nonumber
\bar{X}^{(m)}_t &=& (2 \Pi - {\rm Id})\bigl(  m^{-1}
X_{\lfloor m t \rfloor}^{(m)} + (t- m^{-1} \lfloor mt \rfloor) f(m^{-1} X_{\lfloor m t \rfloor}^{(m)},\xi_{\lfloor mt \rfloor},U_{\lfloor mt \rfloor}) \bigr)
\;, \quad t \geq 0 \\ \label{eq:last}
&=&  m^{-1}
X_{\lfloor m t \rfloor}^{(m)} + (t- m^{-1} \lfloor mt \rfloor) g(m^{-1} X_{\lfloor m t \rfloor}^{(m)},\xi_{\lfloor mt \rfloor},U_{\lfloor mt \rfloor}) 
\;, \quad t \geq 0
\end{eqnarray}
We note that the hyperbolic scaling is different from the diffusive
scaling in \cite{CDS1}. The process $(\bar{X}^{(m)}_t)_{t \geq 0}$ is continuous and $\bar{X}^{(m)}_{k/m} = X_k^{(m)}$ for any integer $k \in {\mathbb N}$.
\subsection{Heuristics for the Non-reflected Walk}
We first look, for pedagogical reasons, at the non-reflected case. We thus consider
 $$\bar{S}^{(m)}_t = m^{-1}
S_{\lfloor m t \rfloor}^{(m)} + (t - m^{-1} \lfloor mt \rfloor) f(m^{-1} S_{\lfloor mt \rfloor},\xi_{\lfloor mt \rfloor},U_{\lfloor mt \rfloor})
\;,\quad {t \geq 0} \; .$$ (As for $X$, we indicate the dependence on $m$ in $S$.)
In light of Assumptions {\bf (A.1--3)}, we expect
the global effect of the environment process $(\xi_n)_{n \geq 0}$ to
reduce for large time
to a deterministic one. More precisely,
if the initial position is such that
$\bar{S}^{(m)}_0 \to x$ as $m \to \8$, we expect
$(\bar{S}^{(m)})_{m \geq 1}$ to converge in probability, uniformly on compact sets, to the solution  $x_\cdot$
of the (averaged)
ordinary differential equation
\begin{equation}
\label{eq:ODE}
\dot x_t = \bar f (x_t)\;, \quad x_0=x\;,
\end{equation}
that is $\lim_{m \to \8} \rho_{0,T}(\bar{S}^{(m)}_\cdot , x_\cdot) =0$ for all $T>0$, where
 $\rho_{0,T}(\phi, \psi) =
\sup\{ |\phi_t-\psi_t|; \, t \in [0,T]\}$ denotes the distance in supremum norm
on the space ${\mathcal C}([0,T];\R^d)$ of continuous functions from $[0,T]$ into $\R^d$.

Loosely speaking, the Large Deviations Principle (LDP in short) for $(\bar{S}^{(m)}_{\cdot})_{m \geq 1}$ follows from the Freidlin and Wentzell theory
\cite[Chapter 7]{FV}, or at least from a variant of it
as explained below. The idea is the following.
The irreducible Markov chain $(\xi_n)_{n \geq 0}$ with a finite state space
obeys a LDP (see \cite[Theorem 3.1.2, Exercise 3.1.4]{DZ}). In particular,
 the function $H$ defined for  $x, \a \in \bR^d$ by
\begin{eqnarray} \label{def:H}
H(x,\a)&=& \lim_{n \to \8} \frac{1}{n} \ln \bbE^i
\exp \big\langle \a , \sum_{k=1}^n f(x, \xi_k, U_k) \rangle
\\ \nonumber
&=& \lim_{n \to \8} \frac{1}{n} \ln \bbE^i \prod_{k=1}^n
\left[ \sum_{v \in  {\mathcal V}} e^{ \langle \a, v \rangle} p(x,\xi_k,v) \right]\;,
\end{eqnarray}
exists and is independent of the starting point $\xi_0=i \in E$.
Here, $ \bbE^i$ denotes expectation over $(\xi_k,U_k)$ starting with
 $\xi_0=i$, and the last equality is a direct integration on the
 i.i.d. sequence $(U_n)_{n \geq 1}$.
From assumption {\bf (A.1)} and finiteness of $E$,
the limit is uniform in $x \;, \alpha$ on compact
subsets of $ {\mathbb R}^d$ and in $i \in E$.

In fact, $H(x,\a)$ is equal to the logarithm of the
Perron-Frobenius eigenvalue
(e.g., \cite[Theorem 3.1.1, Exercise 3.1.4]{DZ}) of the matrix
\begin{equation}
  \label{eq:perron}
Q(x,\a) = \big[P(i,j) \bbE e^{\langle \a, f(x,i,U) \rangle}\big]_{(i,j) \in E\times E}.
\end{equation}
Since the entries of the
above matrix are regular
and the leading eigenvalue is simple,
$H$ is continuous in $x$
and infinitely differentiable in $\a$.
For $x, v \in \bR^d$, the Legendre transform of $H(x,\cdot)$
\begin{equation}
  \label{eq:defL}
L(x,v)= \sup \{ \langle \a , v \rangle - H(x,\a); \, \a \in \bR^d\}
\end{equation}
is non-negative and  convex in $v$. It is even strictly convex, in view
of the differentiability of $H(x,\cdot)$ (see \cite[Chapter 5, (1.8)]{FV}).
In particular for all $x \in \R^d$, $\bar{f}(x)= \nabla_\alpha  H(x,0)$ 
is the unique zero of $L(x,\cdot)$.
Since $|H(x,\alpha)| \leq |\alpha|$ for all $x, \, \alpha \in \R^d$, we have $(v \in \R^d, \, |v|>1) \Rightarrow L(x,v) = + \infty$.

In some sense, the convergence in \eqref{def:H} corresponds to \cite[Chapter 7, Lemma 4.3]{FV}. By the regularity of $H$, we expect \cite[Chapter 7, Theorem 4.1]{FV} to hold in our
framework. For $x \in \R^d$ and for a sequence $(x_m)_{m \geq 1}$ converging towards $x$, with $mx_m \in \Z^d$ for all $m \geq 1$, we expect $(\bar{S}^{(m)}_{\cdot})_{m \geq 1}$
to satisfy a LDP with $m$ as normalizing coefficient and with the following action functional
\begin{equation*}
\begin{split}
\I0T (\phi) &= \int_0^T L(\phi_s, \dot \phi_s)ds
\quad {\rm if} \ \phi_0=x \ {\rm and} \ \phi \ {\rm is \ absolutely \ continuous}
\\
&=  \8 \ {\rm otherwise}.
\end{split}
\end{equation*}

\subsection{Large Deviations Principle for the Reflected Walk}
We now prove the LDP for the reflected walk. Generally speaking, it follows from the LDP for the process $\bar{S}^{(m)}_{\cdot}$ and from the
contraction principle (see e.g. \cite[Theorem 4.2.1, p. 126]{DZ}). For this reason, we have first to make rigorous the previous paragraph. In what follows, we will see that the
theory of Freidlin and Wentzell cannot be applied in a straight way. Indeed, for our own purpose (see the next section for the application to the deadlock time problem), we are
seeking for uniform large deviations bounds with respect to the starting point. In \cite{FV}, the authors obtain uniform bounds for systems driven by a Lipschitz continuous
field $f$. Since our own $f$ takes its values in a discrete set, it cannot be continuous.

To overcome the lack of regularity of $f$, we follow the approach of Dupuis \cite{dupuis:88}. The idea is to use a ``uniform''
version of the G\"artner-Ellis theorem to obtain
uniform bounds (see \cite[Theorem 2.3.6]{DZ} for the original version of the G\"artner-Ellis theorem). More precisely, we follow Section 5 in \cite{dupuis:88}. In this framework,
we emphasize that $(
\bar{X}_t^{(m)})_{t \geq 0}$ is 1-Lipschitz continuous (in time)
and adapted to the filtration $({\mathcal G}_t^{(m)} = \sigma(\xi_k,U_k, \, k \leq \lfloor tm \rfloor))_{t \geq 0}$.
We consider the non-projected and projected versions
\begin{equation*}
\begin{split}
&Y_t^{(m)} = Y_{k/m}^{(m)} + (t-k/m) f(\bar{X}_{k/m}^{(m)},\xi_k,U_k) \; ,
\\
&Z_t^{(m)} = \Pi \bigl( Z_{k/m}^{(m)} + (t-k/m) f(\bar{X}_{k/m}^{(m)},\xi_k,U_k) \bigr) \; ,
\end{split}
\end{equation*}
with $Y_0^{(m)} = Z_0^{(m)} =\bar{X}^{(m)}_0$.
 They are also 1-Lipschitz continuous in time and adapted to $({\mathcal G}_t^{(m)})_{t \geq 0}$. We let the reader check that, for all $t \geq 0$, $|Z_t^{(m)}-\bar{X}_t^{(m)}|_{\infty} \leq 1/m$.
Moreover, for $t \in [k/m,(k+1)/m)$,
 \begin{equation*}
 \begin{split}
&{Z}_t^{(m)} = {Z}_{k/m}^{(m)} + Y_t^{(m)} - Y_{k/m}^{(m)} - \bigl[ K_t^{(m)} - K_{k/m}^{(m)} \bigr],
\\
&K_{t}^{(m)} - K_{k/m}^{(m)} =
 \bigl({Z}_{k/m}^{(m)} + Y_t^{(m)} - Y_{k/m}^{(m)} \bigr)
- \Pi\bigl({Z}_{k/m}^{(m)} + Y_t^{(m)} - Y_{k/m}^{(m)} \bigr),
 \end{split}
\end{equation*}
with $K_0^{(m)}=0$. Summing over $k$, we have
$
{Z}_t^{(m)} = Y_t^{(m)} - K_t^{(m)}$.
The process $K^{(m)}$ is of bounded variation on compact sets.
If ${Z}_{k/m}^{(m)} \in (0,1)^d$, $K_t^{(m)} - K_{k/m}^{(m)} = 0$ for $k/m$ $\leq t <(k+1)/m$.
Otherwise, ${Z}_{k/m}^{(m)} \in \partial [0,1]^d$ and $K_t^{(m)} - K_{k/m}^{(m)} \in \R_+ {\mathcal N}({Z}_{k/m}^{(m)}) = \R_+ {\mathcal N}({Z}_t^{(m)})$.
We deduce that ${Z}^{(m)}$ is nothing but $\Psi(Y^{(m)})$ ($\Psi$ being the Skorohod mapping).
Since $Z^{(m)}$ and $\bar{X}^{(m)}$ are close,
it is sufficient to establish the LDP for $Y^{(m)}$ and to conclude by the contraction principle,.

The LDP for $Y^{(m)}$ follows from \cite[Theorem 3.2]{dupuis:88} (up to a slight modification of the proof). Indeed, we can write
\begin{equation}
\label{increment:Y}
Y_t^{(m)} = Y_{k/m}^{(m)} + (t-k/m) f\bigl(\Psi(Y^{(m)})_{k/m}+\varepsilon_{k/m}^{(m)},\xi_k,U_k \bigr) \; , \; k/m \leq t < (k+1)/m \; ,
\end{equation}
with $|\varepsilon_{k/m}^{(m)}|_{\infty} \leq 1/m$.
This form is the analogue of the writing obtained in \cite[p. 1532]{dupuis:88} for $\tilde{X}_n^{\varepsilon}$.
In \eqref{increment:Y}, we can choose an arbitrary initial condition $y \in [0,1]^d$ for $Y^{(m)}$ (it is not necessary to assume that $my \in \Z^d$). Similarly, we choose an
arbitrary starting point $i \in E$ for $\xi$. To establish the LDP, we have to check Assumptions A1 and A3 in \cite{dupuis:88}.
In our framework, A1 is
clearly satisfied. We investigate A3. We first prove that $H$ is $x$-Lipschitz continuous, {\it uniformly} 
in $\alpha$ (so that $L$ is also $x$-Lipschitz continuous, uniformly in
$v$). For $x \in \R^d$, $\alpha \in \R^d$ and $j \in
E$, we set
\begin{equation*}
H_0(x,\alpha;j) = \ln \bigl( \sum_{v \in {\mathcal V}} \exp(\langle \alpha,v \rangle) p(x,j,v) \bigr).
\end{equation*}
By {\bf (A.1)} and {\bf (A.2)}, $H_0$ is $x$-Lipschitz continuous (uniformly in $\alpha$ and $j$). The Lipschitz constant is denoted by $K'$. By \eqref{def:H},
 we obtain for $x, y$ in $\R^d$ and $\alpha \in \R^d$
\begin{equation*}
\begin{split}
&H(x,\alpha) - H(y,\alpha)
\\
&= \lim_{n \rightarrow + \infty} \frac{1}{n} \bigl\{ \ln {\mathbb E}^j \bigl[ \exp \bigl( \sum_{k=1}^n H_0(x,\alpha;\xi_k) \bigr) \bigr] -
\ln {\mathbb E}^j \bigl[ \exp \bigl( \sum_{k=1}^n H_0(y,\alpha;\xi_k) \bigr) \bigr] \bigr\}
\leq K' |x-y|.
\end{split}
\end{equation*}
It remains to estimate the conditional law of the increments of $Y^{(m)}$ given the past. For a given $t>0$,
we consider a 1-Lipschitz continuous
function $\phi \in {\mathcal C}([0,t];\R^d)$, with $\phi_0 \in [0,1]^d$.
From Subsection \ref{subsection:skorohod}, we know that $\Psi$ is $1/2$-H\"older continuous on compact subsets of ${\mathcal C}_{[0,1]^d}([0,t];\R^d)$, so that we can find a
constant $\gamma >0$ such that $\rho_{0,t}(\Psi(Y^{(m)}),\Psi(\phi)) \leq \gamma \rho_{0,t}^{1/2}(Y^{(m)},\phi)$.
For $\alpha \in \R^d$, $\delta, \Delta >0$ and $A \in {\mathcal G}^m_t$,
with ${\mathbb
P}(A) \not= 0$ and $A \subset \{\rho_{0,t}(Y^{(m)},\phi) \leq \delta\}$ (so that $A \subset \{\rho_{0,t}(\Psi(Y^{(m)}),\Psi(\phi)) \leq \gamma \delta^{1/2} \}$), we have
\begin{equation*}
\begin{split}
&{\mathbb E}^i \bigl[ \exp \bigl( m \langle \alpha, Y^{(m)}_{t+\Delta} - Y^{(m)}_t \rangle \bigr) |A \bigr]
\\
&\leq e^{2|\alpha|} {\mathbb E}^i \biggl[ \exp \bigl( \sum_{k=\lfloor tm \rfloor +1}^{\lfloor (t+ \Delta)m \rfloor } \langle \alpha, f(\bar{X}_{k/m}^{(m)},\xi_k,U_k) \rangle \bigr) |A \biggr]
\\
&= e^{2|\alpha|} {\mathbb E}^i \biggl[ \exp \bigl( \sum_{k=\lfloor tm \rfloor +1}^{\lfloor (t+\Delta)m \rfloor -1} \langle \alpha,
f(\bar{X}_{k/m}^{(m)},\xi_k,U_k) \rangle \bigr)
\exp \bigl( H_0(
\bar{X}_{\lfloor (t+\Delta)m \rfloor/m}^{(m)},\alpha;\xi_{\lfloor (t+\Delta) m \rfloor}) \bigr) |A \biggr]
\\
&\leq e^{2|\alpha|} \exp( K' \Delta) {\mathbb E}^i \biggl[ \exp \bigl( \sum_{k=\lfloor tm \rfloor +1}^{\lfloor (t+\Delta)m \rfloor - 1} \langle \alpha, f(\bar{X}_t^{(m)},\xi_k,U_k) \rangle \bigr)
\exp \bigl( H_0(\bar{X}_t^{(m)},\alpha;\xi_{\lfloor (t+\Delta) m \rfloor}) \bigr) |A \biggr]
\end{split}
\end{equation*}
By iterating the procedure, we obtain
\begin{equation*}
\begin{split}
&{\mathbb E}^i \bigl[ \exp \bigl( m \langle \alpha, Y^{(m)}_{t+\Delta} - Y^{(m)}_t \rangle \bigr) |A \bigr]
\\
&\leq e^{2|\alpha|} \exp(K' \Delta^2 m) {\mathbb E}^i \biggl[ \exp \bigl( \sum_{k=\lfloor tm \rfloor +1}^{\lfloor (t+\Delta)m \rfloor}  H_0(
\bar{X}_t^{(m)},\alpha;\xi_{k}) \bigr) |A \biggr]
\\
&\leq e^{2|\alpha|} \exp(K' \Delta (\Delta^2 m+1))
{\mathbb E}^i \biggl[ \exp \bigl( \sum_{k=\lfloor tm \rfloor +1}^{\lfloor (t+\Delta)m \rfloor}  H_0(Z_t^{(m)},\alpha;\xi_{k}) \bigr) |A \biggr]
\\
&\leq e^{2|\alpha|} \exp( K' \Delta (\Delta + \gamma \delta^{1/2}+1/m) m) {\mathbb E}^i \biggl[ \exp \bigl( \sum_{k=\lfloor tm \rfloor +1}^{\lfloor (t+\Delta)m \rfloor}  H_0(
\Psi(\phi)_t,\alpha;\xi_{k}) \bigr) |A \biggr]
\\
&\leq e^{4|\alpha|} \exp(K' \Delta (\Delta + \gamma \delta^{1/2}+1/m) m) \sup_{j \in E} {\mathbb E}^j  \biggl[ \exp \bigl( \sum_{k=0}^{\lfloor  \Delta m \rfloor  } \langle \alpha,
f(\Psi(\phi)_t,\xi_k,U_k) \rangle \bigr) \biggr].
\end{split}
\end{equation*}
We deduce that, uniformly in the starting points $y$ and $i$, uniformly in $\alpha$ on compact subsets and uniformly in $(A,\phi)$ satisfying $A \subset \{\rho_{0,t}(Y^{(m)},\phi) \leq \delta\}$
\begin{equation*}
\limsup_{m \rightarrow + \infty} 1/(m\Delta)
\ln {\mathbb E}^i \bigl[ \exp \bigl( m \langle \alpha, Y^{(m)}_{t+\Delta} - Y^{(m)}_t \rangle \bigr) |A \bigr] \leq H(\Psi(\phi)_t,\alpha)
+ 2 K' (\Delta + \gamma \delta^{1/2}).
\end{equation*}
Similarly, we can prove a lower bound for the liminf. Even if written in a different manner (because of the Skorohod mapping $\Psi$
and because of the conditioning -- we give a precise sense to the right-hand side in \cite[A3, (3.6), (3.7)]{dupuis:88} -- ),
these two bounds correspond to those required in Assumption
A3 in \cite{dupuis:88} (see the discussion on this point in \cite[Section 5]{dupuis:88}).

We deduce that the sequence $(Y^{(m)})_{m \geq 1}$ satisfies on ${\mathcal C}_y([0,T];\R^d)$ ($T>0$) a LDP with the
normalizing factor $m$ and
with the action functional
$\I0T^y : \phi \mapsto \int_0^T
L(\Psi(\phi)_t,\dot{\phi}_t)dt$ if $\phi_0=y$ and $\phi$ is absolutely continuous and $\infty$ otherwise. We let the reader check that this
action functional is lower semicontinuous on ${\mathcal C}_{[0,1]^d}([0,T];\R^d)$ and that its level sets are compact for the supremum norm topology.
By the ``robust'' version of the G\"artner-Ellis proved in \cite{dupuis:88}, the LDP is uniform in $y \in [0,1]^d$.

The uniformity of the LDP with respect to the initial condition is crucial. By the regularity of $L$ in $x$ (it is Lipschitz continuous, uniformly in $\alpha$), it is plain to
deduce that for any $x \in [0,1]^d$, for any closed subset $F \in {\mathcal C}_{[0,1]^d}([0,T];\R^d)$ and any open subset $G \in {\mathcal C}_{[0,1]^d}([0,T];\R^d)$
\begin{equation}
\label{LDP:unif}
\begin{split}
& \lim_{\delta \searrow 0} \sup_{|y-x|< \delta}\limsup_{m \rightarrow + \infty} m^{-1} \ln {\mathbb P} \bigl\{Y^{(m),y} \in F \bigr\} \leq - \inf_{\phi \in F} \I0T^x(\phi),
\\
&\lim_{\delta \searrow 0} \inf_{|y-x|< \delta}\liminf_{m \rightarrow + \infty}
 m^{-1} \ln {\mathbb P} \bigl\{Y^{(m),y} \in G \bigr\} \geq - \inf_{\phi \in G} \I0T^x(\phi),
\end{split}
\end{equation}
where the notation $Y^{(m),y}$ indicates that $Y^{(m)}$ starts from $y$ (i.e. $Y^{(m),y}_0=y$).

By the contraction principle (see e.g. \cite[Theorem 4.2.1, p. 126]{DZ}),  for any $y \in [0,1]^d$, $(\Psi(Y^{(m)}))_{m \geq 1}$ satisfies on ${\mathcal C}_y([0,T];\R^d)$ a LDP with $m$ as
normalizing factor and with the following action functional
\begin{equation}
\label{eq:cost}
\J0T^y (\phi) = \inf \biggl\{ \int_0^T L(\Psi(\psi)_s,\dot \psi_s)ds, \, \Psi(\psi) = \phi \biggr\} =
\inf \biggl\{ \int_0^T L(\phi_s,\dot \psi_s)ds, \, \Psi(\psi) = \phi
\biggr\},
\end{equation}
if $\phi_0=y$ and there is an absolutely continuous path $\psi$ such that $\Psi(\psi) = \phi$, and $\J0T^y(\phi) = \8$ otherwise.

Let us mention at this point that an alternative, more explicit expression of $\J0T^y$
will be given below. Again, the action functional $\J0T^y$ is lower semicontinuous on the
set ${\cal C}_{[0,1]^d}([0,T];$ $[0,1]^d)$. The proof is rather standard and is left to the reader. We can also prove that the level sets ${\cal J}_{y,T}(a) =
\{ \phi \in  {\cal C}_y([0,T];\R^d): \, \J0T^y(\phi) \leq
a\}$, for $y \in [0,1]^d$, are compact in the supremum norm topology. Moreover, \eqref{LDP:unif} yields for any $x \in [0,1]^d$
\begin{equation}
\label{LDPref:unif}
\begin{split}
&
\lim_{\delta \searrow 0} \sup_{|y-x|< \delta}\limsup_{m \rightarrow + \infty}
m^{-1} \ln {\mathbb P} \bigl\{\Psi(Y^{(m),y}) \in F \bigr\} \leq - \inf_{\phi \in F} \J0T^x(\phi),
\\
&
\lim_{\delta \searrow 0} \inf_{|y-x|< \delta}\liminf_{m \rightarrow + \infty}
m^{-1} \ln {\mathbb P} \bigl\{\Psi(Y^{(m),y}) \in G \bigr\} \geq - \inf_{\phi \in G} \J0T^x\phi),
\end{split}
\end{equation}

Now, we can come back to the sequence $(\bar{X}^{(m)})_{m \geq 1}$. For a sequence $(x_m)_{m \geq 1}$ of initial conditions in $[0,1]^d$,
with $mx_m \in \Z^d$ and $x_m \rightarrow x$, we have 
$|\bar{X}^{(m),x_m}_t - \Psi(Y^{(m),x_m}_t)|_{\infty} \leq 1/m$ for all $t$. By \eqref{LDPref:unif}, we deduce
\begin{thm} \label{th:LDPreflect}
Assume that {\bf (A.1--3)} are in force and
consider $T>0$, $x \in [0,1]^d$ and a sequence $(x_m)_{m \geq 1}$  converging
towards $x$, with $m x_m \in [0,m]^d \cap {\mathbb Z}^d$ for all $m
\geq 1$. Then, the sequence $(\bar{X}^{(m)})$ satisfies on ${\mathcal C}([0,T];[0,1]^d)$ a LDP with $m$ as normalizing factor and $\J0T^x$ as action functional.
\end{thm}

Following the proof of \cite[Corollary 5.6.15]{DZ}, we deduce from \eqref{LDPref:unif} the following ``robust'' version (the word ``robust'' indicates that the bounds are uniform
with respect to the initial condition)

\begin{prop}
\label{LDP:varadhan}
Assume that {\bf (A.1--3)} are in force and
consider $T>0$ and $K$ a compact subset of $[0,1]^d$. Then, for any closed subset $F$ of $C([0,T],[0,1]^d)$ and any open subset $G$ of $C([0,T],[0,1]^d)$,
\begin{equation*}
\begin{split}
&\limsup_{m \rightarrow + \infty} \bigl[
m^{-1} \ln \sup_{x \in K, mx \in {\mathbb Z}^d} \P_x
\{\bar{X}^{(m)} \in F\} \bigr]
\leq - \inf_{x \in K} \inf_{\phi \in F}
 J_{0,T}^x(\phi),
\\
&\liminf_{m \rightarrow + \infty}
\bigl[ m^{-1} \ln \inf_{x \in K, mx \in {\mathbb Z}^d} \P_x
\{ \bar{X}^{(m)} \in G\} \bigr]
\geq
- \sup_{x \in K} \inf_{\phi \in G} J_{0,T}^x(\phi).
\end{split}
\end{equation*}
\end{prop}

\subsection{Law of Large Numbers for the Reflected Walk}

We discuss now the zeros of the action functional. We first consider the solution $(\chi_t^{x_0})_{t \geq 0}$, $x_0 \in [0,1]^d$, to \eqref{eq_lim_ref}.
Setting
\begin{equation*}
\forall t \geq 0, \ y_t =  x_0+\int_0^t \bar{f}(\chi_s^{x_0}) ds,
\end{equation*}
we have, for $T>0$, $(\chi_t^{x_0})_{0 \leq t \leq T} = \Psi((y_t)_{0 \leq t \leq T})$. Since the path $t \in \R_+  \mapsto
y_t$ is absolutely continuous, we deduce
$$\J0T^{x_0}(\chi^{x_0}) \leq \int_0^T L \bigl(y_s,\bar{f}(y_s) \bigr) ds = 0,$$
so that $\chi^{x_0}$ is a zero of $\J0T^{x_0}$.
In fact, this is the only possible zero for the given initial condition $x_0$.
Consider indeed another path $\phi$ with values in $[0,1]^d$, such that $\J0T^{x_0}(\phi)=0$. The set of absolutely continuous functions $\psi$ such that
$\psi_0=x_0$,
$$\int_0^T L(\phi_s,\dot{\psi}_s) ds \leq 1 \quad {\rm and} \quad \Psi(\psi) = \phi,$$
is compact. Since the functional $\psi \mapsto \int_0^T L(\phi_s,\dot{\psi}_s) ds$ is lower semicontinuous, it attains its infimum on this compact set. Hence,
there exists an absolutely continuous function $\psi$ such that $\psi_0=x_0$ and
$$\int_0^T L(\phi_s,\dot{\psi}_s) ds = 0 \quad {\rm and} \quad \Psi(\psi) = \phi.$$
It is clear that $\dot{\psi}_t = \bar{f}(\phi_t)$. Since $\Psi(\psi) = \phi$, there exists a process $k$ as in \eqref{eq_lim_ref} such that
\begin{equation*}
\forall t \in [0,T] , \ \phi_t = x_0 + \int_0^t \bar{f}(\phi_s)ds - k_t.
\end{equation*}
This proves that $\phi=\chi^{x_0}$ up to time $T$.

A direct consequence is the following
\begin{cor}
\label{conv_pb}
Assume that {\bf (A.1--3)} are in force and consider
a sequence $(x_m)_{m \geq 1}$ in $[0,1]^d$, with $m
  x_m \in {\mathbb Z}^d$ for all $m \geq 1$, such that $x_m \rightarrow x$ as
  $m \rightarrow +\infty$.
Then,
the sequence of random paths
$(\bar{X}^{(m)})_{m \geq 1}$, with $\bar{X}^{(m)}_0 =x$ for all $m
  \geq 1$, converges,  in
probability, uniformly on compact time intervals
to the solution  $(\chi^{x}_t)_{t \geq 0}$
of the (averaged) reflected differential equation \eqref{eq_lim_ref},
  with $\chi^x_0=x$.
\end{cor}

\subsection{A Different Expression for the Action Functional}
Following \cite{dupuis 87}, we write the action functional $\J0T$ in a different way.
We recall that
 ${\mathcal N}(x)$ denotes  the set of unit outward normals to
$\partial [0,1]^d$ at a point $x$ on the boundary. We
define the function $\Lr$ by $\Lr(x,\cdot)=L(x,\cdot)$  for
$x \in (0,1)^d$, and for $x \in \partial [0,1]^d$,
\begin{equation}
  \label{eq:Lref}
\Lr (x,v)= \left\{
  \begin{array}{ll}
   +\8 & \quad{\rm if \ } \exists n \in {\mathcal N}(x) : \langle v,n \rangle >0\\
   L(x,v) &\quad {\rm if \ } \forall n \in {\mathcal N}(x) : \langle v, n \rangle <0\\
    {\mathop{\inf}\limits_{\beta \geq 0, n \in {\mathcal N}(x), n \perp v}}
       L(x,v+\beta n)
        &\quad {\rm otherwise}
   \end{array} \right.
\end{equation}
The last case  occurs when $ \langle v, n \rangle \leq0$
$ \forall n \in {\mathcal N}(x)$ and
$\exists n' \in {\mathcal N}(x): \langle v, n' \rangle=0$. Then,
the motion takes place on the boundary, in the sense that,
for $\epsilon>0$ small enough,
$x + \epsilon v$ remains in the face orthogonal to $n'$.
Observe that, in contrast to $L(x,\cdot)$, the function $\Lr(x,\cdot)$
may be non convex and
discontinuous for $x \in \partial [0,1]^d$.
\begin{thm}
\label{expr_act_func}
Assume that {\bf (A.1--3)} are in force.
If $\phi$ is absolutely continuous it holds
\begin{equation*}
\J0T^{\phi_0}(\phi) = \int_0^T \Lr (\phi_t,\dot{\phi}_t ) dt.
\end{equation*}
If $\phi$ is not absolutely continuous, then $\J0T^{\phi_0}(\phi)= \infty$.
\end{thm}
$\Box$
By Theorem 2.2 in \cite{lions sznitman 1984}, we know that $\Psi(\psi)$ is absolutely continuous if $\psi$ is absolutely continuous.
In particular, if $\phi$ is not absolutely continuous, there cannot exist an absolutely continuous $\psi$ such that $\Psi(\psi)=\phi$.

Assume now that $\phi$ is absolutely continuous. Then, there exists at least one absolutely continuous path $\psi$ such that $\Psi(\psi)=\phi$, namely $\phi$
itself with $k=0$. We thus denote by $\psi$ an absolutely continuous path such that $\phi = \Psi(\psi)$ and set $k = \psi - \phi$. Then
$k$ is also absolutely continuous and $\dot{k}_t = \beta_t n_t$ with $\beta_t = d |k|_t/dt \geq 0$ ($=0$ if $\phi_t \not \in \partial [0,1]^d$)
and $n_t \in {\mathcal N}(\phi_t)$ if $\phi_t \in \partial [0,1]^d$.
Moreover, for a.e. $t$, for all $\ell \in \{1,\dots,d\}$, $(\dot{\phi}_t)_{\ell} {\mathbf 1}_{\{(\phi_t)_{\ell} \in \{0,1\}\}} = 0$ so that $\dot{\phi}_t
\perp \dot k_t$. Hence
\begin{equation*}
\int_0^T L(\phi_t,\dot{\psi}_t)dt \geq \int_0^T
\Lr (\phi_t,\dot{\phi}_t) dt.
\end{equation*}
This proves that
\begin{equation*}
\J0T(\varphi) \geq
\int_0^T
\Lr (\phi_t,\dot{\phi}_t) dt.
\end{equation*}
We investigate the converse inequality. If the right-hand side is infinite, the proof is over. Thus, we can assume that it is finite,
in particular $\Lr (\phi_t,\dot{\phi}_t)<\8$  for almost every $t \in [0,T]$.
It is enough to construct some $\psi$ with  $\Psi(\psi) = \phi$
and $ L(\phi_t,\dot{\psi}_t)=   \Lr (\phi_t,\dot{\phi}_t)$ a.e..
For times $t$'s when $\phi_t \in  \partial [0,1]^d$ and
$\Lr (\phi_t,\dot{\phi}_t)<\8$ is given by the last line of
(\ref{eq:Lref}),  the infimum is achieved at some pair
 $\beta_t \geq 0, n_t \in {\mathcal N}(\phi_t)$ (this pair is unique by the strict convexity of $L$).
Since $H(x,\alpha)$ is bounded by $|\alpha|$, $|v| >1 \Rightarrow L(x,v)=+\infty$. We deduce that $|\dot{\phi}_t + \beta_t n_t|  \leq 1$, so that $|\beta_t| \leq 1 +
|\dot{\phi}_t|$.
For other times $t$, set $\beta_t = 0, n_t$ arbitrary.
The mapping $t \in [0,T]  \mapsto \beta_t$ is clearly measurable and  integrable.
Hence, we can define  $\dot{\psi}_t = \dot{\phi}_t +
\beta_t n_t$,
$\psi_0 = \phi_0$ and $\psi_t=\psi_0+ \int_0^t \dot{\psi}_s ds$.
The function $\psi$ meets all our requirements.
\qed

%
\section{Analysis of the Deadlock Phenomenon} \label{sec-deadlock}
We now investigate the deadlock time of the algorithm. Fixing a real number $\ell \in (0,d)$, we define
\begin{equation} \label{ensG}
G = \{x \in [0,1]^d: \, |x|_1 < \ell\}\;,\quad {\rm and} \quad
\partial G = \{ x \in [0,1]^d: \, |x|_1  = \ell \}
\end{equation}
its boundary relative to
$[0,1]^d$. We also define the discrete counterparts at scale  $m$,
$G^{(m)} = \{x \in (m^{-1} {\mathbb Z}^d) \cap [0,1]^d: \, |x|_1 < m^{-1}
\lfloor m \ell \rfloor \}$, $\bar{G}^{(m)} = \{x \in (m^{-1} {\mathbb Z}^d) \cap [0,1]^d: \, |x|_1 \leq m^{-1}
\lfloor m \ell \rfloor \}$ and $\partial G^{(m)} =
\bar{G}^{(m)} \setminus G^{(m)}
= \{x \in (m^{-1} {\mathbb Z}^d) \cap [0,1]^d: \, |x|_1 = m^{-1}
\lfloor m \ell \rfloor \}$.
 The deadlock time for the process is
$$\tau^{(m)}= \frac{1}{m}  \inf \Big\{n \geq 0: \, |{X}_n|_1= \lfloor m \ell
\rfloor\big\}=\inf \big\{t \geq 0: \, \bar{X}^{(m)}_t \in \partial G^{(m)}\Big\}.
$$

\vspace{5pt}
We consider the following simple situation:
\vspace{5pt}

\noindent
{\bf Assumption (A.4).}
The point 0 is the unique equilibrium point of the RDE (\ref{eq_lim_ref}). It is stable
and attracts the closure $\bar G = G
\cup \partial G$, that is, for all $x_0 \in \bar{G}$ and $t>0$, $\chi_t^{x_0} \in G$ and $\lim_{t \to \8} \chi_t^{x_0} =0$.

Example (\ref{eq:maier}) given below satisfies the previous
assumption provided that $g_1, g_2$ are (strictly) positive on $(0,1]$.

%
%

\medskip

\noindent
 {\bf Quasi-potential.}
The function
$$V(x,y)= \inf \{ \J0T^x(\phi); \, \phi_0=x, \, \phi_T=y, \, T>0\}$$
is called the { quasi-potential}. It describes the cost
for the random path $\bar{X}^{(m)}$ starting from $x$
to reach the point $y \in G$ at some time scaling with $m$ as $m$
becomes large. (We emphasize that, here and below, the notation
$\J0T^x(\phi)$ implicitly assumes that $\phi$ is a function from $[0,T]$ to
$[0,1]^d$.)

\begin{prop}
\label{reg:quasipot}
Under Assumptions {\bf (A.1--3)}, there exists a constant $C>0$, such that, for all $x,y \in [0,1]^d$,
with $\lambda = |x-y|_1 >0$, the
function $\psi : t \in [0,\lambda] \mapsto x + t (y-x)/\lambda$ satisfies
$\psi_0=x$, $\psi_{\lambda}=y$ and $J_{0,\lambda}^x(\psi) \leq C
\lambda$. In particular, $V(x,y) \leq C |x-y|_1$.
\end{prop}

\noindent
$\Box$ Proof. By {\bf (A.2)} and \eqref{def:H}, for all $x \in [0,1]^d$ and $\alpha \in {\mathbb R}^d$,
$
H(x,\alpha) \geq \ln ( c \exp(|\alpha|_{\infty}) ) = \ln(c) + |\alpha|_{\infty}$,
with $c = \inf\{p(z,i,v); \, z \in [0,1]^d, \, i \in E, \, v \in {\mathcal V}\} >0$. Hence, for all $v \in \R^d$,
$L(x,v) \leq \sup_{\alpha} \{\langle \alpha,v \rangle - |\alpha|_{\infty}\} - \ln(c) \leq
\sup_{\alpha} \{ |\alpha|_{\infty}(|v|_1-1)\} - \ln(c) \leq - \ln(c)$ if $|v|_1 \leq 1$. The proof is easily completed. \qed

 \subsection{ Deadlock Time and Exit Points} \label{sec:deadlocktime-points}

We define the minimum value of  the quasi-potential
$V(0,\cdot)$ on the boundary  of $G$ by
$$\bar{V} = \inf \{ \J0T^0(\phi); \, \phi_0=0, \, |\phi_T|_1=\ell, \, T>0\}$$
and the  set of minimizers
\begin{equation}
  \label{eq:setmini}
  \cM=\{y \in \partial G : \, V(0,y)=\bar{V}\}.
\end{equation}
By Proposition \ref{reg:quasipot}, $\bar{V}$ is finite. A consequence of Theorem \ref{th:LDPreflect} and Proposition
\ref{reg:quasipot} is
\begin{thm} \label{th:tmoy}
Assume that {\bf (A.1--4)} are in force and
consider a sequence $(x_m)_{m \geq 1}$ in $G$, with $mx_m \in {\mathbb
  Z}^d$ for all $m \geq 1$, such that
$x_m \rightarrow x \in G$. Then,
\begin{equation}
  \label{eq:tmoy1}
\bbE_{x_m}[ \tau^{(m)}] = \exp \bigl[  m (\bar{V}+o(1)) \bigr]
\end{equation}
as $m \to \8$. Moreover, for all positive $\delta$,
\begin{equation}
  \label{eq:tmoy2}
  \lim_{m \to \8} \bbP_{x_m} \big\{ \exp  [m (\bar{V} - \delta)] <
\tau^{(m)} < \exp  [m (\bar{V} + \delta)] \big\} =1.
\end{equation}
 Finally, for all $\varepsilon>0$, it holds
\begin{equation}
  \label{eq:tmoy3}
  \bbP_{x_m} \big\{ d(\bar{X}^{(m)}_{\tau^{(m)}},{\mathcal M})<\varepsilon \big\} \to 1 \quad
{\rm as\ } m \to \8 \; ,
\end{equation}
where $d(\bar{X}^{(m)}_{\tau^{(m)}},{\mathcal M})$ denotes the distance from $\bar{X}^{(m)}_{\tau^{(m)}}$ to the set ${\mathcal M}$.
\end{thm}
$\Box$ Proof.
The proof follows the standard theory of Markov perturbations of
dynamical systems in \cite[Chapter 6]{FV}. For the sake of
completeness, we provide the main steps according to the very detailed
scheme in \cite[Section 5.7]{DZ} (Section 5.7 is devoted to large
deviations for stochastic differential equations with a small noise).

We define, for $x \in G$, $V(x,\partial G) = \inf \{\J0T^x(\phi); \,
\phi_0=x, \, |\phi_T|_1=\ell, \, T>0\}$, so that $V(0,\partial G)=\bar V$.
We also define the ball $\bar{B}_{\rho}^{(m),+}$
in the lattice orthant of mesh $1/m$,
$$\bar{B}_{\rho}^{(m),+} = \{ z \in (m^{-1} {\mathbb Z}^d) \cap [0,1]^d: \,
|z|_1 \leq m^{-1} \lfloor m \rho \rfloor\}. $$

In the whole proof, we assume that $0 < 2\rho < \ell$.

\begin{lem}
\label{lem:5.7.18}
For
any $\eta>0$ and for any $\rho>0$ small enough, there exists $T_0
< + \infty$, such that
\begin{equation*}
\liminf_{m \rightarrow +\infty} m^{-1} \ln \inf_{x \in
  \bar{B}_{\rho}^{(m),+}}
\P_x \{\tau^{(m)} \leq T_0\} \geq - \bar{V} -\eta.
\end{equation*}
\end{lem}

$\Box$ Proof.
We first fix a small $\eta >0$. By the definition of
$\bar{V}$, we can find $S_0>0$ and $\phi^0 \in {\mathcal
  C}([0,S_0],[0,1]^d)$, with $\phi_0^0=0$, such that $J_{0,S_0}^0(\phi^0) \leq
\bar{V} + \eta$ and $\phi_{S_0} \in \partial G$. By Proposition \ref{reg:quasipot} and by the additive form of $J$,
see Theorem \ref{expr_act_func}, we can extend $\phi$ after $S_0$ to leave
$\bar G$ at low cost, and assume that
$\phi([0,S_0])\cap \partial G \not = \emptyset$ and $\delta = d(\phi_{S_0},\bar{G})>0$.

For $x \in [0,1]^d$, $|x|_1 < 2\rho$, we can find by Proposition \ref{reg:quasipot}
a path $\zeta \in {\mathcal C}([0,2\rho];[0,1]^d)$ such that $\zeta_0 = x$,
$\zeta_{2\rho}=0$ and $J_{0,2\rho}^{x}(\zeta)
\leq C \rho$. By concatenating $\zeta$ and $\phi$, we obtain a path $\phi^x$.
For $\rho \leq \eta/C$, it satisfies
\begin{equation*}
J_{0,T_0}(\phi^x) \leq \bar{V} + 2 \eta \; ,
\end{equation*}
with $T_0 = S_0 + 2\rho$. Now, the set
\begin{equation*}
\Psi = \bigcup_{x \in [0,1]^d, \ |x|_1 < 2\rho} \bigl\{ \psi \in {\mathcal C}([0,T_0];[0,1]^d) : \, \rho_{0T_0}(\psi,\phi^x) < \delta /2 \bigr\},
\end{equation*}
is an open subset of ${\mathcal C}([0,T_0];[0,1]^d)$. By Proposition \ref{LDP:varadhan},
\begin{equation*}
\begin{split}
\liminf_{m \rightarrow + \infty} m^{-1} \inf_{x \in \bar{B}_{\rho}^{(m),+}} \ln \P_x
\{\tau^{(m)} \leq T_0\}
&\geq \liminf_{m \rightarrow + \infty} m^{-1} \inf_{x \in \bar{B}_{\rho}^{(m),+}} \ln \P_x
\{\bar{X}^{(m)} \in \Psi \}
\\
&\geq - \sup_{x \in [0,1]^d \, , \ |x|_1 \leq \rho} \inf_{\psi \in \Psi} J_{0,T_0}^x(\psi)
\\
&\geq - \sup_{x \in [0,1]^d \, , \ |x|_1 \leq \rho}  J_{0,T_0}^x(\phi^x)
\geq - \bar{V} - 2 \eta.
\end{split}
\end{equation*}
This completes the proof. \qed
\begin{lem}
\label{lem:5.7.19}
Let $\sigma_{\rho} = \inf\{t \geq 0: \,
\bar{X}_t^{(m)} \in \bar{B}_{\rho}^{(m),+} \cup \partial
G^{(m)}\}$. Then,
\begin{equation*}
\lim_{t \rightarrow + \infty} \limsup_{m \rightarrow + \infty}
\bigl[ m^{-1} \ln \sup_{x \in G^{(m)}}
\P_x\{\sigma_{\rho}>t\} \bigr] = - \infty.
\end{equation*}
\end{lem}
\noindent
$\Box$ Proof. For $x \in  \bar{B}_{\rho}^{(m),+}$, there is nothing to prove. Now,
as in the proof of \cite[Lemma 5.7.19]{DZ}, we can define for $t \geq 0$ the closed set $\Psi_t =
 \{\phi \in {\mathcal C}([0,t];[0,1]^d): \, \psi_s \in \bar{G}
   \setminus B_{\rho/2}^+, \  \forall s \in [0,t] \}$, where $B_{\rho/2}^+$
 is the ball in the orthant, $B_{\rho/2}^+=\{z \in [0,1]^d: \, |z|_1 < \rho/2\}$. For $\bar{X}_0^{(m)} \in \bar{G}^{(m)}$ and $m$ large,
 $\sigma_{\rho} > t$ implies  $(\bar{X}^{(m)}_s)_{0 \leq s \leq t} \in \Psi_t$. By Proposition \ref{LDP:varadhan},
\begin{equation*}
\begin{split}
&\limsup_{m \rightarrow + \infty} \bigl[ m^{-1} \ln \sup_{x \in G^{(m)} \setminus \bar{B}^{(m),+}_{\rho}}
\P_x \{ \sigma_{\rho}>t \} \bigr]
\\
&\leq \limsup_{m \rightarrow + \infty} \bigl[ m^{-1} \ln \sup_{x \in G^{(m)}
  \setminus \bar{B}^{(m),+}_{\rho} } \P_{x} \{ \bar{X}^{(m)} \in \Psi_t  \} \bigr] \leq - \inf_{x \in \bar{G} \setminus {B}_{\rho/2}^+} \inf_{\psi
  \in \Psi_t} J_{0,t}^x(\psi) = - \inf_{\psi \in \Psi_t} J_{0,t}^{\psi_0}(\psi).
\end{split}
\end{equation*}
Using the stability of the solutions to \eqref{eq_lim_ref} (see Subsection \ref{subsection:skorohod}) and the
additivity of the action functional (see Theorem
\ref{expr_act_func}), we can complete as in \cite{DZ}. \qed

\begin{lem}
\label{lem:5.7.21}
Let $N$ be a closed subset, included in $\partial G$. Then, for every
$\varepsilon >0$,
\begin{equation*}
\lim_{\rho \rightarrow 0} \limsup_{m \rightarrow + \infty}
\bigl[
m^{-1} \ln \sup_{y \in S_{2 \rho}^{(m),+}} \P_y \{{\rm
    dist}(\bar{X}^{(m)}_{\sigma_{\rho}},N) < \varepsilon \} \bigr]\leq
  - \inf_{z \in N} V(0,z) + \delta_{\varepsilon},
\end{equation*}
with $\lim_{\varepsilon \rightarrow 0} \delta_{\varepsilon} =0$.
Here, $S_{ \rho}^{(m),+}=\{z
  \in m^{-1} {\mathbb Z}^d \cap [0,1]^d: \, |z|_1 = m^{-1} \lfloor \rho
  m \rfloor\}$ is the sphere in the lattice orthant with mesh $1/m$.
\end{lem}
$\Box$ Proof.
The proof is the same as in \cite[Lemma 5.7.21]{DZ}, except the application of
Corollary 5.6.15. For $T>0$, we can define, as in \cite{DZ}, $\Phi = \{\phi \in
{\mathcal C}([0,T];[0,1]^d): \, \exists t \in [0,T], \, \phi_t \in
N\}$. If $\sigma_{\rho }\leq T$ and ${\rm dist}(\bar{X}^{(m)}_{\sigma_{\rho}},N) < \varepsilon$,
then $\rho_{0,T}(\bar{X}^{(m)},\Phi) \leq \varepsilon$. So that, Proposition \ref{LDP:varadhan} yields
\begin{equation*}
\begin{split}
&\limsup_{m \rightarrow + \infty} \bigl[
m^{-1} \ln \sup_{y \in S_{2\rho}^{(m),+}} \P_y \{
\sigma_{\rho} \leq T, \ {\rm dist}(\bar{X}^{(m)}_{\sigma_{\rho}},N) <
\varepsilon\} \bigr]
\\
&\leq - \inf_{d(y,S_{2 \rho}^+) \leq \varepsilon} \; \inf_{d(\phi,\Phi)\leq
 \varepsilon} J_{0,T}^y(\phi)
= - \inf \bigl\{J_{0,T}^{\phi_0}(\phi); \,
d(\phi_0,S_{2 \rho}^+) \leq \varepsilon, \ d(\phi,\Phi)\leq  \varepsilon \bigr\},
\end{split}
\end{equation*}
with  $S_{ \rho}^{+}=\{z
  \in [0,1]^d: \, |z|_1 =  \rho \}$ is the sphere in the lattice orthant.
Using the semicontinuity of $J$, the reader can check that (see \cite[Lemma 4.1.6]{DZ})
\begin{equation*}
\lim_{\varepsilon \rightarrow 0} \inf_{d(\phi_0,S_{2 \rho}^+) \leq \varepsilon} \; \inf_{d(\phi,\Phi)\leq  \varepsilon} J_{0,T}^{\phi_0}(\phi)
=
\inf_{\phi_0 \in S_{2 \rho}^+} \inf_{\phi \in \Phi} J_{0,T}^{\phi_0}(\phi).
\end{equation*}
The end of the proof is the same. \qed

\begin{lem}
\label{lem:5.7.22}
Let $K$ be a compact subset of $[0,1]^d$ included in $G^{(m)}$ for $m$ large.
Then,
\begin{equation*}
\lim_{m \rightarrow + \infty} \inf_{x_m \in K} \P_{x_m} \{ \bar{X}^{(m)}_{\sigma_{\rho}}
\in \bar{B}_{\rho}^{(m),+} \} = 1.
\end{equation*}
\end{lem}
$\Box$ Proof.
The proof is
the same as in \cite[Lemma 5.7.22]{DZ}, up to the infimum over the compact set $K$. By  {\bf (A.4)} and by
the regularity of the flow $(t,x) \in \R_+ \times [0,1]^d \mapsto \chi_t^x$, the hitting time $T = \inf\{t \geq 0: \, \forall x \in [0,1]^d,
|\chi_t^x|_1 \leq \rho/2\}$ is finite. Moreover, $\inf_{t \in [0,T], x \in K} d(\chi_t^x,\partial G) >0$. Using Corollary \ref{conv_pb}, it is plain to conclude. \qed

Finally, we have the following obvious result
\begin{lem}
\label{lem:5.7.23}
$\displaystyle \sup_{x \in G^{(m)}} \P_x \{ \sup_{0 \leq t \leq \rho}
|\bar{X}_t^{(m)} -x| \geq 2 \rho\}
=0$.
\end{lem}

It now remains to follow the proof of
\cite[Theorem 5.7.11]{DZ}.
The crucial point to note is the following: $\tau^{(m)}$ and
$\sigma_{\rho}$ take their values in $m^{-1} {\mathbb N}$ and are
stopping times for the filtration $({\mathcal F}^{\xi,X}_{\lfloor mt
  \rfloor})_{t \geq 0}$. In particular, the Markov property (for $(\xi,X)$) applies
quite easily. For example, for $x \in G^{(m)}$ and $s,t$ in ${\mathbb N}^*$ (and thus in $m^{-1} {\mathbb N}^*$),
\begin{equation*}
\begin{split}
&\!\!\!\!\!\!\!\!
{\mathbb P}_x \{\sigma_{\rho} \leq t, \ \bar{X}_{\sigma_{\rho}}^{(m)} \in \bar{B}^{(m),+}_{\rho}, \tau^{(m)}
\leq t+s \}
\\
&\geq {\mathbb P}_x \{ \sigma_{\rho} \leq  t, \ \bar{X}_{\sigma_{\rho}}^{(m)} \in \bar{B}^{(m),+}_{\rho}  \}
\inf_{y \in \bar{B}_{\rho}^{(m),+}, i \in E} {\mathbb P}_y^i \{ \tau^{(m)} \leq s \},
\end{split}
\end{equation*}
so that
\begin{equation*}
\begin{split}
&\!\!\!\!\!\!\!\!
{\mathbb P}_x\{\tau^{(m)} \leq t+s\}
\\
&\geq {\mathbb P}_x \{\sigma_{\rho} \leq t, \ \bar{X}_{\sigma_{\rho}}^{(m)} \in \bar{B}^{(m),+}_{\rho}, \tau^{(m)}
\leq t+s \} + {\mathbb P}_x \{\sigma_{\rho} \leq t, \ \bar{X}_{\sigma_{\rho}}^{(m)} \in \partial G^{(m)} \}
\\
&\geq {\mathbb P}_x \{ \sigma_{\rho} \leq t \}
\inf_{y \in \bar{B}_{\rho}^{(m),+}, i \in E} {\mathbb P}_y^i \{ \tau^{(m)} \leq s \}.
\end{split}
\end{equation*}
This shows that (5.7.24) in \cite{DZ} holds. Similarly, for $t \in {\mathbb N}^*$ and
$k \in {\mathbb N}$,
\begin{equation*}
\P_x \{\tau^{(m)} > (k+1) t\}
\leq \P_x \{\tau^{(m)} > k  t \}
\sup_{y \in G^{(m)}, i \in E} \P_y^i \{ \tau^{(m)} > t  \}.
\end{equation*}
Now, the  upper bounds in \eqref{eq:tmoy1} and \eqref{eq:tmoy2}
can be derived as in  \cite{DZ}.

Turn to the lower bounds. Following \cite{DZ}, we introduce the
following notations (pay attention to
 that $m$ in \cite{DZ}
refers to a complete different parameter than in our case):
\begin{equation} \label{eq:tautheta}
\theta_0 = 0, \ \tau_n = \inf \{ t \geq \theta_n:  \,
  \bar{X}_t^{(m)} \in \bar{B}_{\rho}^{(m),+} \cup \partial G^{(m)} \},
  \ \theta_{n+1} = \inf \{t \geq \tau_n: \,  \bar{X}_t^{(m)} \in S_{2 \rho}^{(m),+}\},
\end{equation}

\begin{figure} [htb]
\begin{center}
\includegraphics[
width=0.4\textwidth,angle=0]
{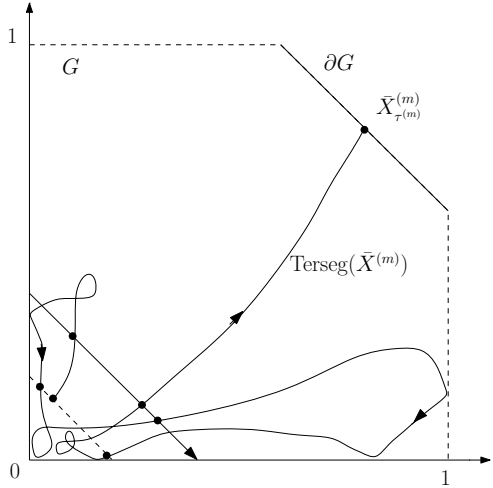}
\caption{The path $\bar X^{(m)}$ up to the deadlock time $\tau^{(m)}$
($d=2, \ell=1.7$).
Spheres $S_\rho, S_{2 \rho}$ are indicated by dashed lines.
The seven large dots on the path are the locations at times
$\theta_0=\tau_0 = 0,  \theta_1, \tau_1,  \theta_2, \tau_2,
 \theta_3, \tau_3= \tau^{(m)}$. The last part of the curve is the terminal 
segment Terseg defined in the proof of Theorem \ref{th:segter}.
}
\label{fig:stopping}
\end{center}
\end{figure}

with $\theta_{n+1}=+\8$ if $\bar{X}_{\tau_n}^{(m)} \in \partial G^{(m)}$.
These stopping times are indicated in Figure \ref{fig:stopping}.
It is plain to obtain (5.7.26) of  \cite{DZ} (with the Markov property and
Lemma \ref{lem:5.7.21}, with $N=\partial G$ and $\varepsilon$ as small
as necessary) as well as (5.7.27) (with Lemma \ref{lem:5.7.23}).
The end of the
proof of the lower bound just follows the strategy in \cite{DZ}.

Turn to the second statement in Theorem \ref{th:tmoy}. This
is a particular case of \textit{b)} in \cite{DZ}.
Set $N = \partial G \cap \{x \in [0,1]^d: \, {\rm dist}(x,{\mathcal M}) \geq
\varepsilon\}$. It
is a closed set. Then, for $\varepsilon'>0$, we can focus on
\begin{equation*}
\sup_{y \in S_{2 \rho}^{(m),+}} \P_y \{ {\rm
  dist}(\bar{X}^{(m)}_{\sigma_{\rho}},N) < \varepsilon'\}.
\end{equation*}
Setting $V_N = \inf_{y \in N} V(0,y)$, we deduce from Lemma
\ref{lem:5.7.21}  that for
$\rho,\varepsilon'>0$ small enough and for $m$ large enough
\begin{equation*}
\sup_{y \in S_{2 \rho}^{(m),+}} \P_y \{ {\rm
  dist}(\bar{X}^{(m)}_{\sigma_{\rho}},N) < \varepsilon'\}
\leq \exp[-m(V_N-\eta)],
\end{equation*}
with $\eta<(V_N-\bar{V})/3 <0$.
Then, we can follow the proof in \cite{DZ} and prove that for $x_m \in G^{(m)}$,
$x_m \rightarrow x \in G$,
\begin{equation*}
\lim_{m \rightarrow + \infty}  \P_{x_m} \{ {\rm
  dist}(\bar{X}^{(m)}_{\tau^{(m)}},N) < \varepsilon'\} = 0.
\end{equation*}
Since ${\rm dist}(\bar{X}^{(m)}_{\tau^{(m)}},\partial G) \leq C/m$, we complete
the proof. \qed

\subsection{Generic Behavior Leading to Deadlock}
\label{sec:deadlock-behavior}

From (\ref{eq:tmoy3}) we observe that when  $\cM$
reduces to a single point $y^*$, the location of the process
$ \bar{X}^{(m)}$ when exiting $G$ converges to $y^*$. We can extend this
observation from the exit point to the path itself before it exits $G$.
To do so, we
first need to  extend the action functional to any
interval of $\R$, which can be done in a trivial way
thanks to
Theorem \ref{expr_act_func}: for any continuous path $(\psi_t)_{t \leq 0}$, with $\lim_{- \infty} \psi =0$, we denote by
$J_{-\infty,0}(\psi)$ the integral of $L^{{\rm ref}}(\psi_t,\dot{\psi}_t)$ from $- \infty$ to $0$.
Since $0$ is a fixed point
for the limit RDE by
Assumption  {\bf (A.4)}, we have $\Lr (0,0)=0$, and then
$$
\inf \{  J_{0,T}^0(\phi); \, \phi_0=0, \, \phi_T=y, \, T>0\}=
\inf \{  J_{-\infty,0}(\psi); \, \lim_{-\infty}\psi=0, \, \psi_0=y\},
$$
for $y \in G$.
Indeed, for all $T, \phi$ as in the left-hand side, the path $\psi$
given by $\psi_t=\phi_{t+T}$ for $t \in [-T,0]$ and $\psi_t=0$
for $t<-T$ is such that $ J_{-\infty,0}(\psi)= J_{0,T}(\phi)$.
This proves that the left-hand side is greater than the right-hand side. Conversely, for a path $\psi$ with $\lim_{-\infty} \psi=0$ and
$\psi_0=y$, we can find, for every $\delta >0$, $T<0$ such that $|\psi_T| < \delta$. By Proposition \ref{reg:quasipot}, we can find a path
$\theta$ from $[0,\delta]$ into $[0,1]^d$, with
$\theta_0=0$ and $\theta_{\delta}=\psi_T$, such that $J_{0,\delta}^0(\theta)\leq C \delta$.
Concatenating this path to
the restriction of the path $\psi$ to $[T,0]$ (up to a trivial change of time in $\psi$), we obtain a new path $\phi$. It is defined on
$[0,T+\delta]$ and satisfies $\phi_0=0$, $\phi_{T+\delta}=y$ and $J_{0,T}^0(\phi) \leq C \delta + J_{-\infty,0}(\psi)$. This proves that the two infimums
are equal.

Now, we can state the convergence result of the exit path.
\begin{thm} \label{th:segter}
Under Assumptions {\bf (A.1--4)},
assume uniqueness of the optimal path to exit $G$ from 0, i.e., assume
that $\cM=\{y^*\}$ and that there is a unique
$\varphi:(-\infty,0] \rightarrow \bar{G}$, $\varphi((-\infty,0)) \subset G$, minimizing
$J_{-\infty,0}(\varphi)$ subject to $\varphi_0 =y^*,
\lim_{t \rightarrow - \infty} \varphi_t = 0$ (in such a case, $\varphi$ is also the unique minimizing path with values in $[0,1]^d$ --
and not only in $\bar{G}$ -- ).
Let $K$ be a compact set, included in $G$, and containing a
neighborhood of the origin. We denote by $\alpha_K^{(m)}$ the last exit time
before $\tau^{(m)}$
of $\bar{X}^{(m)}$ from $K \cap (m^{-1} {\mathbb Z}^d)$. Then, for any
sequence $(x_m)_{m \geq 1}$, $x_m \in G^{(m)}$ and $x_m \rightarrow x \in G$, and any $\varepsilon >0$
\begin{equation*}
\lim_{m \rightarrow + \infty} {\mathbb P}_{x_m} \big\{ \exists t \in
[\alpha_K^{(m)},\tau^{(m)}], \ |\bar{X}_t^{(m)}-
\varphi_{t-\tau^{(m)}}| > \varepsilon \big\} =0.
\end{equation*}
\end{thm}
$\Box$ Proof. Our proof is inspired by \cite[Section 2, Chapter 4]{AZ}. We keep the notations introduced in the proof of
Theorem \ref{th:tmoy}. In addition, we define $\nu = \max \{n \geq
1: \, \theta_{n} < \tau^{(m)} \}$. If $\tau^{(m)} = \tau_0$, we
set
$\nu = 0$. We denote by $\TS(\bar{X}^{(m)})$ the terminal
``segment" of the path $\bar{X}^{(m)}$, that is, the restriction of
$\bar{X}^{(m)}$ to the interval
$[\theta_{\nu},\tau^{(m)}=\tau_{\nu}]$, but shifted in time to the
interval $[0,\tau_{\nu}-\theta_{\nu}]$. More precisely, if we denote
by $\Theta_t$ the shift operator, i.e. $\Theta_t \psi(s) =
\psi(s+t)$, then $\TS(\bar{X}^{(m)})$ is defined as the restriction of
$\Theta_{\theta_{\nu}} (\bar{X}^{(m)})$ to $[0,\tau_{\nu} -
\theta_{\nu}]$.
\medskip

Fix $\varepsilon>0$.
For $y \in \bar{B}_{\rho}^{(m),+}$ and $L \in {\mathbb N}^*$, we have
$\tau_0=0$ and
\begin{equation}
\label{terseg:1}
\begin{split}
&\P_y \{ \rho_{0,\tau_{\nu} -
  \theta_{\nu}}(\TS(\bar{X}^{(m)}),\Theta_{\theta_{\nu}-\tau_{\nu}}
\varphi) \geq  \varepsilon \}
\\
&\leq
 \P_y \{\tau^{(m)} > \tau_{L} \} +
\sum_{k=1}^{L} \P_y \{ \tau^{(m)} = \tau_k,
\rho_{\theta_k,\tau_k}(\bar{X}^{(m)},\Theta_{-\tau_k} \varphi)
\geq \varepsilon \}
\\
&\leq  \P_y \{\tau^{(m)} > \tau_{L} \} + \sum_{k=1}^{L}
\P_y \{ \bar{X}^{(m)}_{\tau_k} \in \partial G^{(m)},
\rho_{\theta_k,\tau_k}(\bar{X}^{(m)},\Theta_{-\tau_k} \varphi) \geq
\varepsilon \}.
\end{split}
\end{equation}
Focus on the second term. The Markov property yields
\begin{equation}
\label{terseg:2}
\begin{split}
&\sum_{k=1}^{L} \P_y \{ \bar{X}^{(m)}_{\tau_k} \in \partial
G^{(m)}, \rho_{\theta_k,\tau_k}(\bar{X}^{(m)},\Theta_{-\tau_k}
\varphi) \geq \varepsilon \}
\\
&\hspace{15pt} \leq L \sup_{z\in S_{2 \rho}^{(m),+}, i \in E} \P_z^i
\{\bar{X}^{(m)}_{\sigma_{\rho}} \in \partial G^{(m)},
\rho_{0,\sigma_{\rho}}(\bar{X}^{(m)},\Theta_{-\sigma_{\rho}}\varphi)
\geq \varepsilon \}.
\end{split}
\end{equation}
For $T > 0$, we can bound the last quantity as follows
\begin{equation}
\label{terseg:3}
\begin{split}
&\P_z \bigl\{ \bar{X}_{\sigma_{\rho}}^{(m)} \in \partial G^{(m)},
\rho_{0,\sigma_{\rho}}(\bar{X}^{(m)},\Theta_{-\sigma_{\rho}}
\varphi) \geq \varepsilon
 \bigr\}
\\
&\leq \P_z \bigl\{ \sigma_{\rho} \geq T \bigr\} + \P_z \bigl\{
\bar{X}_{\sigma_{\rho}}^{(m)} \in \partial G^{(m)}, \sigma_{\rho}
\leq T, \rho_{0,\sigma_{\rho}}(\bar{X}^{(m)},\Theta_{-\sigma_{\rho}}
\varphi) \geq \varepsilon \bigr\}.
\end{split}
\end{equation}
Now set, for $T,r>0$, $\Gamma_T(r) = \{ \psi \in {\mathcal
C}([0,T];[0,1]^d): \, \psi([0,T]) \cap
\partial G \not = \emptyset, \, \rho_{-T,0}(\Theta_T \psi,\varphi) \geq
r \}$. We then recall the following result in \cite{AZ} (see Lemma
2.8, p. 105, the proof relies on the uniqueness of $\varphi$ and is
exactly the same in our setting, except (4), p. 106, which has to be
read $\liminf_{k \rightarrow + \infty}
d_{-T_k,0}(\Theta_{T_k}g^k,\varphi)>0$):
\begin{equation*}
\forall r>0, \ \exists \alpha >0, \ \forall T>0, \ \inf_{\psi \in
\Gamma_T(r), \psi_0=0} J_{0,T}(\psi) > \bar{V} + \alpha.
\end{equation*}
We now consider $T,r>0$ and $\psi \in \Gamma_T(r)$ with $|\psi_0|
\leq 2 \rho$. We then prove that the above lower bound still holds
for $\rho$ small enough. Indeed, we can consider a path
$\tilde{\psi}$, with $\tilde{\psi}_0 = 0$, $\tilde{\psi}_S = \psi_0$
and $\tilde{\psi}_{t+S} = \psi_t$ for $t \in [0,T]$. Using
Proposition \ref{reg:quasipot}, we can assume that $S \leq C \rho$
and that $J_{0,T+S}(\tilde{\psi}) \leq C \rho + J_{0,T}(\psi)$. We
choose $C \rho \leq \alpha/2$. Since
$\rho_{0,T+S}(\tilde{\psi},\Theta_{-(T+S)} \varphi) \geq r$, we have
$J_{0,T+S}(\tilde{\psi}) > \bar{V} + \alpha$. Finally,
$J_{0,T}(\psi)
> \bar{V} + \alpha/2$.

\medskip

We now choose $r=\varepsilon/2$. For the corresponding $\alpha >0$,
we choose $C\rho \leq \alpha/2$ as above. Then, by means of Lemma
\ref{lem:5.7.19}, we can pick $T$ large enough so that for $m$ large enough
\begin{equation}
\label{terseg:4}
\sup_{z \in S_{2\rho}^{(m),+}} \P_z \bigl\{ \sigma_{\rho} \geq T
\bigr\} \leq \exp(-m(\bar{V}+1)).
\end{equation}
Now, for $0 < \varepsilon'< \varepsilon/2$,
\begin{equation*}
\begin{split}
&\sup_{z \in S_{2\rho}^{(m),+}} \P_z
\bigl\{\bar{X}_{\sigma_{\rho}}^{(m)} \in
\partial G^{(m)},    \sigma_{\rho} \leq T,
\rho_{0,\sigma_{\rho}}(\bar{X}^{(m)},\Theta_{-\sigma_{\rho}}
\varphi) \geq \varepsilon \bigr\}
\\
&\leq \sup_{z \in S_{2 \rho}^{(m),+}} \P_z \bigl\{
\rho_{0,T}(\bar{X}^{(m)},A_T(\varepsilon/2,2\rho)) \leq  \varepsilon'
\bigr\},
\end{split}
\end{equation*}
where $A_T(\varepsilon/2,2\rho)$ stands for the set of continuous
functions from $[0,T]$ into $[0,1]^d$, with $|\psi_0| \leq 2 \rho$,
for which we can find $t \in [0,T]$ such that the restriction of
$\psi$ to $[0,t]$ belongs to $\Gamma_t(\varepsilon/2)$. This is a
closed set. Hence, Proposition \ref{LDP:varadhan} yields for $\varepsilon'$ small enough and
$m$ large enough
\begin{equation}
\label{terseg:5}
\begin{split}
&\sup_{z \in S_{2\rho}^{(m),+}} \P_z
\bigl\{\bar{X}_{\sigma_{\rho}}^{(m)} \in
\partial G^{(m)},  \sigma_{\rho} \leq T,
\rho_{0,\sigma_{\rho}}(\bar{X}^{(m)},\Theta_{-\sigma_{\rho}}
\varphi) \geq \varepsilon \bigr\}
\\
&\leq \exp[ - m ( \inf_{z \in B_{2\rho}^+} \inf_{d(\phi,A_T(\varepsilon/2,2\rho)) \leq \varepsilon'}
J_{0,T}^{z}(\phi) - \alpha/12)]
\\
&\leq \exp[- m(\inf_{d(\phi,A_T(\varepsilon/2,2\rho))\leq  \varepsilon'}J_{0,T}^{\phi_0}(\phi)
- \alpha/12)] \leq \exp[- m(\inf_{\phi \in A_T(\varepsilon/2,2\rho)}J_{0,T}^{\phi_0}(\phi) - \alpha/6)],
\end{split}
\end{equation}
the last inequality following from \cite[Lemma 4.1.6]{DZ}.
For all $\phi \in A_T(\varepsilon/2,2\rho)$, there exists $t \in [0,T]$ such that the restriction of $\phi$ to $[0,t]$ belongs to $\Gamma_t(\varepsilon/2)$. We deduce that
$J_{0,T}^{\phi_0}(\phi) \geq J_{0,t}^{\phi_0}(\phi) \geq \bar{V} + \alpha/2$.
Finally, by \eqref{terseg:1}, \eqref{terseg:2}, \eqref{terseg:3}, \eqref{terseg:4} and \eqref{terseg:5},
\begin{equation*}
\begin{split}
\P_y \{ \rho_{0,\tau_{\nu} -
  \theta_{\nu}}(\TS(\bar{X}^{(m)}),\Theta_{\theta_{\nu}-\tau_{\nu}}
\varphi) \geq  \varepsilon \} \leq \P_{y}\{\tau^{(m)} > \tau_{L}\}
+ 2 L \exp(-m(\bar{V}+\alpha/3)).
\end{split}
\end{equation*}
We can conclude as in the proof of \cite[Theorem 5.7.11, (b)]{DZ}.
We can find a constant $C$ such that
\begin{equation*}
\begin{split}
&\sup_{y \in B_{\rho}^{(m),+}} \P_y \{ \rho_{0,\tau_{\nu} -
  \theta_{\nu}}(\TS(\bar{X}^{(m)}),\Theta_{\theta_{\nu}-\tau_{\nu}}
\varphi) \geq  \varepsilon \}
\\
&\hspace{15pt} \leq C L^{-1}
\exp(m(\bar{V}+\alpha/6)) + 2 L \exp(-m(\bar{V}+\alpha/3)).
\end{split}
\end{equation*}
We then choose $L = \lfloor \exp(m(\bar{V}+\alpha/4)) \rfloor$.
For an arbitrary initial condition in $G$, we conclude as in
the proof of \cite[Theorem 5.7.11, (b)]{DZ} by means of Lemma
\ref{lem:5.7.22} (and the Markov property). \qed
\medskip

\subsection{Exponential Limit Law for Deadlock Time}
\label{sec:deadlock-exp}

Since the exponential law is the generic distribution for rare events,
it appears naturally in the following refinement of 
Theorem \ref{th:tmoy} (see e.g. \cite[Theorem 5.21]{OV}).
\medskip

\begin{thm} \label{th:tloi} In addition to {\bf (A.1--4)}, assume that the matrix $P^2$ is irreducible and that there exists a constant
$\kappa>0$ such that for all $x,y \in [0,1]^2$ and
$i \in E$
\begin{equation}
\label{hyp:tloi}
\sum_{u \in \Lambda, u  \perp x-y} \bigl| p(x,i,u)-p(y,i,u) \bigr|
+ \sum_{u \in \Lambda} \bigl(p(x,i,u) - p(y,i,u)\bigr)
{\rm sgn}\bigl(\langle x-y,u \rangle \bigr) \leq -\kappa |x-y|_1
\end{equation}
(As usual, ${\rm sgn}(\cdot)$ denotes the sign function, with ${\rm sgn}(u)=u/|u|$ for $u \neq 0$ and
${\rm sgn}(0)=0$.). Define $T_m^i=\min\{t>0: \,
\bbP_0^i ( \tau^{(m)} > t)\leq e^{-1}\}$ for $i \in E$ and $m \geq 1$.
Then, for any sequence of starting points $(x_m)_{m \geq 1}$ in G, with $x_m \rightarrow x \in G$
as $m \rightarrow + \infty$,
\begin{center}
the law of $\tau^{(m)} / T_m^i$ under $\P^i_{x_m}$ weakly converges to an 
exponential law of mean 1.
\end{center}
\end{thm}
In what follows, we will prove that, for any $i,j \in E$, $T_m^i/T_m^j \rightarrow 1$ as $m$ tends to $+ \infty$. In particular, the image law
$(\tau^{(m)} / T_m^j)(\P^i_{x_m})$ weakly converges to an exponential law of mean 1 for any $i,j \in E$.

Condition \eqref{hyp:tloi} is not empty : Example
(\ref{eq:maier}) given below fulfills \eqref{hyp:tloi} if $g_1,g_2$ are strictly increasing with $g_1', g_2' \geq \kappa'$ a.e. for some $\kappa'>0$.

$\Box$ Proof. The following result
is the analogue of \cite[Lemma 5.22]{OV}.
Its proof is deferred to Section
\ref{sec:preuve-couplage},
\begin{lem}
\label{couplage}
There exists $\delta >0$, such that, for all $i \in E$ and $S>0$,
\begin{equation*}
\lim_{m \rightarrow + \infty} \sup_{|x|_1,|y|_1 \leq \delta m,|x-y|_1 \in
  2{\mathbb N}} \sup_{t \geq S} |\P_{x/m}^i \{ \tau^{(m)} > tm^2 \} -
\P_{y/m}^i \{\tau^{(m)} > tm^2 \} |=0.
\end{equation*}
\end{lem}

With this lemma at hand, we can prove

\begin{lem}
\label{lem:chgt_depart}
For all $\eta >0$ and $S>0$, we can find a sequence $(\varepsilon_m)_{m \geq 1}$ of positive reals, tending to 0 as $m \rightarrow+ \infty$, such that
for all $i, j \in E$,
\begin{equation*}
\forall t\geq S, \ {\mathbb P}_0^i \{\tau^{(m)} > t m^2\} \leq {\mathbb P}_0^j \{\tau^{(m)} > tm^2 - \eta m\} + \varepsilon_m
\end{equation*}
\end{lem}

\noindent $\Box$ Proof of Lemma \ref{lem:chgt_depart}. For $i \in
E$, we set $\vartheta_i= \inf \{ k \in 2 {\mathbb N}: \, X_k = i\}$.
Since $P^2$ is assumed to be irreducible, it is a finite stopping
time. For $\delta$ as in Lemma \ref{couplage} and $\eta < \delta$,
\begin{equation*}
\begin{split}
\P_0^i \{\tau^{(m)} > t m^2 \}
&\leq \P_0^i \{\tau^{(m)} > t m^2, \vartheta_j <
\eta m \} + \P_0^i \{\vartheta_j \geq \eta m\}
\\
&\leq \sup_{\{|x|_1 \leq \delta m, |x|_1\in 2 {\mathbb N}\}}
\P_{x/m}^j \{\tau^{(m)}
> t m^2 - \eta m\}+ \P_0^i \{\vartheta_j \geq \eta m\}.
\end{split}
\end{equation*}
It is clear that $\lim_{m \rightarrow + \infty} \P_0^i \{\vartheta_j \geq \eta m\}=0$.
By Lemma \ref{couplage}, the proof is easily completed.\qed

We now complete the proof of Theorem
\ref{th:tloi}. We keep the notations introduced in the proof of Theorem \ref{th:tmoy}. Following \cite[Lemma 5.23]{OV}, we can set for $i \in E$
\begin{equation*}
\forall t \geq 0, \ F^{(m),i}(t) = \P_0^i\{\tau^{(m)} > t T_m^i\} =
\P_0^i\{\tau^{(m)} > m^{-1} \lfloor m t T_m^i \rfloor \}.
\end{equation*}
By Theorem \ref{th:tmoy}, for every $\delta >0$, we have $\lim_{m
\rightarrow + \infty} T_m^i \exp[-m(\bar{V}-\delta)] = + \infty$ and
$\lim_{m \rightarrow + \infty} T_m^i \exp[-m(\bar{V}+\delta)]=0$.
Moreover, by the Markov
property, for $j \in E$, $\rho < \ell$ and $t >0$,
\begin{equation*}
\begin{split}
&\sup_{x_m \in G^{(m)}} \P_{x_m}^j \{ \tau^{(m)} > m^{-1} \lfloor m t T_m^i \rfloor,
\sigma_{\rho} < m^{-1} \lfloor T_m^i \rfloor \}
\\
&\leq \sup_{x_m \in G^{(m)}}\P_{x_m}^j \{\sigma_{\rho} < m^{-1} \lfloor T_m^i \rfloor \}
\sup_{y \in
 \bar{B}_{\rho}^{(m),+}, k \in E}
\P_y^k \{ \tau^{(m)} > m^{-1}(\lfloor m t T_m^i \rfloor - \lfloor T_m^i \rfloor) \}
\\
&\leq  \sup_{y \in \bar{B}_{\rho}^{(m),+},k \in E}
\P_y^k \{\tau^{(m)} >  m^{-1}(\lfloor m t T_m^i \rfloor - \lfloor T_m^i \rfloor)
 \}.
\end{split}
\end{equation*}
In the above supremum, 
we aim at applying Lemma \ref{couplage} to the starting points $0$ and $y$ ($\rho$ being small enough).
There is no difficulty if $|y|_1 \in 2 m^{-1} {\mathbb N}$.
If $|y|_1 \in (2 m^{-1}{\mathbb N} + m^{-1})$, the Markov property yields
$\P_y^k \{\tau^{(m)} >m^{-1} (\lfloor m t T_m^i \rfloor - \lfloor T_m^i \rfloor)
 \} \leq \sup_{|z-y|_1=1/m, k' \in E} \P_z^{k'} \{\tau^{(m)} >m^{-1} (\lfloor m t T_m^i \rfloor - \lfloor T_m^i \rfloor -1)
 \}$, so that we can still apply Lemma \ref{couplage}. By Lemma \ref{lem:5.7.19}, we deduce that we can choose $\rho$ small enough and find
 some sequence
$(\delta_{m})_{m \geq 1}$ with $\lim_{m \rightarrow + \infty} \delta_{m} = 0$
such that
\begin{equation}
\label{couplage_fin_0}
\begin{split}
\sup_{x_m \in G^{(m)}} \P_{x_m}^j \{ \tau^{(m)} > m^{-1} \lfloor m t
T_m^i \rfloor  \} &\leq \sup_{k \in E} \P_0^k \{\tau^{(m)} >  m^{-1}(\lfloor m t T_m^i
\rfloor - \lfloor T_m^i \rfloor -1)
 \} + \delta_m
 \\
 &\leq \P_0^i \{\tau^{(m)} >  m^{-1}(\lfloor m t T_m^i
\rfloor - 2 \lfloor T_m^i \rfloor)
 \} + \delta_m,
\end{split}
\end{equation}
the second line following from Lemma \ref{lem:chgt_depart}.
The Markov property yields for $t,s > 0$,
\begin{equation}
\label{couplage_fin_2}
\begin{split}
&{\mathbb P}_0^i\{\tau^{(m)} > m^{-1} \lfloor m(t+s) T_m^i \rfloor \}
\\
&\leq {\mathbb P}_0^i  \{\tau^{(m)} > m^{-1} (\lfloor m (t+s)
T_m^i \rfloor - \lfloor mt T_m^i \rfloor)\}
\sup_{x_m \in G^{(m)},j \in E} \P_{x_m}^j \{\tau^{(m)} > m^{-1}\lfloor m t
T_m^i \rfloor  \}
\\
&\leq  {\mathbb P}_0^i  \{\tau^{(m)} > m^{-1} (\lfloor m (t+s)
T_m^i \rfloor - \lfloor mt T_m^i \rfloor)\} \P_0^i \{\tau^{(m)} >  m^{-1}(\lfloor m t T_m^i
\rfloor - 2\lfloor T_m^i \rfloor)
 \}
\\
&\hspace{15pt} + \delta_m.
\end{split}
\end{equation}
We can prove the converse inequality in a similar way. For any compact subset $K \subset G$, we
deduce from Lemmas \ref{lem:5.7.19}, \ref{lem:5.7.22},
\ref{couplage} and \ref{lem:chgt_depart} that, for $tm >1$, (up to a modification from line to line of the sequence $(\delta_m)_{m
\geq 1}$ -- which may depend on $K$ -- )
\begin{equation}
\begin{split}
\label{couplage_fin_1} &\inf_{x_m \in K \cap G^{(m)}} \P_{x_m}^j \{ \tau^{(m)} > m^{-1} \lfloor m t
T_m^i \rfloor  \}
\\
&\geq \inf_{x_m \in K \cap G^{(m)}} \P_{x_m}^j \{ \tau^{(m)} > m^{-1} \lfloor m t T_m^i \rfloor,
\sigma_{\rho} < m^{-1} \lfloor T_m^i \rfloor,
\bar{X}_{\sigma_{\rho}}^{(m)} \in \bar{B}_{\rho}^{(m),+} \}
\\
&\geq  \inf_{x_m \in K \cap G^{(m)}} \P_{x_m}^j \{\sigma_{\rho} < m^{-1} \lfloor T_m^i \rfloor,
\bar{X}_{\sigma_{\rho}}^{(m)} \in \bar{B}_{\rho}^{(m),+} \}
\\
&\hspace{15pt} \times \inf_{y
\in \bar{B}_{\rho}^{(m),+},k \in E} \P_y^k \{ \tau^{(m)} > m^{-1} \lfloor m t T_m^i
\rfloor  \}
\\
&\geq (1-\delta_{m}) \bigl[ \P_0^i \{\tau^{(m)} > m^{-1}(\lfloor m t T_m^i
\rfloor + \lfloor T_m^i \rfloor)
 \} - \delta_m \bigr]
 \\
 &\geq \P_0^i \{\tau^{(m)} > m^{-1}(\lfloor m t T_m^i
\rfloor + \lfloor T_m^i \rfloor) \} - \delta_m.
\end{split}
\end{equation}
Now, for $t,s>0$, \eqref{couplage_fin_1} yields
\begin{equation}
\label{exp_law:1}
\begin{split}
&{\mathbb P}_0^i\{\tau^{(m)} > m^{-1} \lfloor m(t+s) T_m^i \rfloor \}
\\
&\geq {\mathbb P}_0^i  \{\tau^{(m)} > m^{-1} (\lfloor m (t+s)
T_m^i \rfloor - \lfloor mtT_m^i \rfloor),{\bar X}^{(m)}_{ m^{-1} (\lfloor m(t+s) T_m^i \rfloor - \lfloor mt T_m^i\rfloor )} \in K
\}
\\
&\hspace{15pt} \times \inf_{x_m \in K \cap G^{(m)},j \in E} \P_{x_m}^j
\{\tau^{(m)} > m^{-1} \lfloor mt T_m^i \rfloor  \}
\\
&\geq {\mathbb P}_0^i  \{\tau^{(m)} > m^{-1} (\lfloor m (t+s)
T_m^i \rfloor - \lfloor mtT_m^i \rfloor),{\bar X}^{(m)}_{ m^{-1} (\lfloor m(t+s) T_m^i \rfloor - \lfloor mt T_m^i\rfloor )} \in K
\}
\\
&\hspace{15pt} \times \P_0^i \{\tau^{(m)} > m^{-1}(\lfloor m t T_m^i
\rfloor + \lfloor T_m^i \rfloor) \} - \delta_m.
\end{split}
\end{equation}
By {\bf (A.4)}, for any starting point $x_0 \in \bar{G}$, $\chi^{x_0}_t \in G$  for $t>0$. In particular, $d(\chi^{x_0}_1,\partial G)>0$. By the stability property for RDEs driven
by Lipschitz continuous coefficients, we have $\inf_{x_0 \in \bar{G}} d(\chi^{x_0}_1,\partial G) >0$. In other words, we can find a compact subset $K_0 \subset G$ such that
$\chi^{x_0}_1 \in K_0$ for any $x_0 \in \bar{G}$. We denote by $\varepsilon = d(K_0,\partial G) >0$ the distance from $K_0$ to $\partial G$. By
Corollary \ref{conv_pb},
\begin{equation}
\label{exp_law:2}
\lim_{m \rightarrow + \infty} \sup_{x_m \in G^{(m)},j \in E} \P_{x_m}^j \{ d(\bar{X}_1^{(m)},\partial G)  \leq \varepsilon/2 \} =0.
\end{equation}
By the Markov property,
\begin{equation*}
\begin{split}
&{\mathbb P}_0^i  \{\tau^{(m)} > m^{-1} (\lfloor m (t+s)
T_m^i \rfloor - \lfloor mt T_m^i \rfloor), d({\bar X}^{(m)}_{ m^{-1} (\lfloor m(t+s) T_m^i \rfloor - \lfloor mtT_m^i \rfloor)},\partial G) \geq \varepsilon/2
\}
\\
&= {\mathbb P}_0^i  \{\tau^{(m)} > m^{-1} (\lfloor m (t+s) T_m^i \rfloor - \lfloor mt T_m^i \rfloor) \}
\\
&\hspace{15pt} - {\mathbb P}_0^i \{
\tau^{(m)} > m^{-1} (\lfloor m (t+s)
T_m^i \rfloor - \lfloor m t T_m^i \rfloor)
,
d({\bar X}^{(m)}_{ m^{-1} (\lfloor m(t+s) T_m^i \rfloor - \lfloor mtT_m^i \rfloor)},\partial G) < \varepsilon /2\}
\\
&\geq {\mathbb P}_0^i  \{\tau^{(m)} > m^{-1} (\lfloor m (t+s)
T_m \rfloor - \lfloor mt T_m\rfloor)\} - \sup_{x_m \in G^{(m)},j \in E} \P_{x_m}^j \{ d(\bar{X}_1^{(m)},\partial G)  \leq \varepsilon/2 \}.
\end{split}
\end{equation*}
We can plug $K = \{z \in G : \, d(z,\partial G) \geq \varepsilon/2\}$ in \eqref{exp_law:1}. By \eqref{exp_law:2} and the above inequality,
\begin{equation}
\label{exp_law:3}
\begin{split}
&{\mathbb P}_0^i\{\tau^{(m)} > m^{-1} \lfloor m(t+s) T_m^i \rfloor \} \geq
{\mathbb P}_0^i\{\tau^{(m)} > m^{-1} (\lfloor m(t+s) T_m^i \rfloor-\lfloor mt T_m^i \rfloor) \}
\\
&\hspace{15pt}
\times {\mathbb P}_0^i\{\tau^{(m)} > m^{-1} (\lfloor mt T_m^i \rfloor + \lfloor T_m^i \rfloor) \} - \delta_m.
\end{split}
\end{equation}
By
\eqref{couplage_fin_2} and \eqref{exp_law:3},
\begin{equation*}
F^{(m),i}(s + \delta_m) F^{(m),i}(t+ \delta_m) - \delta_m  \leq
F^{(m),i}(t+s) \leq F^{(m),i}(t-\delta_m)  F^{(m),i}(s-\delta_m) + \delta_m,
\end{equation*}
so that $\limsup_{m \rightarrow + \infty} F^{(m),i}(k+\varepsilon) \leq e^{-k}$ for $k \in {\mathbb N}$ and $\varepsilon >0$.
In particular, the sequence
$\tau^{(m)}/T_m^i$ is tight. Up to a subsequence, it converges in
law. The limit distribution function is denoted by $F$.
Up to a countable subset of $(0,+\infty)$, $F^{(m),i}(t)$ converges to $F(t)$.
Hence, we can pass to the limit in the above inequality. For all
$\eta >0$,
\begin{equation*}
F(t+\eta) F(s+\eta) \leq F(t+s) \leq F(t-\eta) F(s-\eta).
\end{equation*}
It is plain to deduce that the limit distribution is the exponential
law with mean one. By \eqref{couplage_fin_0} and
\eqref{couplage_fin_1}, this is true for any starting point.
Moreover, for all $j \in E$, $(\tau^{(m)}/T_m^i)({\mathbb P}^j_0)$ weakly converges to the exponential law with mean one.
Since $(\tau^{(m)}/T_m^j)({\mathbb P}^j_0)$ weakly converges to the same distribution, we deduce that $T_m^i/T_m^j \rightarrow 1$ as $m \rightarrow + \infty$.
\qed

\subsection{Hamilton-Jacobi Equation for the Quasi-Potential}
\label{sec:comput-quasipot}

In practice, it is important to compute the quasi-potential
$V(0,x)$ as well as  the optimal paths. (In what follows, we write, for the sake of simplicity, $V(x)=V(0,x)$.)

In \cite[Chapter 5, Theorem 4.3]{FV}
and \cite[Exercise 5.7.36]{DZ}, it is shown that the quasi-potential
is characterized through a Hamilton-Jacobi equation of the form
$$
H\big(x,\nabla V(x)\big)=0.
$$

Loosely speaking, the equation for the quasi-potential has the same
structure in our setting. However, due to the reflection phenomenon, it satisfies some specific
boundary condition.

{\bf Form of the Equation.} Here, we specify both the
equation and the boundary condition in the viscosity sense, the
notion of viscosity solutions being, in a general way, particularly
well adapted to optimal control problems. (See for example
\cite{bardi:capuzzo} or \cite{barles} for a review on this connection.) Indeed, the
quasi-potential is nothing but the value function of some optimal
control problem. In the formula \eqref{eq:cost}, $L(\phi_s,\psi_s)$
may be interpreted as some instantaneous cost at time $s$ when the
trajectory $\phi$ is driven by the control $\psi$. The controlled
dynamical system obeys the rule: $\forall t \geq 0$, $\phi_t =
\psi_t - k_t$, with $k$ as in \eqref{Skorokhod}.

\begin{prop}
\label{prop:eqHJB}
We assume that {\bf (A.1--3)} are in force. Then, for every $x \in (0,1)^d$ and every continuously differentiable function $\theta$ on
a neighborhood $U \subset (0,1)^d$
of $x$,
\begin{equation}
\label{eq:eqHJB}
\begin{split}
&H\bigl(x,\nabla \theta(x) \bigr)
\leq  0 \ {\rm if} \ V-\theta \ {\rm has \ a \ local \ maximum \ at} \ x,
\\
&H\bigl(x,\nabla \theta(x) \bigr)
= 0 \ {\rm if} \ V-\theta \ {\rm has \ a \ local \ minimum \ at} \ x.
\end{split}
\end{equation}
Moreover, for every $x \in \partial [0,1]^d$ and every continuously differentiable function $\theta$ on $U \cap [0,1]^d$, $U$ being a neighborhood of $x$,
\begin{equation}
\label{eq:eqHJBound}
\begin{split}
&H\bigl(x,\nabla \theta(x) \bigr)
\geq 0 \ {\rm if} \
\left\{
\begin{array}{l}
\forall n \in {\mathcal N}(x), \  \langle \nabla \theta(x),n \rangle \geq 0,
\\
V-\theta \ {\rm has \ a \ local \ minimum \ at} \ x \ {\rm on} \ U \cap [0,1]^d,
\end{array}
\right.
\\
&H\bigl(x,\nabla \theta(x) \bigr)
\leq 0 \ {\rm if} \
\left\{
\begin{array}{l}
\forall n \in {\mathcal N}(x), \  \langle \nabla \theta(x),n \rangle \leq 0,
\\
V-\theta \ {\rm has \ a \ local \ minimum \ at} \ x \ {\rm on} \ U \cap [0,1]^d.
\end{array}
\right.
\end{split}
\end{equation}
\end{prop}

The asymmetry between the two conditions in \eqref{eq:eqHJB} is
standard in the theory of optimal control. The first line says that
$V$ is a viscosity subsolution of the Hamilton-Jacobi equation in
$(0,1)^d$, the second one that $V$ is a bilateral supersolution.
Generally speaking, $V$ is also a bilateral subsolution at $x  \in
(0,1)^d$, i.e. $H(x,\nabla \theta(x))=0$ if $V-\theta$ has a local maximum at $x$,
if there exists an optimal path reaching $x$. We refer the
reader to \cite[\S 2.3, Chapter III]{bardi:capuzzo} for more
details.

The boundary condition \eqref{eq:eqHJBound} is a boundary condition
of Neumann type. This Neumann condition expresses the reflected
structure of the controlled dynamical system. The viscosity
formulation of the Neumann boundary condition has been introduced in
\cite{lions:85}. In what follows, we will explain the link between
this weak formulation and the standard Neumann condition.
\\

\noindent $\Box$ {Proof.} The proof is standard.
We first give a suitable version of the Bellman dynamic programming principle for the quasi-potential
$V$. Then, we will deduce Proposition \ref{prop:eqHJB}.

\begin{lem}
\label{lem:DPP}
For all $x \in [0,1]^d$, for all $t > 0$,
\begin{equation}
\label{eq:dpp0} V(x) = \inf \bigl\{ V(y) + \int_0^t
L(\phi_s,\dot{\psi}_s) ds; \ (y,\phi,\psi) : \phi_0=y, \ \phi_t=x, \
\phi=\Psi(\psi) \bigr\}.
\end{equation}
\end{lem}
(In the above formula, we can assume that $|\dot{\psi}_s| \leq 1$
for a.e. $s \in [0,t]$ since $L(x,v)=+\infty$ for $|v| > 1$. In
particular, we can assume that $|\dot{\phi}_s| \leq 1$ and $|\dot{\psi}_s -
\dot{\phi}_s| \leq 1$ for a.e. $s \in [0,t]$. Indeed, $|\dot{\psi}_s|^2
=|\dot{\phi}_s|^2 + |\dot{\psi}_s - \dot{\phi}_s|^2$ for a.e. $s \in [0,t]$.)

The proof of Lemma \ref{lem:DPP} is left to the reader. Details may
be found in \cite[Proposition 2.5, Chapter III]{bardi:capuzzo}.

With the Bellman dynamic programming principle at hand, it is
standard to prove that $V$ is both a subsolution and a supersolution
at $x \in (0,1)^d$, i.e. $H(x,\nabla \theta(x)) \leq 0$ if
$V-\theta$ has a local maximum at $x$ and $H(x,\nabla \theta(x))
\geq 0$ if $V-\theta$ has a local minimum at $x$. (See for example
the proof of
 \cite[Proposition 2.8, Chapter III]{bardi:capuzzo}.)

We now investigate the first boundary condition.

For a given $x \in \partial [0,1]^d$, we assume that there exists a
continuously differentiable function $\theta$ on $U \cap [0,1]^d$,
$U$ being a neighborhood of $x$, such that $V-\theta$ has a local
minimum at $x$ on $U \cap [0,1]^d$. Without loss of generality, we
can assume that $\theta(x)=V(x)$ and that the minimum is global on
$U \cap [0,1]^d$ so that $V(y)-\theta(y) \geq 0$ for all $y \in U
\cap [0,1]^d$. We also assume  $\langle \nabla \theta(x),n \rangle
\geq 0$ for all $n \in {\mathcal N}(x)$.

For $t$ small, we can  assume that $y \in U$ in the dynamic
programming principle. We deduce that, for all $t$ small,
\begin{equation*}
\theta(x) \geq \inf \bigl\{ \theta(y) + \int_{0}^t L(\phi_s,\dot{\psi}_s)ds\bigr\},
\end{equation*}
the infimum being taken over the same triples as above. Developing
$\theta(x)- \theta(y)$, we can write
\begin{equation}
\label{eq:visco1}
\begin{split}
\sup \bigl\{ \int_{0}^t \langle \nabla \theta (\phi_s),\dot{\psi}_s
\rangle ds - \int_{0}^t \langle \nabla \theta (\phi_s),\dot{\psi}_s
&- \dot{\phi}_s  \rangle ds - \int_{0}^t L(\phi_s,\dot{\psi}_s)ds;
\\
&\hspace{15pt}  \phi_0=y, \ \phi_t=x, \ \phi=\Psi(\psi) \bigr\} \geq
0.
\end{split}
\end{equation}
Having in mind that $\dot{\psi}_s - \dot{\phi}_s \in {\mathbb R}_+
{\mathcal N}(\phi_s)$  (with ${\mathbb R}_+ {\mathcal N}(\phi_s)=
\{0\}$ if $\phi_s \in (0,1)^d$) and $|\dot{\psi}_s - \dot{\phi}_s|
\leq 1$ for a.e. $s \in [0,t]$, we deduce
\begin{equation*}
\begin{split}
\sup \bigl\{ \int_{0}^t &\sup_{n \in {\mathcal N}(\phi_s)} \max
\bigl( 0, - \langle \nabla \theta (\phi_s),n \rangle \bigr) +
\int_{0}^t H\bigl(\phi_s,\nabla \theta(\phi_s) \bigr) ds;
\\
&\phi : [0,t] \rightarrow [0,1]^d, \ \phi_t=x, \
|\dot{\phi}_s| \leq 1 \ {\rm for \ a.e.} \ s \bigr\} \geq 0,
\end{split}
\end{equation*}
Despite the lack of regularity of the boundary of $[0,1]^d$, we can
prove that, for $|z-x|$ small enough, ${\mathcal N}(z) \subset
{\mathcal N}(x)$. Since $\nabla \theta$ and $H$ are continuous, we
deduce
\begin{equation*}
 \sup_{n \in {\mathcal N}(x)} \max \bigl( 0, - \langle \nabla \theta (x),n \rangle \bigr)
+  H\bigl(x,\nabla \theta(x) \bigr)  + \varepsilon_t  \geq 0,
\end{equation*}
with $\varepsilon_t \rightarrow 0$ as $t$ tends to 0.
By assumption, the first term in the above left-hand side is zero. This completes the proof.

We now prove that $V$ is a bilateral supersolution in $(0,1)^d$ and
satisfies the second boundary condition. The idea follows from
\cite[\S 2.3, Chapter III]{bardi:capuzzo} and consists in reversing
the dynamic programming principle. This permits to write $x$ as the
initial condition of the controlled trajectory $\phi$.

We let the reader check that for all $x \in [0,1]^d$ and for all $t>0$,
\begin{equation}
\label{eq:dpp00} V(x) \geq \sup \bigl\{ V(y) - \int_0^t
L(\phi_s,\dot{\psi}_s) ds; \ (y,\phi,\psi) : \phi_0=x, \ \phi_t=y, \
\phi=\Psi(\psi) \bigr\}.
\end{equation}
(Pay attention: there is no equality in \eqref{eq:dpp00} at this stage of the paper. 
Equality holds if there exists an optimal path from $0$ to $x$.
This is the reason why we are not able to prove that $V$ is a bilateral subsolution of the
Hamilton-Jacobi equation.)

Following \cite[Proposition 2.8]{bardi:capuzzo}, this shows that $V$ is a bilateral supersolution of the Hamilton-Jacobi equation in $(0,1)^d$.

We now prove the second boundary condition. As above,
 we assume that there exists a continuously differentiable function $\theta$ on $U \cap [0,1]^d$, $U$ being
a neighborhood of $x$, such that $\theta(x)=V(x)$ and
$V(y)-\theta(y) \geq 0$ for all $y \in U \cap [0,1]^d$. We also
assume $\langle \nabla \theta(x),n \rangle \leq 0$ for all $n \in
{\mathcal N}(x)$.

We choose a control $\psi$ with a constant speed. For $\alpha \in
\R^d$, we choose $\psi_s=x + \alpha s$ for all $s \in [0,t]$. We
then define $\phi = \Psi(\psi)$. By \eqref{Skorokhod}, we can write
$\phi_s=x + \alpha s - k_s$, with $\dot{k}_s \in {\mathbb
R}_+{\mathcal N}(\phi_s)$. For $t$ small enough, $\phi_t$ is in $U$
and \eqref{eq:dpp00} yields
\begin{equation*}
\theta(x) \geq \theta(\phi_t) -  \int_0^t L\bigl(\phi_s, \alpha \bigr)ds.
\end{equation*}
Developing $\theta(\phi_t) - \theta(x)$ as in \eqref{eq:visco1}, we obtain
\begin{equation*}
\int_0^t \langle \nabla \theta(\phi_s),\alpha \rangle ds - |\alpha| \int_0^t
\sup_{n \in {\mathcal N}(\phi_s)} \max(0,\langle \nabla
\theta(\phi_s),n \rangle ) ds - \int_0^t L\bigl(\phi_s, \alpha
\bigr)ds \leq 0.
\end{equation*}
As above, we obtain
\begin{equation*}
\langle \nabla \theta(x),\alpha \rangle  - L(x,\alpha) -  |\alpha| \sup_{n \in {\mathcal N}(x)} \max(0,\langle \nabla \theta(x),n \rangle)
 \leq  0.
\end{equation*}
By assumption, $\sup_{n \in {\mathcal N}(x)} \max(0,\langle \nabla \theta(x),n \rangle)=0$. We deduce $H(x,\nabla \theta(x)) \leq 0$. \qed

We now explain the form of the equation when the quasi-potential is
continuously differentiable on $[0,1]^d \setminus \{0\}$. (We
exclude $0$ from the set of differentiable points because there is a 
boundary condition of Dirichlet type in 0:
 $V(0)=0$. Anyhow, as seen in the next section, there are
specific examples in which $V$ is continuously differentiable on the
whole $[0,1]^d$.) To this end, we introduce a modification of the
gradient at the boundary. Assuming that $\nabla V$ exists at $x \in
\partial [0,1]^d \setminus \{0\}$, we set
\begin{equation*}
\forall i \in \{1,\dots,d\}, \ \bigl( \nabla_{+} V(x) \bigr)_i =
\left\{
\begin{array}{l}
\bigl[ \partial V/\partial x_i \bigr](x) \ {\rm if} \ 0 < x_i<1,
\\
\min \bigl(\bigl[ \partial V/\partial x_i \bigr](x),0 \bigr) \ {\rm
if} \ x_i= 0,
\\
\max \bigl(\bigl[ \partial V/\partial x_i \bigr](x),0
\bigr) \ {\rm if} \ x_i= 1.
\end{array}
\right.
\end{equation*}
Similar modifications of the gradient of the quasi-potential appear
in \cite[Section II]{lions:85}. Following the notations introduced
there in, we give another writing for $\nabla_{+}V(x)$. We denote by
$\nabla_T V(x)$ the tangential part of $\nabla V(x)$, i.e.
\begin{equation*}
\forall u \perp {\mathcal N}(x), \ \langle \nabla_T V(x),u \rangle =
\langle \nabla V(x),u \rangle, \quad \forall n \in {\mathcal N}(x),
\ \langle \nabla_T V(x),n \rangle= 0.
\end{equation*}
We also denote by ${\bf e}(x)$ the set ${\mathcal N}(x) \cap
{\mathcal V}$, so that ${\bf e}(x)$ is an orthonormal basis of the
cone generated by ${\mathcal N}(x)$. (It satisfies $\langle e,n
\rangle \geq 0$ for all $e \in {\bf e}(x)$ and $n \in {\mathcal
N}(x)$.) Then, $\nabla V(x)$ may be expressed as
\begin{equation}
\label{eq:dgrad1} \nabla V(x) = \nabla_T V(x) + \sum_{e \in {\bf
e}_+(x)}   \langle \nabla V(x),e \rangle  e +  \sum_{e \in {\bf
e}_-(x)} \langle \nabla V(x),e \rangle  e,
\end{equation}
with ${\bf e}_+(x) = \{e \in {\bf e}(x), \ \langle \nabla V(x),e
\rangle > 0\}$ and ${\bf e}_-(x) = \{e \in {\bf e}, \ \langle \nabla
V(x),e \rangle < 0\}$. (In what follows, we will also make use of
${\bf e}_0(x) = \{e \in {\bf e}(x), \ \langle \nabla V(x),e \rangle
= 0\}$.) With these notations at hand, we have
\begin{equation}
\label{eq:dgrad2} \nabla_+ V(x) = \nabla_T V(x) + \sum_{e \in {\bf
e}_+(x)}   \langle \nabla V(x),e \rangle  e.
\end{equation}
The above expression justifies the notation $\nabla_+ V(x)$. We are
now ready to state:
\begin{prop}
\label{prop:HJBclassique} Assume {\bf (A.1--3)}. If the quasi-potential $V$ is continuously
differentiable on $[0,1]^d \setminus \{0\}$, then it satisfies
\begin{equation}
\label{eq:HJBclass1}
\forall x \in (0,1)^d, \ H\bigl(x,\nabla V(x) \bigr)=0,
\end{equation}
with the boundary condition
\begin{equation}
\label{eq:HJBclass2} \forall x \in \partial [0,1]^d \setminus \{0\},
\ H \bigl(x,\nabla_{+} V(x) \bigr)=0.
\end{equation}
\end{prop}
\noindent

By continuity of $\nabla V$, we notice that \eqref{eq:HJBclass1}
holds for all $x \in [0,1]^d \setminus \{0\}$. Moreover, we
emphasize that \eqref{eq:HJBclass2} is a boundary condition of
Neumann type. If $\nabla V$ satisfies the standard Neumann
condition, i.e. $\langle \nabla V(x),n \rangle =0$ for all $n \in
{\mathcal N}(x)$, at some $x \in \partial [0,1]^d \setminus \{0\}$,
then $\nabla_+ V(x)$ and $\nabla V(x)$ are equal. In this case,
\eqref{eq:HJBclass2} follows from \eqref{eq:HJBclass1}.

As explained in  \cite[Section II]{lions:85}, Hamilton-Jacobi
equations under the standard Neumann condition, i.e. $\langle \nabla
V(x),n \rangle =0$ for all $n \in {\mathcal N}(x)$ and $x \in
\partial [0,1]^d \setminus \{0\}$, may not be well-posed. This explains why a weaker formulation
of the boundary condition may be necessary. Anyhow,
\eqref{eq:HJBclass2} is slightly different from the Neumann
condition given in \cite[Section II]{lions:85} since the original
formulations in terms of viscosity solutions are different. (The
optimal control problems are a bit different.) Moreover, the
existence of ``angles'' along the hypercube $[0,1]^d$ induces
additional difficulties in our framework. (In comparison, the
boundary is assumed to be smooth in \cite[Section II]{lions:85}.)
\\

\noindent $\Box$ Proof. The proof is obvious inside the domain.
(Choose $\theta=V$ in the statement.)

To prove the boundary condition, we characterize the continuously
differentiable functions $\theta$ such that $V-\theta$ has a local
minimum at $x \in \partial [0,1]^d \setminus \{0\}$. Following the
proof of \cite[Lemma 1.7, Chapter II]{bardi:capuzzo}, for a given $p \in
{\mathbb R}^d$, there exists a continuously differentiable function
$\theta$ (on a neighborhood of $x$) such that $V-\theta$ has a local
minimum at $x$ and $\nabla \theta(x)=p$ if and only if 
the tangential part $p_T$ of $p$ is equal to $\nabla_T V(x)$
and $\langle p,n
\rangle \geq \langle \nabla V(x),n \rangle$ for all $n \in {\mathcal
N}(x)$.

In what follows, the typical value of $p$ is $p=\nabla_{+} V(x)$.
Indeed, $(\nabla_+ V(x))_T = \nabla_T V(x)$ and
$\langle \nabla_{+}V(x),n \rangle \geq \langle \nabla V(x),n
\rangle$ for all $n \in {\mathcal N}(x)$.

If $\langle \nabla V(x),n \rangle \leq 0$ for all $n \in {\mathcal
N}(x)$, then $\langle \nabla_{+} V(x),n \rangle =0$ for all $n \in {\mathcal N}(x)$. (See
\eqref{eq:dgrad2}.) Hence, we can apply both conditions in
\eqref{eq:eqHJBound}. We deduce that
\begin{equation}
\label{eq:HJBclass3} H \bigl(x,\nabla_{+} V(x) \bigr) = 0.
\end{equation}

On the contrary, if $\langle \nabla V(x),n \rangle \geq 0$ for all $n \in
{\mathcal N}(x)$, the result is obvious. Indeed, $\nabla_{+} V(x)=
\nabla V(x)$ in this case. Since $\nabla V$ is continuous, the
Hamilton-Jacobi equation \eqref{eq:HJBclass2} is true up to the
boundary.

The intermediate cases may be treated by a similar argument of
continuity. With the expressions \eqref{eq:dgrad1} and
\eqref{eq:dgrad2} at hand, we set $n_+(x) = \sum_{e \in {\bf e}_+(x)
\cup {\bf e}_0(x)} e$ and, for $\varepsilon>0$, $y_{\varepsilon} = x
- \varepsilon n_+(x)$. For $\varepsilon$ small enough, we have
${\mathcal N}(y_{\varepsilon})= {\bf e}_-(x)$ (if ${\mathbf e}_-(x)$
is empty then $y_{\varepsilon} \in (0,1)^d$ and ${\mathcal
N}(y_{\varepsilon})$ is also empty) and, by continuity of $\nabla
V$, $\langle \nabla V(y_{\varepsilon}),e \rangle < 0$ for all $e \in
{\bf e}_-(x)$. By \eqref{eq:HJBclass3},
\begin{equation*}
H \bigl(y_{\varepsilon},\nabla_{+} V(y_{\varepsilon}) \bigr) = 0.
\end{equation*}
As $\varepsilon$ tends 0, $\nabla_{+} V(y_{\varepsilon})$ tends to
$\nabla_{+} V(x)$. Indeed, by \eqref{eq:dgrad1} and
\eqref{eq:dgrad2},
\begin{equation*}
\begin{split}
\nabla_+ V(y_{\varepsilon}) = \nabla V(y_{\varepsilon}) - \sum_{e
\in {\bf e}_-(y_{\varepsilon})} \langle \nabla V(y_{\varepsilon}),e
\rangle e &= \nabla V(y_{\varepsilon}) - \sum_{e \in {\bf e}_-(x)}
\langle \nabla V(y_{\varepsilon}),e \rangle e
\\
&\underset{\varepsilon\rightarrow 0}{\longrightarrow} \nabla V(x) -
\sum_{e \in {\bf e}_-(x)} \langle \nabla V(x),e \rangle e = \nabla_+
V(x).
\end{split}
\end{equation*} This completes the
proof. \qed

The boundary conditions are not formulated in a complete way in
Proposition \ref{prop:HJBclassique}. As stated below,
\eqref{eq:eqHJBound} implies additional conditions on the
derivatives $\nabla_{\alpha} H(x,\nabla V(x))$ and $\nabla_{\alpha}
H(x,\nabla_{+} V(x))$. In \cite[Section II]{lions:85}, these
additional conditions are formulated in a different way: the
formulation used there in is about the signs of $H(x,\nabla
V(x)+\lambda n)$ and $H(x,\nabla_{+}V(x)+\lambda n)$ for $\lambda
\in {\mathbb R}$ and $n \in {\mathcal N}(x)$. We let the reader see
how to pass from one formulation to another. Our formulation will be
more convenient for the sequel of the paper.

\begin{prop}
\label{prop:HJBclassique2} Under the assumptions of Proposition
\ref{prop:HJBclassique}, for all $x \in \partial [0,1]^d \setminus
\{0\}$,
\begin{equation}
\label{eq:dHJBclass1} \forall e \in {\mathbf e}_+(x) \cup {\mathbf
e}_0(x), \ \langle \nabla_{\alpha} H\bigl(x,\nabla V(x) \bigr),e
\rangle \geq 0 \ {\rm and} \ \langle \nabla_{\alpha}
H\bigl(x,\nabla_+ V(x) \bigr),e \rangle \geq 0,
\end{equation}
and,
\begin{equation}
\label{eq:dHJBclass2} \forall e \in {\mathbf e}_-(x), \ \langle
\nabla_{\alpha} H\bigl(x,\nabla V(x) \bigr),e \rangle \leq 0 \  {\rm
and } \   \langle \nabla_{\alpha} H\bigl(x,\nabla_+ V(x) \bigr),e
\rangle \geq 0.
\end{equation}
\end{prop}

\noindent $\Box$ Proof. We fix $x \in \partial [0,1]^d \setminus
\{0\}$. We start by proving that $\langle \nabla_{\alpha}
H(x,\nabla_+ V(x)),e \rangle \geq 0$ for all $e \in {\mathbf e}(x)$.

We know that $p = \nabla_+ V(x)$ satisfies $p_T = \nabla_T V(x)$, $\langle p,n \rangle \geq
\langle \nabla V(x),n \rangle$ and $\langle p,n \rangle \geq 0$ for
all $n \in {\mathcal N}(x)$. For $e \in {\mathbf e}(x)$ and $\lambda
>0$, the same is true when replacing $\nabla_+ V(x)$ by $\nabla_+ V(x)
+ \lambda e$. (Indeed, $\langle e,n\rangle \geq 0$.) According to
the discussion led in the proof of Proposition
\ref{prop:HJBclassique}, we can find an admissible $\theta$ such
that $\nabla \theta(x)=\nabla_+ V(x) + \lambda e$ in
\eqref{eq:eqHJBound}. We deduce that $H(x,\nabla_+ V(x) + \lambda e)
\geq 0$. Since $H(x,\nabla_+ V(x))=0$, we obtain $\langle
\nabla_{\alpha} H(x,\nabla_+ V(x)),e \rangle \geq 0$.

As a by-product, the first inequality in \eqref{eq:dHJBclass1} is
true when $\langle \nabla V(x),n \rangle \geq 0$ for all $n \in {\mathcal N}(x)$, i.e. 
when ${\mathbf e}_-(x)$ is empty. In this case, $\nabla V(x) =
\nabla_+ V(x)$.

We now prove the first inequality in \eqref{eq:dHJBclass2} when
$\langle \nabla V(x),n\rangle < 0$ for all $n \in {\mathcal N}(x)$,
i.e. ${\mathbf e}_+(x) = {\mathbf e}_0(x)=\emptyset$. Then, for $e
\in {\mathbf e}_-(x)$ and $0 < \lambda < |\langle \nabla
V(x),e\rangle |$, $p=\nabla V(x) + \lambda e$ satisfies $\langle
\nabla V(x), n \rangle \leq \langle p,n \rangle <0$ for all $n \in
{\mathcal N}(x)$.  By \eqref{eq:eqHJBound}, we deduce $H(x,\nabla
V(x) + \lambda e) \leq 0$. Since $H(x,\nabla V(x))=0$, we obtain
$\langle \nabla_{\alpha} H(x,\nabla V(x)),e \rangle \leq 0$.

We finally prove the first inequalities in \eqref{eq:dHJBclass1} and
\eqref{eq:dHJBclass2} without the assumptions ${\mathbf e}_-(x) =
\emptyset$ or ${\mathbf e}_+(x)={\mathbf e}_0(x)=\emptyset$. For $e
\in {\mathbf e}(x)$ and $\varepsilon >0$, we set $y_{\varepsilon} =
x - \varepsilon \sum_{e' \in {\mathbf e}(x),e' \not = e} e'$. For
$\varepsilon >0$ small enough, ${\mathcal N}(y_{\varepsilon})= e$.
If $e \in {\mathbf e}_-(x)$, then $\langle \nabla
V(y_{\varepsilon}),e \rangle <0$ for $\varepsilon$ small enough. By
the above analysis, $\langle \nabla_{\alpha}
H(y_{\varepsilon},\nabla V(y_{\varepsilon})),e \rangle \leq 0$.
Letting $\varepsilon$ tend to zero, we deduce that $\langle
\nabla_{\alpha} H(x,\nabla V(x)),e \rangle \leq 0$. If $e \in
{\mathbf e}_+(x) \cup {\mathbf e}_0(x)$, we know, by the above
analysis, that $\langle \nabla_{\alpha} H(y_{\varepsilon},\nabla_+
V(y_{\varepsilon})),e \rangle \geq 0$. As $\varepsilon$ tends to 0,
$\nabla_+ V(y_{\varepsilon}) \rightarrow \nabla V(x)$. (To prove it, it is sufficient to check that
$\langle \nabla_+ V(y_{\varepsilon}),e \rangle \rightarrow \langle \nabla V(x),e \rangle = 0$. 
Since $\langle \nabla_+ V(y_{\varepsilon}),e \rangle = \langle \nabla V(y_{\varepsilon}),e \rangle
{\mathbf 1}_{\{\langle \nabla V(y_{\varepsilon}),e \rangle \geq 0\}}$, this is true.) In the limit, we obtain $\langle \nabla_{\alpha}
H(x,\nabla  V(x)),e \rangle \geq 0$. \qed

{\bf Uniqueness of the Solution.} The above results provide the
typical form, both in the viscosity and in the classical senses, of
the Hamilton-Jacobi equation satisfied by the quasi-potential. A practical
question is to identify the quasi-potential with a
known solution of the Hamilton-Jacobi equation.

Generally speaking, we are not able to prove that there is a unique
continuous viscosity solution $u$ satisfying both $u(0)=0$ and
\eqref{eq:eqHJB} and \eqref{eq:eqHJBound}. By adapting the
techniques exposed in \cite{lions:90}, we can only prove, under
additional assumptions on $H$, that there exists at most one
bilateral subsolution $u$ to the Hamilton-Jacobi equation inside
$(0,1)^d$ satisfying at the same time $u(0)=0$, \eqref{eq:eqHJB} and
\eqref{eq:eqHJBound}. (Recall that $u$ is a bilateral subsolution at $x \in (0,1)^d$ if $H(x,\nabla \theta(x))=0$ for any
continuously differentiable $\theta$ such that $u-\theta$ has a local maximum at $x$.)
We won't perform the proof in the paper since
we do not whether the quasi-potential is a bilateral subsolution of
the Hamilton-Jacobi equation inside $(0,1)^d$.

Indeed, as already explained, the only thing we know is: if there
exists an optimal path from 0 to $x \in (0,1)^d$, then the
quasi-potential is a bilateral subsolution of the Hamilton-Jacobi
equation. Proving the existence of optimal paths for general
quasi-potentials may be very difficult. (See e.g. \cite[\S 2.5,
Chapter III]{bardi:capuzzo}.)

Anyhow, if the quasi-potential is assumed to continuously
differentiable, finding optimal paths may be easier. (See e.g.
\cite{FV} for a general result concerning the non-reflected case.)
For this reason, we feel simpler to provide a uniqueness result to
the Hamilton-Jacobi equation, but just for classical solutions. More
specifically,  we provide below a uniqueness result in which we both
identify the quasi-potential with a known classical solution of the
Hamilton-Jacobi equation and build optimal paths as solutions of a
suitable backward reflected differential equation.

We start with the necessary form of the optimal paths, if exist. To
this end, we extend $\nabla_+ V$ to the whole $[0,1]^d \setminus
\{0\}$ by setting $\nabla_+ V(y) = \nabla V(y)$ if $y \in (0,1)^d$.

\begin{prop}
\label{prop:optimnecessary} Under {\bf (A.1--3)}, assume
that the quasi-potential $V$ is
continuously differentiable on $[0,1]^d \setminus \{0\}$. Let 
 $x \in [0,1]^d \setminus \{0\}$ and $(\varphi_t)_{t
\leq 0}$ be a path satisfying $\varphi_0=x$ and $\lim_{t \rightarrow -
\infty} \varphi_t=0$ and achieving the infimum in the
definition of $V(x)$. Then,  $(\varphi_t)_{t  \leq 0}$ 
is absolutely continuous and verifies the  backward reflected differential 
equation
\begin{equation}
\label{eq:optimal} \dot{\varphi}_t = \nabla_{\alpha} H \bigl(
\varphi_t,\nabla_+ V(\varphi_t) \bigr) - \dot{k}_t \ {\rm a.e. \ on
\ the \ set} \ \{t \leq 0 : \, \varphi_t \not = 0\},
\end{equation}
$k$ being as in \eqref{Skorokhod}, i.e. $\dot{k}_t \in {\mathcal N}(\varphi_t)$ if 
$\varphi_t \in \partial [0,1]^d$ and $\dot{k}_t=0$ otherwise, and satisfying the compatibility
condition
\begin{equation}
\label{eq:compa} \langle \dot{k}_t,\nabla_+ V(\varphi_t) \rangle
 =0 \ {\rm a.e. \ on \
the \ set} \ \{t \leq 0 : \, \varphi_t \not = 0\}.
\end{equation}
(We emphasize that $\{t \leq 0: \varphi_t \not = 0\}$ is an interval. Indeed, if $\varphi_t =0$ for some $t \leq 0$, then 
$\varphi_s = 0$ for $s \leq t$.)
\end{prop}

\noindent $\Box$ Proof. We admit for the moment the following
\begin{lem}
\label{lem:optim} For every compact subset $\kappa \subset [0,1]^d
\setminus \{0\}$, there exists a constant $c_{\kappa} >0$ such that
for all $y \in \kappa$ and $v \in \R^d$, $|v| \leq 1$,
\begin{equation}
\label{eq:uniquenessoptimal} L(y,v) \geq \langle v,\nabla_+ V(y)
\rangle + c_{\kappa} \bigl|\nabla_{\alpha} H \bigl(y,\nabla_+ V(y)
\bigr) - v \bigr|^2.
\end{equation}
\end{lem}

We then consider a path $(\phi_t)_{t \leq 0}$ with $\phi_0 = x \not =0$,
$\lim_{t \rightarrow - \infty} \phi_t =0$ and $J_{-\infty,0}(\phi) <
+ \infty $ (so that, without loss of generality, $|\dot{\phi}_t|
\leq 1$ for a.e. $t \leq 0$). By \eqref{eq:Lref}, we can find a
measurable mapping $t \in (-\infty,0] \mapsto (\beta_t,n_t) \in
{\mathbb R}_+ \times {\mathcal N}(\phi_t)$ such that for a.e. $t
\leq 0$
\begin{equation}
\label{eq:optim100} L^{\rm ref}(\phi_t,\dot{\phi}_t) =
L(\phi_t,\dot{\phi}_t + \beta_t n_t).
\end{equation}
(In the above formula, $\beta_t=0$ if $\phi_t \in (0,1)^d$ or
$\phi_t \in \partial [0,1]^d$ and $\langle \dot{\phi}_t,n \rangle
<0$ for all $n\in{\mathcal N}(\phi_t)$. We refer to the proof of
Theorem \ref{expr_act_func} for the measurability property.
We also note that $|\dot{\phi}_t + \beta_t n_t | \leq 1$ for a.e. $t \leq 0$ since 
$L(\phi_t,\dot{\phi}_t + \beta_t n_t) < + \infty$.)

For a given compact subset $\kappa \subset [0,1]^d \setminus \{0\}$
containing $x$, we set $T_{\kappa} = \inf \{T \geq 0 : \,
\varphi_{-T} \not \in \kappa \}$. Lemma \ref{lem:optim} and
\eqref{eq:optim100} yield for a.e. $t \in [-T_{\kappa},0]$
\begin{equation*}
\label{eq:path1} L^{\rm ref}(\phi_t,\dot{\phi}_t) \geq \big\langle
\dot{\phi}_t + \beta_t n_t,\nabla_+ V(\phi_t) \big\rangle
 + c_{\kappa}
\bigl|\nabla_{\alpha} H\bigl(\phi_t,\nabla_+ V(\phi_t)\bigr) -
(\dot{\phi}_t + \beta_t n_t) \bigr|^2.
\end{equation*}
We let the reader check that, for $i\in \{1,\cdots,d\}$, the
Lebesgue measure of the set $\{t \leq 0: \, (\phi_t)_i \in\{0,1\},
(\dot{\phi}_t)_i \not=0\}$ is zero. 
(Indeed, the path $\phi$ is a.e. differentiable.)
Hence, $\langle
\dot{\phi}_t,\nabla_+ V(\phi_t) \rangle = [d/dt](V(\phi_t))$ for
a.e. $t \leq 0$. We deduce that for a.e. $t \in [-T_{\kappa},0]$
\begin{equation*}
\label{eq:path2} L^{\rm ref}(\phi_t,\dot{\phi}_t) \geq
[d/dt]\bigl(V(\phi_t) \bigr) + \big\langle \beta_t n_t,\nabla_+
V(\phi_t) \big\rangle + c_{\kappa} \bigl|\nabla_{\alpha} H
\bigl(\phi_t,\nabla_+ V(\phi_t) \bigr) - (\dot{\phi}_t + \beta_t
n_t) \bigr|^2.
\end{equation*}
We deduce that $\phi$ satisfies
\begin{equation}
\label{eq:path3}
\begin{split}
J_{-\infty,0}(\phi) &\geq V(\phi_{-T_{\kappa}}) +
\int_{-T_{\kappa}}^0 L^{\rm ref}(\phi_t,\dot{\phi}_t) dt
\\
&\geq V(\phi_0) + \int_{-T_{\kappa}}^0 \langle \beta_t n_t, \nabla_+
V(\phi_t) \rangle dt
\\
&\hspace{15pt} + c_{\kappa} \int_{- T_{\kappa}}^0
\bigl|\nabla_{\alpha} H \bigl(\phi_t,\nabla_+ V(\phi_t) \bigr) -
(\dot{\phi}_t + \beta_t n_t) \bigr|^2 dt.
\end{split}
\end{equation}
Noting that  $\langle \beta_t n_t,\nabla_+ V(\phi_t) \rangle \geq 0$
for all $t \leq 0$, we complete the proof. \qed

\noindent $\Box$ Proof of Lemma \ref{lem:optim}. For $y \in \kappa$,
$v \in {\mathbb R}^d$, $|v| \leq 1$, and $\varepsilon \in (-1,1)$,
\begin{equation*}
\begin{split}
L(y,v) &\geq \big\langle v,\nabla_+ V(y) - \varepsilon
\bigl[\nabla_{\alpha} H \bigl(y,\nabla_+ V(y) \bigr) - v \bigr]
\bigr\rangle
\\
&\hspace{5pt} - H \Bigl(y,\nabla_+ V(y) - \varepsilon
\bigl[\nabla_{\alpha} H \bigl(y,\nabla_+ V(y) \bigr) - v \bigr]
\Bigr).
\end{split}
\end{equation*}
By Proposition \ref{prop:HJBclassique}, we know that $H(y,\nabla_+
V(y)) = 0$. Applying Taylor's formula, in zero, to the function
\begin{equation*}
\varepsilon \in (-1,1) \mapsto H \Bigl(y,\nabla_+ V(y) - \varepsilon
\bigl[ \nabla_{\alpha} H \bigl(y,\nabla_+ V(y) \bigr)- v \bigr]
\Bigr),
\end{equation*}
we obtain
\begin{equation*}
\begin{split}
L(y,v) &\geq \big\langle v,\nabla_+ V(y) - \varepsilon
\bigl[\nabla_{\alpha} H\bigl(y,\nabla_+ V(y) \bigr) -v \bigr]
\big\rangle
\\
&\hspace{15pt} + \bigl\langle \nabla_{\alpha} H \bigl(y,\nabla_+
V(y) \bigr),\varepsilon \bigl[\nabla_{\alpha} H\bigl(y,\nabla_+ V(y)
\bigr) -
 v \bigr] \bigr\rangle
\\
&\hspace{15pt}
 - (C/2) \varepsilon^2 \bigl|\nabla_{\alpha} H \bigl(y,\nabla_+
V(y) \bigr) - v \bigr|^2,
\end{split}
\end{equation*}
with $C=\sup\{| \nabla_{\alpha,\alpha}^2H(z,\nabla_+ V(z) - \eta
[\nabla_{\alpha} H(z,\nabla_+ V(z)) - v])|; \, z \in \kappa, \, v \in
{\mathbb R}^d, \, |v| \leq 1, \, \eta \in [-1,1]\}$. By the
regularity of $H$ and $\nabla V$, the constant $C$ is finite. Hence,
\begin{equation*}
L(y,v) \geq \langle v,\nabla_+ V(y) \rangle + \bigl(\varepsilon -
(C/2) \varepsilon^2 \bigr) \bigl|\nabla_{\alpha} H \bigl(y,\nabla_+
V(y) \bigr) - v \bigr|^2.
\end{equation*}
Without loss of generality, we can assume that $C>1$ and choose
$\varepsilon = 1/C$ in the above formula. This completes the proof.
\qed

In light of Proposition \ref{prop:optimnecessary}, we understand that the boundary conditions in
Proposition \ref{prop:HJBclassique2} describe the shape of the optimal paths (if exist) at the boundary.

In what follows, we explain more specifically what happens in dimension two. For
example, we consider $x$ on the boundary with $x_1=0$ and $x_2 \in (0,1)$. In this case
${\mathbf e}(x)=\{-e_1\}$.

If $[\partial V/\partial x_1](x) > 0$,
then $\langle \nabla V(x),-e_1 \rangle <0$ and 
 $-e_1 \in {\mathbf e}_-(x)$.
By Proposition
\ref{prop:HJBclassique2}, we know that 
 $\langle \nabla_{\alpha} H(x,\nabla V(x)),-e_1
\rangle \leq 0$, i.e. 
 $\langle \nabla_{\alpha} H(x,\nabla V(x)),e_1
\rangle \geq 0$. Assume to simplify that 
 $\langle \nabla_{\alpha} H(x,\nabla V(x)),e_1
\rangle > 0$.
 By continuity, $\langle \nabla_{\alpha} H(y,\nabla
V(y)),e_1 \rangle > 0$ for $y$ in a neighborhood of $x$. If there exists an
optimal path $(\varphi_t)_{t \leq 0}$ reaching $x$ at $t=0$, we understand
from \eqref{eq:optimal}
 that
$(\varphi_t)_{t \leq 0}$ has to hit the boundary before reaching $x$. (Otherwise,
there exists $\varepsilon >0$ such that $(\varphi_t)_1 >0$ for $t \in [-\varepsilon,0)$,
so that $(\dot{\varphi}_t)_1 = \langle \nabla_{\alpha} H(\varphi_t,\nabla V(\varphi_t)),e_1 \rangle > 0$, and,
the path cannot reach $x$.)
 This is illustrated by Figure \ref{figu:optimal1} below.
\begin{figure}[htb]
\begin{center}
\psfrag{x}{$x$}
\includegraphics[
width=0.06\textwidth,angle=0] {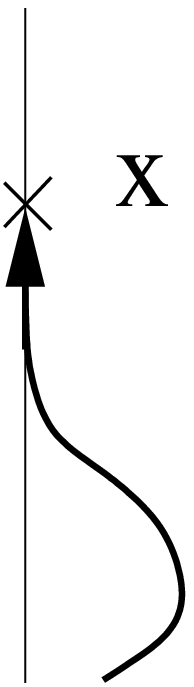} \caption{Typical optimal
path: $[\partial V/\partial x_1](x)>0$ and  $\langle \nabla_{\alpha} H(x,\nabla V(x)),e_1
\rangle > 0$}
\label{figu:optimal1}
\end{center}
\end{figure}

Similarly, if $[\partial V/\partial x_1](x) <0$, i.e. $-e_1 \in {\mathbf e}_+(x)$, we know from
Proposition
\ref{prop:HJBclassique2} that $\langle \nabla_{\alpha} H(x,\nabla V(x)),e_1
\rangle \leq 0$. We assume to simplify that
$\langle \nabla_{\alpha} H(x,\nabla V(x)),e_1
\rangle < 0$. For $y$ in a neighborhood of $x$, $\langle \nabla_{\alpha} H(y,\nabla V(y)),e_1
\rangle < 0$. Since $[\partial V/\partial x_1](x)<0$, we also have
$[\partial V/\partial x_1](y)<0$ and thus
$\nabla_+ V(y) = \nabla V(y)$ for $y$ close to $x$. Thus,
$\langle \nabla_{\alpha} H(y,\nabla_+ V(y)),e_1
\rangle < 0$ for $y$ in a neighborhood of $x$.
Then, the first coordinate of
$\varphi_t$, i.e. $(\varphi_t)_1$, is non-increasing as $t$ grows up to 0.
In particular, if $\varphi_{-\varepsilon}=0$ for some small $\varepsilon >0$, the path remains on the boundary
from time $-\varepsilon$ to time $0$. In such a case, $\dot{k}_t =  \langle \nabla_{\alpha} H(\varphi_t,\nabla V(\varphi_t)),e_1
\rangle e_1$ for a.e. $t \in [-\varepsilon,0]$ so that $\langle \nabla_+ V(\varphi_t),\dot{k}_t \rangle = [\partial V/\partial x_1](\varphi_t)
\langle \nabla_{\alpha} H(\varphi_t,\nabla V(\varphi_t)),e_1
\rangle >0$. This violates the compatibility condition \eqref{eq:compa}.
We deduce that the optimal path cannot hit the boundary in a small neighborhood of
$x$ before reaching $x$.
This is illustrated by Figure \ref{figu:optimal2} below.
\begin{figure}[htb]
\begin{center}
\psfrag{x}{$x$}
\includegraphics[
width=0.1\textwidth,angle=0] {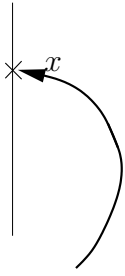} \caption{Typical optimal
path: $[\partial V/\partial x_1](x)<0$ and $\langle \nabla_{\alpha} H(x,\nabla V(x)),e_1 \rangle <0$}
\label{figu:optimal2}
\end{center}
\end{figure}

The case where $[\partial V/\partial x_1](x) =0$ leads to too many different possibilities to make a general
comment. (Anyhow, an example is provided in the next section.)

Proposition \ref{prop:optimnecessary} shows that, if optimal paths
exist, the reflected differential equation \eqref{eq:optimal} is
solvable. We emphasize that \eqref{eq:optimal} is not a reflected
differential equation of standard type since the boundary condition
is given by the terminal value of the trajectory. In particular,
solving \eqref{eq:optimal} is more intricate than solving a standard
Skorohod problem. As shown below,  the boundary conditions
\eqref{eq:dHJBclass1} and \eqref{eq:dHJBclass2} play a crucial role
in the solvability of the equation \eqref{eq:optimal}.

\begin{prop}
\label{prop:optimsuff} Assume {\bf (A.1--3)} and 
that there exists a 
function $W \in {\mathcal C}([0,1]^d,{\mathbb R})$, continuously
differentiable on $[0,1]^d \setminus \{0\}$, such that, for all $x
\in (0,1)^d$, $ H(x,\nabla W(x)) = 0$, for all $x \in \partial
[0,1]^d \setminus \{0\}$, $H( x, \nabla_+ W(x)) =0$, and, for all $e
\in {\bf e}(x)$, $\langle \nabla_{\alpha} H (x,\nabla_+ W(x) ),e
\rangle \geq 0$.


In addition, assume that, for all $x \in \partial [0,1]^d \setminus \{0\}$
and for all $e \in {\bf e}_0^W(x)= \{e \in {\mathbf e}(x), \ \langle
\nabla_{\alpha} H(x,\nabla_+ W(x)),e \rangle = 0\}$, there exists a
neighborhood $U$ of $x$ such that the sign of $\langle \nabla W,e
\rangle$ is constant on the intersection of $U$ with the face
orthogonal to $e$, i.e. either
\begin{equation}
\label{eq:optimsuff01} \forall y \in U \cap \partial [0,1]^d, \ e
\in {\mathcal N}(y) \Rightarrow \langle \nabla W(y),e \rangle \leq
0,
\end{equation}
or,
\begin{equation}
\label{eq:optimsuff02} \forall y \in U \cap \partial [0,1]^d, \ e
\in {\mathcal N}(y) \Rightarrow \langle \nabla W(y),e \rangle \geq
0.
\end{equation}
Then, for any $x \in [0,1]^d \setminus \{0\}$, there exist  an
absolutely continuous path $(\varphi_t)_{t \leq 0}$ and a real
$i(\varphi) \in (0,+\infty]$ such that $\varphi_t=0$ for all $t \leq
-i(\varphi)$ if $i(\varphi) < + \infty$, $\varphi_0=x$ and
\begin{equation}
\label{eq:optimsuff1}  \dot{\varphi}_t = \nabla_{\alpha} H \bigl(
\varphi_t, \nabla_+ W(\varphi_t) \bigr)  - \dot{k}_t, \ {\rm for \
a.e.} \ t \in (-i(\varphi),0],
\end{equation}
$k$ being as in \eqref{Skorokhod}, i.e. $\dot{k}_t \in {\mathcal N}(\varphi_t)$ if 
$\varphi_t \in \partial [0,1]^d$ and $\dot{k}_t=0$ otherwise, and satisfying the compatibility
condition
\begin{equation}
\label{eq:optimsuff2} \langle \dot{k}_t,\nabla_+ W(\varphi_t)
\rangle = 0 \ {\rm for \ a.e.} \ t \in (-i(\varphi) ,0].
\end{equation}
(Above, $\nabla_+ W(x)=\nabla W(x)$ for $x \in (0,1)^d$.)
\end{prop}

The additional conditions \eqref{eq:optimsuff01} and \eqref{eq:optimsuff02} permit to avoid degenerate situations in which the sign of
$\langle \nabla W,e \rangle$ changes at $x$ for some $e \in {\mathcal N}(y)$. Having in mind Figures \ref{figu:optimal1} and \ref{figu:optimal2},
this permits to determine, {\it a priori}, the shape of the optimal paths reaching $x$.

We emphasize that no assumption is necessary on the sign of 
$\langle \nabla_{\alpha} H(x,\nabla W(x)),e \rangle$, $e \in {\mathbf e}(x)$.
(In fact, using the convexity of $H$, we could prove that 
all the inequalities in \eqref{eq:dHJBclass1} and \eqref{eq:dHJBclass2}
hold for $W$ under the assumptions of Proposition \ref{prop:optimsuff}.)
\\

\noindent $\Box$ Proof. It is sufficient to prove that, for all $x \in
[0,1]^d \setminus \{0\}$, there exist a real $\varepsilon >0$ and an
absolutely continuous path $(\varphi_t)_{- \varepsilon \leq t \leq
0}$ such that $\varphi_0 =x $ and \eqref{eq:optimsuff1} and
\eqref{eq:optimsuff2} hold on $[-\varepsilon,0]$. (By concatenating
the local solutions, we obtain a global solution. When the resulting
path hits the origin, the concatenation procedure stops. In this
case, $i(\varphi)$ is finite. If the path doesn't hit the origin, $i(\varphi)$ is
infinite.
In the next theorem, we will prove under additional assumptions on $W$ that the path tends to 0 as $t$ tends to $-\infty$ if $i(\varphi)=+ \infty$.)

If $x \in (0,1)^d$, the proof is trivial. (It is sufficient to
solve, locally, the backward differential equation
\begin{equation*}
\dot{\varphi}_t = \nabla_{\alpha} H \bigl( \varphi_t, \nabla
W(\varphi_t) \bigr), \ t \leq 0,
\end{equation*}
with the boundary condition $\varphi_0=x$. 
Since $\nabla_{\alpha}H$ is bounded by 1 and $\nabla
W$ is continuous, this is possible.)

If $x \in \partial [0,1]^d \setminus \{0\}$, the idea still consists
in solving a backward differential equation, without reflection,
but along a face of the hypercube.

We first specify the choice of the face. By \eqref{eq:optimsuff01}
and \eqref{eq:optimsuff02}, there exists a neighborhood $U$ of $x$
such that, for all $e \in {\mathbf e}^W(x)$, the
sign of $\langle \nabla W,e \rangle$ is constant on the intersection
of $U$ with the face orthogonal to $e$. (If $e \in {\mathbf
e}(x) \setminus {\mathbf e}_0^W(x)$, this is trivial by continuity of
$\nabla W$. If $e \in {\mathbf e}_0^W(x)$, this follows from \eqref{eq:optimsuff01}
and \eqref{eq:optimsuff02}.) We then consider the (largest) face ${\it f}$
containing $x$ and orthogonal to ${\mathbf e}_f(x)$, with
\begin{equation*}
{\mathbf e}_f(x) = \bigl\{e \in {\mathbf e}(x) :
\forall y \in U \cap \partial [0,1]^d, \ e \in {\mathcal N}(y)
\Rightarrow \langle \nabla W(y),e \rangle \leq 0\bigr\}.
\end{equation*}
We denote by $\delta$ the dimension of $f$. We can find a
subset $F \subset \{1,\dots,d\}$, the cardinal of $F$ being
equal to $\delta$, such that the family $(e_j)_{j \in F}$ is a basis of the
plane generated by $f$. We then consider the system of differential equations
\begin{equation}
\label{eq:sys1}
\begin{split}
&(\dot{\varphi}_t)_j =  \langle \nabla_{\alpha} H\bigl(
\varphi_t,\nabla_+ W(\varphi_t) \bigr), e_j \rangle,  \ t \leq 0,\ j \in F,
\\
&(\dot{\varphi}_t)_j =  0, \ t \leq 0, \ j \in \{1,\dots,d\} \setminus F,
\end{split}
\end{equation}
with the boundary condition $\varphi_0=x$.

A priori, this problem isn't well-posed, even in a small time duration. Indeed,
$\varphi_t$ may
leave the hypercube in a zero time
so that $\nabla_{\alpha} H(\varphi_t,\nabla_+ W(\varphi_t))$ may not be defined.
(Recall that $x_i$ may be 0 or 1 for some $i \in F$.)
To obtain a well-posed problem,
we consider the following version
\begin{equation}
\label{eq:sys2}
\begin{split}
&(\dot{\varphi}_t)_j =  \langle \nabla_{\alpha} H\bigl(\Pi[
\varphi_t],\nabla_+ W(\Pi[\varphi_t]) \bigr), e_j \rangle, \ t \leq 0,\ j \in F,
\\
&(\dot{\varphi}_t)_j =  0, \ t \leq 0,\ j \in \{1,\dots,d\} \setminus F,
\end{split}
\end{equation}
with the same boundary condition as above, where $\Pi$ denotes the projection
on the hypercube. In the above system, either $e_j$ or $-e_j$ belongs to ${\mathbf
e}_f(x)$ for $j \in \{1,\dots,d\} \setminus F$. (That is $x_j = 0$ or $1$.)
Since $(\dot{\varphi}_t)_j=0$ for such
$j$'s, ${\mathbf e}_f(x) \subset {\mathcal N}(\Pi[\varphi_t])$.
For $t$ close to zero,
$\Pi[\varphi_t] \in U$, so that 
$\langle \nabla W(\Pi[\varphi_t]),e \rangle \leq 0$ for $e \in {\mathbf e}_f(x)$. As a
by-product, $\langle \nabla_+ W(\Pi[\varphi_t]),e_j \rangle$ is equal to
0 for $j \not \in F$. For $j \in F$, either $(\Pi[\varphi_t])_j \in (0,1)$ 
or $(\Pi[\varphi_t])_j \in \{0,1\}$. In the first case,
$\langle \nabla_+ W(\Pi[\varphi_t]),e_j \rangle$ is equal to
$\langle \nabla W(\Pi[\varphi_t]),e_j \rangle$. In the second case,
either $e_j$ or $-e_j$ is a normal vector at $\Pi[\varphi_t]$ and belongs to 
${\mathbf e}(x) \setminus {\mathbf
e}_f(x)$, so that $\langle \nabla_+ W(\Pi[\varphi_t]),e_j \rangle$ is still equal to
$\langle \nabla W(\Pi[\varphi_t]),e_j \rangle$.
We deduce 
that $ \nabla_+ W(\Pi[\varphi_t])$ may be expressed as $\sum_{j \in F} \langle \nabla W(\Pi[\varphi_t]),e_j \rangle e_j$ in the above system.
Thus, the coefficients of the system are continuous in the neighborhood of 
the boundary condition,
so that the problem admits a solution on
some interval $[-\varepsilon,0]$, $\varepsilon >0$.

We now show that we can get rid of $\Pi$, at least for $\varepsilon$ small enough. To do so, it is enough
to prove that $\varphi_t$ belongs to $[0,1]^d$, or, equivalently,
that $(\varphi_t)_j \in [0,1]$ for $j \in
\{1,\dots,d\}$. For $j \not \in F$, this is obvious since $(\dot{\varphi}_t)_j = 0$. We thus assume $j \in F$.
 If $(\Pi[\varphi_t])_j >0$, then
$(\varphi_t)_j =(\Pi[\varphi_t])_j >0$. If
$(\Pi[\varphi_t])_j=0$, then $-e_j$ is a normal vector to the
hypercube at $\Pi[\varphi_t]$. By the boundary conditions satisfied by $W$, this implies
$\langle \nabla_{\alpha} H (\Pi[\varphi_t],\nabla_+
W(\Pi[\varphi_t])),-e_j\rangle \geq 0$. In this case, $
(\dot{\varphi}_t)_j \leq 0$. As $t$ decreases on
$[-\varepsilon,0]$, $(\varphi_t)_j$ cannot go below 0. Similarly, it cannot go beyond 1.
We deduce that, for $\varepsilon$ small enough, \eqref{eq:sys1} holds true.

We finally prove that \eqref{eq:optimsuff1} holds on $[-\varepsilon,0]$. 
%
We can always write
\begin{equation*}
\begin{split}
\dot{\varphi}_t &= \nabla_{\alpha} H\bigl( \varphi_t,\nabla_+
W(\varphi_t) \bigr)  - \sum_{j \in \{1,\dots,d\} \setminus F}
\langle \nabla_{\alpha} H\bigl( \varphi_t,\nabla_+ W(\varphi_t)
\bigr), e_j\rangle e_j
\\
&=\nabla_{\alpha} H\bigl( \varphi_t,\nabla_+
W(\varphi_t) \bigr)  - \dot{k}_t,
\end{split}
\end{equation*}
with
\begin{equation*}
\dot{k}_t = \sum_{j \in \{1,\dots,d\} \setminus F} \langle
\nabla_{\alpha} H\bigl( \varphi_t,\nabla_+ W(\varphi_t) \bigr),
e_j\rangle e_j =  \sum_{e \in {\mathbf e}_f(x)} \langle
\nabla_{\alpha} H\bigl( \varphi_t,\nabla_+ W(\varphi_t) \bigr), e
\rangle e.
\end{equation*}
Since $e_f(x) \subset {\mathcal N}(\varphi_t)$, 
$\langle \nabla_{\alpha} H(\varphi_t,\nabla_+ W(\varphi_t)),e \rangle \geq 0$ for
all $e \in {\mathbf e}_f(x)$. We deduce that $(k_t)_{- \varepsilon\leq t \leq 0}$ (with
$k_0=0$) satisfies \eqref{Skorokhod}. The
compatibility condition is obviously true. \qed

We are now in position to state an identification property for the
quasi-potential.
\begin{thm}
\label{prop:identification} In addition to {\bf (A.1--4)}, assume
that, for all $x \in [0,1]^d \setminus \{0\}$, $\langle x,\bar{f}(x)
\rangle <0$. Assume also that there exists a function $W$ satisfying
the conditions of Proposition \ref{prop:optimsuff} such that
$W(0)=0$. Then $W$ is equal to the quasi-potential  and the infimum
in the quasi-potential is attained at $(\varphi_t)_{t \leq 0}$ given
by Proposition \ref{prop:optimsuff}. (We show below that such a
path satisfies $\lim_{t \rightarrow - \infty}\varphi_t = 0$).
\end{thm}

In the proof, we use the following lemma (the proof is given in
Appendix, see Subsection \ref{sec:convex}). \vspace{5pt}

\begin{lem}
\label{lem:convex} Under {\bf (A.1--3)}, for any $x \in [0,1]^d$,
the mapping $\alpha \in \R^d \mapsto H(x,\alpha)$ is strictly convex
at 0, i.e. the matrix $([\partial^2 H/\partial \alpha_i \partial
\alpha_j](x,0))_{i,j \in \{1,\dots,d\}}$ is positive definite.
\end{lem}

\noindent $\Box$ Proof of Theorem \ref{prop:identification}.
We first prove that $V \geq W$. For a given $x \in [0,1]^d \setminus
\{0\}$, we can consider a path $(\psi_t)_{t\leq 0}$ from $0$ to $x$, i.e.
$\lim_{t \rightarrow - \infty} \psi_t = 0$ and $\psi_0=x$,
such that $J_{-\infty,0}(\psi) \leq V(x) + \delta$ for some $\delta
>0$. Then, $\psi$ is absolutely continuous.
For a.e. $t \leq 0$ such that $\psi_t \in (0,1)^d$, we have
\begin{equation}
\label{eq:id001}
L^{\rm ref}(\psi_t,\dot{\psi}_t) = L (\psi_t,\dot{\psi}_t)
\geq \langle \dot{\psi}_t,\nabla W(\psi_t)\rangle - H (\psi_t,\nabla
W(\psi_t)) = [d/dt](W(\psi_t)),
\end{equation}
since $W$ satisfies the Hamilton-Jacobi equation. The same holds for
$t$ satisfying $\psi_t \in
\partial [0,1]^d \setminus \{0\}$ and $\langle \dot{\psi}_t,n \rangle <0$ for all $n \in {\mathcal
N}(\psi_t)$. For $t$ satisfying
$\psi_t \in
\partial [0,1]^d \setminus \{0\}$ and $\exists n \in{\mathcal N}(\psi_t)$ such
that $\langle \dot{\psi}_t,n \rangle =0$, we claim
\begin{equation*}
\begin{split}
L^{\rm ref}(\psi_t,\dot{\psi}_t) &= \inf_{\beta >0, n \in {\mathcal
N}(\varphi_t), n \perp \dot{\psi}_t} L(\psi_t,\dot{\psi}_t + \beta
n)
\\
&\geq \inf_{\beta >0, n \in {\mathcal N}(\psi_t), n \perp
\dot{\psi}_t} \bigl[ \langle \nabla_+ W(\psi_t),\dot{\psi}_t + \beta
n \rangle - H \bigl(\psi_t,\nabla_+ W(\psi_t) \bigr) \bigr]
\\
&= \inf_{\beta >0, n \in {\mathcal N}(\psi_t), n \perp \dot{\psi}_t}
\bigl[ \langle \nabla_+ W(\psi_t),\dot{\psi}_t + \beta n \rangle
\bigr],
\end{split}
\end{equation*}
by the boundary condition of the Hamilton-Jacobi equation. By
definition of $\nabla_+ W$, we have $\langle \nabla_+ W(\psi_t),n
\rangle \geq 0$ for all $n \in {\mathcal N}(\psi_t)$. Hence,
for $t$ satisfying
$\psi_t \in
\partial [0,1]^d \setminus \{0\}$ and $\exists n \in{\mathcal N}(\psi_t)$ such
that $\langle \dot{\psi}_t,n \rangle =0$, we have
\begin{equation}
\label{eq:id002}
L^{\rm ref}(\psi_t,\dot{\psi}_t) \geq \langle \nabla_+
W(\psi_t),\dot{\psi}_t \rangle.
\end{equation}
For every $i\in \{1,\dots,d\}$, the Lebesgue measure of the set $\{t
\leq 0: \, (\psi_t)_i \in \{0,1\}, \, (\dot{\psi}_t)_i \not = 0\}$
is zero. Hence, we can replace $\langle \nabla_+
W(\psi_t),\dot{\psi}_t \rangle$ by $\langle \nabla
W(\psi_t),\dot{\psi}_t\rangle = [d/dt](W(\psi_t))$ in the above inequality. By
\eqref{eq:id001} and \eqref{eq:id002}, we have
\begin{equation*}
L^{\rm ref}(\psi_t,\dot{\psi}_t) \geq [d/dt](W(\psi_t)) \ {\rm for \ a.e.} \ t \ {\rm such \ that} \ \psi_t \not = 0.
\end{equation*}
Setting $i(\psi) = \inf\{ T \geq 0, \ \psi_{-T} =0\}$ and integrating from $-i(\psi)$ to $0$
($i(\psi)$ being possibly equal to $+\infty$), we deduce that $ V(x) + \delta
\geq W(x)$. Letting $\delta$ tend to $0$, we deduce that $V(x) \geq W(x)$.

We now prove that $V \leq W$. We consider consider the path
$(\varphi_t)_{t\leq 0}$ given by Proposition \ref{prop:optimsuff}.
Recall from
\cite[Chapter 5, (1.5)]{FV} that 
$L(y,\nabla_{\alpha} H(y,v)) = \langle v,\nabla_{\alpha} H(y,v) \rangle
- H(y,v)$ for all $(y,v) \in [0,1]^d \times \R^d$. By the Hamilton-Jacobi equation satisfied by $W$ and by the compatibility condition
\eqref{eq:optimsuff2}, we obtain, for a.e. $-i(\varphi) < t \leq 0$,
\begin{equation}
\label{eq:id1}
\begin{split}
L(\varphi_t,\dot{\varphi}_t + \dot{k}_t)  &= \langle \nabla_+
W(\varphi_t),\dot{\varphi}_t + \dot{k}_t \rangle - H
\bigl(\varphi_t,\nabla_+ W(\varphi_t) \bigr)
\\
&= \langle \nabla_+ W(\varphi_t),\dot{\varphi}_t \rangle =
[d/dt]\bigl(W(\varphi_t) \bigr),
\end{split}
\end{equation}
the last equality following from the same observation as above: for every $i\in \{1,\dots,d\}$, the Lebesgue measure of the set $\{t
\leq 0: \, (\varphi_t)_i \in \{0,1\}, \, (\dot{\varphi}_t)_i \not = 0\}$ is zero.
Hence, for any $T>0$, $T \geq i(\varphi)$,
\begin{equation*}
V(\varphi_{-T},x) \leq \int_{-T}^0 L^{\rm ref}
(\varphi_t,\dot{\varphi}_t) dt 
\leq \int_{-T}^0 L(\varphi_t,\dot{\varphi}_t + \dot{k}_t) dt 
\leq W(x) - W(\varphi_{-T}).
\end{equation*}
If $i(\varphi) < + \infty$, the proof is over by choosing $T=i(\varphi)$.
Otherwise, we have to prove that 0 is an accumulation point of the path
$(\varphi_t)_{t \leq 0}$.

Assume for a while that there exists $\varepsilon >0$ such that, for
all $t \leq 0$, $|\varphi_t| > \varepsilon$.
(In particular, $i(\varphi)= + \infty$.)
By assumption, we know
that, for all $z \in [0,1]^d \setminus \{0\}$, $\langle z,
\nabla_{\alpha} H(z,0) \rangle = \langle z, \bar{f}(z) \rangle < 0$.
(Recall that $\bar{f}(z) = \nabla_{\alpha} H(z,0)$.) By continuity of
$\nabla_{\alpha} H$, we can find a real $\eta >0$ such that
\begin{equation}
\label{borsupLref} \inf \bigl\{ \langle z, \nabla_{\alpha} H(z,v)
\rangle; \, z \in [0,1]^d, \, |z|
 \geq \varepsilon, \, v \in \R^d, \, |v| \leq \eta \bigr\}  <0.
\end{equation}
Moreover, it is plain to see that for a.e. $t \leq 0$
\begin{equation*}
\begin{split}
[d/dt] [|\varphi_t|^2] &= 2 \langle \varphi_t,\nabla_{\alpha} H
\bigl(\varphi_t,\nabla_+ W(\varphi_t) \bigr) \rangle - 2 \langle \varphi_t,\dot{k}_t \rangle
\\
&\leq 2 \langle \varphi_t,\nabla_{\alpha} H
\bigl(\varphi_t,\nabla_+ W(\varphi_t) \bigr) \rangle.
\end{split}
\end{equation*}
(Indeed, if $(\varphi_t)_i<1$, then $(\varphi_t)_i(\dot{k}_t)_i = 0$, and,
if $(\varphi_t)_i=1$, then $(\dot{k}_t)_i \geq 0$.)
By \eqref{borsupLref}, we deduce that there exists a constant $c \in
(0,1)$ such that
\begin{equation}
\label{borsupLref2} -[d/dt] \bigl[|\varphi_t|^2 \bigr] \geq c {\mathbf
1}_{\{|\nabla_+ W(\varphi_t)| \leq \eta\}}  - c^{-1} {\mathbf
1}_{\{|\nabla_+ W(\varphi_t)| > \eta \}}.
\end{equation}
By \eqref{eq:id1}, for a.e. $t \leq 0$,
\begin{equation*}
[d/dt] \bigl[W(\varphi_t)\bigr] = L\bigl(\varphi_t,\nabla_{\alpha}
H(\varphi_t,\nabla_+ W(\varphi_t)) \bigr).
\end{equation*}
By the strict convexity of $L$, for all $z \in [0,1]^d \setminus
\{0\}$, $L(z, \nabla_{\alpha} H(z,\nabla_+ W(z)))=0$ if and only if
$\nabla_{\alpha} H(z,\nabla_+ W(z))= \bar{f}(z) = \nabla_{\alpha}
 H(z,0)$. By the strict convexity of $H(z,\cdot)$ at $0$, this is equivalent to
$\nabla_+ W(z) =0$. We deduce that
 \begin{equation*}
\inf \bigl\{ L(z,\nabla_{\alpha} H(z,\nabla_+ W(z))); \ z \in
[0,1]^d, \ | z| \geq \varepsilon, \ |\nabla_+ W(z)| \geq \eta
\bigr\} >0,
\end{equation*}
if not empty (i.e. $\exists z \in [0,1]^d, \ |z| \geq \varepsilon, \
|\nabla_+ W(z)| \geq \eta$). Up to a modification of $c$, we have
\begin{equation}
\label{borsup:2000} [d/dt] \bigl[W(\varphi_t)\bigr] \geq c {\mathbf
1}_{\{|\nabla_+ W(\varphi_t)| > \eta\}} .
\end{equation}
We deduce that $|\{t \leq 0: \, |\nabla_+ W (\varphi_t)| > \eta\}| <
+ \infty$. Hence, $|\{t \leq 0: \, |\nabla_+ W(\varphi_t)| \leq
\eta\}| = + \infty$. By \eqref{borsupLref2}, there is a
contradiction. We deduce that $0$ is an accumulation point of
$(\varphi_t)_{t \leq 0}$. Hence, $W(x) \geq V(x)$ so that
$W(x)=V(x)$.

Actually, we can prove that $\lim_{t \rightarrow - \infty} \varphi_t
= 0$. Indeed, by \eqref{eq:id1}, $(W(\varphi_t))_{t \leq 0}$ is
nondecreasing (and bounded). We deduce that $\lim_{t \rightarrow -
\infty} W(\varphi_t)=0$ since $0$ is an accumulation point of the
sequence $(\varphi_t)_{t \leq 0}$. Hence, every accumulation point $a$
of the sequence $(\varphi_t)_{t \leq 0}$ satisfies $W(a)=0$. Assume
that there exists another accumulation point $a \not = 0$. Since $0$ is an accumulation point, we can find
two decreasing sequences $(t_n)_{n \geq 0}$ and $(s_n)_{n \geq 0}$,
converging to $- \infty$, such that $t_{n+1} < s_{n+1} < t_n < s_n$
for all $n \geq 0$, $|\varphi_{t_n}|=|a|/2$  for all $n \geq 0$,
$|\varphi_r| \geq |a|/2$ for all $r \in [t_n,s_n]$ and $n \geq 0$,
and $|\varphi_{s_n}-a| \rightarrow 0$. By \eqref{borsupLref2} and
\eqref{borsup:2000}, we can find some constant $C>0$ (depending on
$a$) such that $t \mapsto -|\varphi_t|^2 + C W(\varphi_t)$ is
nondecreasing on each $[t_n,s_n]$, $n \geq 0$. Hence, $-|a|^2/4 + C
W(\varphi_{t_n}) \leq - |a-\varphi_{s_n}|^2 + C W(\varphi_{s_n})$.
Letting $n$ tend to $+ \infty$, we obtain a contradiction. \qed

\section{Two-Stacks Model}
\label{sec:maier}

In this section, we consider a special case. It is a generalization
of an interesting example introduced by Maier \cite{Maier}. With
$d=2$, $E=\{1,2\}$ and $x=(x_1,x_2)$, let
\begin{equation}
  \label{eq:maier}
  p(x,i,v) = \left\{
  \begin{array}{ll}
\ud \l_i  [1-g_1(x_1)]\;, & v=e_1\\
\ud \l_i  [1+g_1(x_1)]\;, & v=-e_1\\
\ud (1-\l_i)  [1-g_2(x_2)]\;, & v=e_2\\
\ud (1-\l_i)  [1+g_2(x_2)]\;, & v=-e_2
  \end{array}
\right.
\end{equation}
 with some $\l_i \in (0,1)$ for all $i \in E$, and some Lipschitz
 continuous functions
$g_1, g_2: [0,1] \to [0,1)$, $g_j(z)>0$ for $z > 0$.


When $g_1=g_2$ and $\l_i=1/2,
i \in E$, this example reduces to that of Maier (see (4) in
\cite{Maier}). Here, the random environment $\xi$ governs the
probability for each coordinate to jump, but not the jump
distribution itself.
Our treatment below is quite different from \cite{Maier},
being more direct and leading to more general results.


From (\ref{def:barf}) we compute
$$
\bar f (x) = - \left(
  \begin{array}{c}
 \l g_1(x_1) \\ (1-\l) g_2(x_2)
  \end{array}
\right)
 \;,
\qquad \l= \sum_{i \in E}
\l_i \mu(i).
$$
In this example, all the assumptions {\bf (A.1--4)} are satisfied.
The assumption of Theorem \ref{th:tloi} holds if $P^2$ is irreducible and $g_1',g_2' \geq \kappa'$ for some constant $\kappa'>0$.

If both $g_1(0)$ and $g_2(0)$ are equal to zero, then $\bar f (0)=0$ and the reflected differential
equation \eqref{eq_lim_ref} is simply the ordinary  differential
equation inside $G$. In this case, the hitting time of the stable equilibrium 0 is infinite.
If, on contrary,  $g_j(0)>0$ for some $j \in \{1,2\}$,
then the solution to the
RDE \eqref{eq_lim_ref} feels the reflection when hitting the $j$-th axis. After hitting the boundary, it moves towards the origin
along the $j$-th axis.


The function $H$ can be expressed in terms of
$$
H_j(x_j,\a_j)= \ln \big[ \cosh \a_j - g_j(x_j) \sinh \a_j\big]\;,\quad
j=1,2.
$$
From (\ref{eq:perron}), $H(x, \a)$ is the logarithm of the largest
eigenvalue
of the matrix
\begin{equation}
  \label{eq:Q}
Q(x,\a)=
\left[ P(i,j) \big\{ \l_i e^{H_1(x_1, \a_1)} +
(1-\l_i) e^{H_2(x_2, \a_2 )} \big\}
 \right]_{i,j \in E}.
\end{equation}
Recall that $E=\{1,2\}$. By solving the characteristic equation, we find,
with shorthand notations  $ P(i,j)= P_{ij}$,
\begin{equation*}
\begin{split}
&H(x, \a) =
\ln \ud \left( P_{11}A_1+P_{22}A_2 + \sqrt{ \big(P_{11}A_1-P_{22}A_2\big)^2+4
P_{12}P_{21}A_1A_2} \right)
\\
&{\rm with} \left\{ \begin{array}{l}
A_1(x,\a)= \l_1 e^{H_1(x_1, \a_1)} +
(1-\l_1) e^{H_2(x_2, \a_2 )} \; ,
\\
A_2(x,\a)= \l_2 e^{H_1(x_1, \a_1)} +
(1-\l_2) e^{H_2(x_2, \a_2 )}\;.
\end{array}
\right.
\end{split}
\end{equation*}
\subsection{Identification of the Quasi-potential}
Although its expression does not look very explicit, the
quasi-potential is quite simple. It can be guessed
by observing that the discrete walk $X_n$ has an invariant measure,
which obeys a large deviations principle: in view of \cite[Chapter 4, Theorem 4.3]{FV}, the rate function -- which is explicit here --
should be the quasi-potential.

In Maier's paper, the quasi-potential was identified by
a Lagrangian approach and using the special structure of the
separable Hamiltonian \cite[p.397]{Maier}. Our approach here is an
alternative yielding to a much shorter route for more general Hamiltonians.
\vspace{5pt}

We start to look for the invariant measure.
The Markov chain on $\{0,1/m,\dots,1\}$ with nearest neighbor
transitions $(1/2) [1\mp g_1(x_1)]$ from $x_1$ to $x_1 \pm 1/m$ (pay attention to the change of sign between $\mp$ and $\pm$) with reflection
at 0 and 1 has
an invariant (even reversible) measure given for $z=k/m $ by
$$
 \pi_1^{(m)}(k/m)=
 \frac{1}{1+g_1(\frac{k}{m})}
\prod_{l=0}^{k-1} \frac{1-g_1(\frac{l}{m})}{1+g_1(\frac{l}{m})}\;, \ 0 <k<m,
$$
and $\pi_1^{(m)}(0)=(1-g_1(0))/[2(1+g_1(0))]$ and $\pi_1^{(m)}(m) =
(1/2)
\prod_{l=0}^{m-1} (1-g_1(\frac{l}{m}))/(1+g_1(\frac{l}{m}))$.
When the function $g_1$ is
Lipschitz
continuous, we obtain for large $m$ and $z \in (0,1)$,
\begin{eqnarray}
\nonumber
  \pi_1^{(m)}(z)&=&
\exp \left\{ \sum_{l=1}^{[mz]} \ln
\frac{1-g_1(\frac{l}{m})}{1+g_1(\frac{l}{m})}
+ \cO (1)\right\}
\\ \label{eq:eqpi}
&=&
\exp \left\{ -2m \int_0^z \tanh^{-1} (g_1(y))  dy + o (m)\right\}
\end{eqnarray}
since $\tanh^{-1}(t)=(1/2) \ln [ (1+t)/(1-t)]$.
We define $\pi_2^{(m)}$
similarly, with $g_2$ instead of $g_1$. The second observation is that
the measure
\begin{equation}
  \label{eq:inv}
 \nu^{(m)}(x,i)=\pi_1^{(m)}(x^{1}) \pi_2^{(m)}(x^{2}) \mu(i)
\end{equation}
is invariant for our Markov chain $(X_n^{(m)}/m,\xi_n)_{n \geq 0}$. Indeed,
invariance of $\pi_1^{(m)}$ for the corresponding transition
implies
$$
\forall y \in [0,1]^2 \cap (m^{-1} {\mathbb Z}^2), \
\sum_{x_{1} \in \{0,1/m,\dots,1\}} \pi_1^{(m)}(x_1) q\bigl((x_1,y_2),i,(y_1-x_1)e_1 \bigr) = \l_i \pi_1^{(m)}(y_1).
$$
Hence, for all $j \in E$ and $y \in [0,1]^2 \cap (m^{-1} {\mathbb Z}^2)$,
\begin{equation*}
\begin{split}
&\sum_{i \in E} \sum_{x \in \{0,1/m,\dots,1\}^2} \pi_1^{(m)}(x_{1})
  \pi_2^{(m)}(x_{2})
\mu(i) P(i,j) q(x,i,y-x)
\\
&= \sum_{i \in E} \sum_{|x_1-y_1|=1,x_2=y_2 \ {\rm or}
|x_2-y_2|=1,x_1=y_1}
\pi_1^{(m)}(x_{1})
  \pi_2^{(m)}(x_{2})
\mu(i) P(i,j) q(x,i,y-x)
\\
&=  \sum_{i \in E}  P(i,j) \mu(i) [\l_i + (1-\l_i)]
\pi_1^{(m)}(y_{1}) \pi_2^{(m)}(y_{2})
\\
&=  \pi_1^{(m)}(y_{1})  \pi_2^{(m)}(y_{2})
\mu(j)
\end{split}
\end{equation*}
As a by product, the first marginal $\nu^{(m)}_0$ of $\nu^{(m)}$, i.e.
$ \nu^{(m)}_0(y)= \pi_1^{(m)}(y_{1})
  \pi_2^{(m)}(y_{2}) $,
is itself invariant for $(X_n^{(m)}/m)_n$ (which is not
a  Markov chain). From the relation (\ref{eq:eqpi})
it is clear that
this new measure satisfies a large deviations principle,
with rate function
\begin{equation}
  \label{eq:W}
  W(x)= 2 \int_0^{x_{1}} \tanh^{-1} (g_1(y))  dy +2
\int_0^{x_{2}} \tanh^{-1} (g_2(y))  dy
\end{equation}
By \cite[Chapter 4, Theorem 4.3]{FV}, we then expect $W$ to be the quasi-potential.
By Proposition \ref{prop:identification}, we prove that 
this equality indeed holds.

\begin{thm} \label{th:V=W}
The function $W$ coincides with the quasi-potential. Moreover, for any point $x \in [0,1]^d \setminus \{0\}$,
there is one and only one optimal path $(\varphi_t)_{t \leq 0}$ from $0$ to $x$. The time reversed path $(\varphi_{-t})_{t \geq 0}$ is the unique
solution to the reflected differential equation given by the law of large numbers (see Corollary \ref{conv_pb}), i.e. $\varphi_t =\chi_{-t}^x$ for all $t \leq 0$.
\end{thm}

\noindent $\Box$ Proof. We check that all the assumptions of Proposition \ref{prop:identification} are fulfilled.

\textit{First Step. Hamilton-Jacobi Equation.}
The function $W$ is clearly smooth. The gradient is given by
\begin{equation}
\label{eq:2003}
\nabla W(x)=2 \bigl(\tanh^{-1}(g_1(x_1)),
\tanh^{-1}(g_2(x_2)) \bigr)
\end{equation}
On the boundary, $\nabla_+ W(x) = 2(0,\tanh^{-1} (g_2(x_2)))$ for $x_1=0$ and $x_2 \in (0,1]$, $\nabla_+ W(x) = 2(\tanh^{-1}(g_1(x_1)),0)$ for
$x_1 \in (0,1]$ and $x_2=0$, $\nabla_+ W(x) = \nabla W(x)$ for $x_1=1$ and $x_2 \in (0,1]$ and for $x_1 \in (0,1]$ and $x_2=1$.
For $x=0$, we have $\nabla_+ W(x)=0$.

We recall the hyperbolic trigonometric identities
$$
\tanh a = \frac{2g}{1+g^2}\;,\quad \sinh a = \frac{2g}{1-g^2}\quad
{\rm for\ } a=2 \tanh^{-1} (g)\;
$$
For $j=1,2$, the quantity
$\exp(H_j(x_j,\a_j))= \cosh (\a_j) - g_j(x_j) \sinh (\a_j)$ is equal to 1 iff
$\a_j=0$ or $\a_j = 2 \tanh^{-1}g_j(x_j)$. From \eqref{eq:Q}, we deduce that
$Q(x,\nabla_+ W(x))=P$ for every $x \in [0,1]^2$.
(With $\nabla_+ W(x) = \nabla W(x)$ for $x \in (0,1)^2$.)
 Hence, the largest eigenvalue of
$Q(x,\nabla_+ W(x))$ is 1. We deduce that $W$ satisfies the Hamilton-Jacobi equation (\ref{eq:HJBclass1})
--(\ref{eq:HJBclass2}).
\vspace{5pt}

\textit{Second Step. Identification of $W$.}
We first compute the gradient of $H$, with respect to $\alpha$, in $(x,\nabla_+ W(x))$, $x \in [0,1]^d$.
Since, for all $j \in \{1,2\}$,
\begin{eqnarray}\nonumber
  \frac{\partial H_j}{  \partial \a_j}
(x_j,\a_j)&=& \frac{\sinh(\a_j) - g_j(x_j) \cosh (\a_j)}
{ \cosh (\a_j) - g_j(x_j) \sinh (\a_j)}\\\nonumber
&=& g_j(x_j) \qquad \qquad {\rm for \ } \a_j= [\partial W/\partial x_j](x)\;,
\end{eqnarray}
we have
\begin{equation}
  \label{eq:gradQ}
\frac{\partial Q}{  \partial \a_1}
(x, \nabla W(x)) = \left[ P(i,j)  \l_i g_1(x_1)
 \right]_{i,j \in E}, \ \frac{\partial Q}{  \partial \a_2}
(x, \nabla W(x)) = \left[ P(i,j) (1- \l_i) g_2(x_2)
 \right]_{i,j \in E}.
\end{equation}
By simplicity of the top eigenvalue we know that $H(x,\cdot)$ is
differentiable. For the same reason, the associated eigenvector
$v(x,\alpha)$ is smooth in $\alpha$. We thus differentiate the equation
$Q(x,\a) v(x,\a)=  \exp(H(x,\a)) v(x,\a)$ at $\alpha= \nabla W(x)$.
At such a point,  $Q=P$, $H=0$ and $v={\1}=(1,\ldots,1)^t$, so that
\begin{equation}
  \label{eq:1142}
\frac{\partial Q}{  \partial \a_1} (x, \nabla W(x)) {\1}
+ P \frac{\partial v}{  \partial \a_1} (x, \nabla W(x))
=
\frac{\partial H}{  \partial \a_1} (x, \nabla W(x)) {\1}
+ \frac{\partial v}{  \partial \a_1}(x, \nabla W(x)).
\end{equation}
From (\ref{eq:gradQ}) we have $[\partial Q /\partial \a_1]
(x, \nabla W(x)) {\1} = g_1(x_1) {\vec \lambda}$ with
$ {\vec \lambda}= (\l_i)_{i \in E}$, and by
multiplying (\ref{eq:1142}) by the invariant measure $\mu$ on the left,
we get $[\partial H/\partial \a_1] (x, \nabla W(x))
= \l g_1(x_1)$. With a similar computation for the
partial derivative with respect to $ \a_2$, we finally obtain
\begin{equation}
  \label{eq:1921}
\nabla_\a  H(x, \nabla W(x)) = \left(
  \begin{array}{c}
\l g_1(x_1)\\ (1-\l) g_2(x_2)
  \end{array}
\right)   = - \bar f (x)\;.
\end{equation}
Repeating the computations from \eqref{eq:2003} to \eqref{eq:1921}, we have
\begin{equation}
\label{eq:d2verif1}
\nabla_\a  H \bigl(x,\nabla_+ W(x) \bigr)
= \left\{
\begin{array}{l}
 \left(
  \begin{array}{c}
\bar{f}_1(0)\\ - \bar{f}_2(x_2)
  \end{array}
\right), \ {\rm for} \ x_1 = 0, \ x_2 \in (0,1],
\\
 \left(
  \begin{array}{c}
- \bar{f}_1(x_1)\\  \bar{f}_2(0)
  \end{array}
\right), \ {\rm for} \ x_1 \in (0,1], \ x_2=0.
\end{array}
\right.
\end{equation}
It is plain to check that the assumptions of Proposition \ref{prop:optimsuff} are fulfilled.
Therefore, $W$ is the quasi-potential.

\textit{Third Step. Optimal Paths.} For a terminal value $x \in [0,1]^2 \setminus \{0\}$, we have to prove that the time reversed path $(\chi_{-t}^x)_{t \leq 0}$
satisfies \eqref{eq:optimsuff1} as well as \eqref{eq:optimsuff2}.
It is sufficient to prove it locally: we prove that, for any $x$, $(\chi_{-t}^x)_{-\varepsilon \leq t \leq 0}$  satisfies both
\eqref{eq:optimsuff1} and \eqref{eq:optimsuff2} on a small interval $[-\varepsilon,0]$ for some $\varepsilon >0$.
By \eqref{eq:1921}, this is easily checked if the terminal point $x$ belongs to $(0,1)^2$.
If the terminal point $x$
belongs to the boundary, several cases are to be considered.

If $x_1=0$ and $x_2 \in (0,1]$, the path $(\chi_t^x)_{t\geq 0}$ remains on $\{0\} \times [0,1]$, so that 
$(\dot{\chi}_t^x)_1=0$ for $t \geq 0$. For some $\varepsilon>0$, we have
$(\chi_t^x)_2 \in (0,1]$ for $t \in [0,\varepsilon]$. By \eqref{eq:d2verif1}, the second coordinate 
satisfies
$(\dot{\chi}_t^x)_2 = - \langle \nabla_{\alpha} H(\chi_t^x,\nabla_+ W(\chi_t^x)),e_2\rangle$ for $t \in [0,\varepsilon]$.
Setting $\varphi_t = \chi_{-t}^x$ for all $t \in [-\varepsilon,0]$, we have
\begin{equation*}
\dot{\varphi}_t = \nabla_{\alpha} H\bigl(\varphi_t,\nabla_+ W(\varphi_t) \bigr) - \bar{f}_1(0) e_1,
\end{equation*}
so that $(\varphi_t)_{\varepsilon \leq t \leq 0}$ satisfies \eqref{eq:optimsuff1}.
Since $\langle\nabla_+ W(0,y),e_1 \rangle = 0$ for all $y \in [0,1]$, the compatibility condition is fulfilled.
The same holds if $x_2=0$ and $x_1 \in (0,1]$.

If $x_1=1$ and $x_2 \in (0,1]$, then the path $(\chi_t^x)_{t \geq 0}$ leaves the boundary immediately: for $t>0$ (and $t$ small), $\chi_t^x
\in (0,1)^2$. Reversing the path,
we conclude as above. The same holds if $x_2=1$ and $x_1 \in (0,1]$.

\textit{Fourth Step. Uniqueness of the Optimal Path.}
It remains to verify that the solutions to \eqref{eq:optimsuff1} are unique. Again, it is sufficient to prove that uniqueness holds locally for
any starting point in $[0,1]^2 \setminus \{0\}$. If the
starting point is in $(0,1)^2$, this is obvious by time reversal. If $x_1=0$ and $x_2 \in (0,1]$, we have
$\langle \nabla_{\alpha} H(x,\nabla W(x)),e_1 \rangle = - \bar{f}_1(0) \geq 0$.
Assume for the moment that $-\bar{f}_1(0) >0$. Then, by Figure \ref{figu:optimal1}, any solution $(\varphi_t)_{t \geq 0}$
to \eqref{eq:optimsuff1}
touches the boundary before reaching $x$. Hence, there exists $\varepsilon >0$ such that $(\varphi_t)_1=0$ and $(\dot{\varphi}_t)_2= -\bar{f}_2((\varphi_t)_2)$ for all $t \in
[-\varepsilon,0]$. Local uniqueness easily follows. Assume now that $\bar{f}_1(0)=0$. Then, for all $y$ in the neighborhood of $x$, with $y_1>0$,
$\langle \nabla_{\alpha} H(y,\nabla W(y)),e_1 \rangle = - \bar{f}_1(y_1) > 0$. Again,
any solution $(\varphi_t)_{t \geq 0}$
to \eqref{eq:optimsuff1}
has to touch the boundary before reaching $x$ (otherwise, it cannot reach the boundary) and we can repeat the argument.
The same holds for $x_2=0$ and $x_1 \in (0,1]$. The case where $x_1=1$ and $x_2 \in (0,1)$
corresponds (up to a symmetry) to Figure \ref{figu:optimal2} and local uniqueness is proved in a similar way. The cases where 
$x_2=1$ and $x_1 \in (0,1)$ and where $(x_1,x_2)=(1,1)$ are similar.
\qed

\medskip

\subsection{Deadlock Phenomenon for the Two-Stacks Model}
We discuss the deadlock phenomenon for the two-stacks model, that is for the domain $G$ from (\ref{ensG}) with $\ell =1$.
Our results should
be compared to Section 5 in \cite{Maier}. In view of Theorems \ref{th:tmoy}
and \ref{th:V=W},
the set of exit points $\cM$
relates to the simple, one-dimensional, variational problem
$$
\bar V = \min\{ W(z,1-z); \, z \in [0,1]\}.
$$
Then, $x \in \cM$ if and only if $x=(z,1-z)$ with $z$ minimizing the
above problem.
Observing that
$[d/dz](W(z,1-z))=2 (\tanh^{-1} (g_1(z))-\tanh^{-1} (g_2(1-z)))$ has
the same sign as $ g_1(z)-g_2(1-z)$,
we distinguish a few remarkable different regimes (some of them being
discussed in \cite{Maier})
for the set $\cM$ of  deadlock configurations and for the shape of the optimal
paths (which describe the  typical course of a deadlock).
\medskip

{\bf Qualitative shape of  optimal paths.}  For $x=(x_1,x_2)
\in (0,1)^2$, we discuss
the optimal path $(\varphi_t)_{t \leq 0}$ from 0 to $x$. By
Theorem \ref{th:V=W}, $\varphi_t = \chi_{-t}^x$ for all $t \leq 0$.  As long as the $k$-th coordinate ($k=1,2$) of $\varphi_t$
is positive, it satisfies $(\dot{\varphi}_t)_k= - \bar{f}_k((\varphi_t)_k)$. Hence, the
time needed to make  the $k$-th coordinate move from 0 to $x_k$ is
$$
t_k=t_k(x_k)= - \int_0^{x_k} \frac{dr}{\bar{f}_k(r)} \in (0, \8].
$$
Note that  $t_k$ is finite if the continuous function $\bar{f}_k$ (or equivalently $g_k$) is non zero at 0,
but  $t_k$ is infinite if $\bar{f}_k(0)=0$ (since $\bar{f}_k$ is Lipschitz continuous,
$|\bar{f}_k(r)|\leq Cr$ for $r>0$ in this case).
In general, the duration of the instanton $\varphi$ from 0 to $x$
is equal to $\max\{t_1,t_2\}$.
\begin{enumerate}
\item {\bf (Case A).} $g_1(0)= g_2(0)=0$. Then $\varphi$ has an infinite duration.
It never hits the
boundary and does not feel the
reflection. When $g_1=g_2$, $x_1=x_2$ and $\lambda=1/2$, the  optimal path
is the line
segment $[x,0]$. But in general, the  optimal path is not a line.
\item  {\bf (Case B).} $g_1(0)>0, g_2(0)>0$. Then,  the  optimal path
 has a finite duration. There is a smooth curve of points $x$'s such that
the  reversed path  from $x$ to $0$ does not hit the axis (strictly) before
0: the curve is in fact defined by $t_1(x_1)=t_2(x_2)$.
For $x$'s such that $t_1(x_1)<t_2(x_2)$,  $\chi^x$  hits the vertical axis (strictly)
above 0, and later on, moves down towards 0 along this axis.
\item  {\bf (Case C).} $g_1(0)>0, g_2(0)=0$. For all $x \in (0,1)^2$,
$\varphi$ hits the vertical axis in a finite time, and later on, moves down towards
0 along this axis  reaching it in infinite time.
\end{enumerate}
Some optimal paths are shown in Figures \ref{fig:cas1} and \ref{fig:cas1_5}
below.

\begin{figure} [htb]
\begin{center}
\includegraphics[
width=0.4\textwidth,angle=0]
{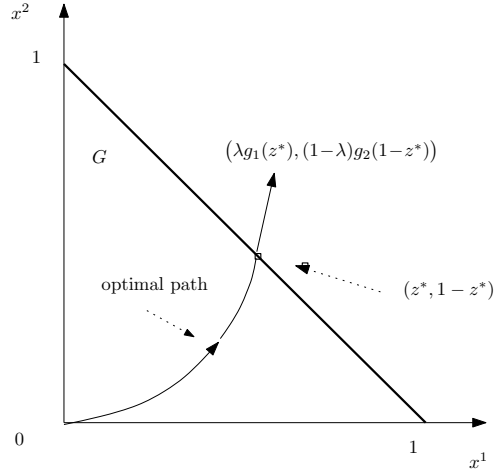}
\caption{Optimal  deadlock point and path, Case 1 with $\l<1/2$, $g_1(0)=
g_2(0)=0$}
\label{fig:cas1}
\end{center}
\end{figure}

\begin{figure} [htb]
\begin{center}
\includegraphics[
width=0.4\textwidth,angle=0]
{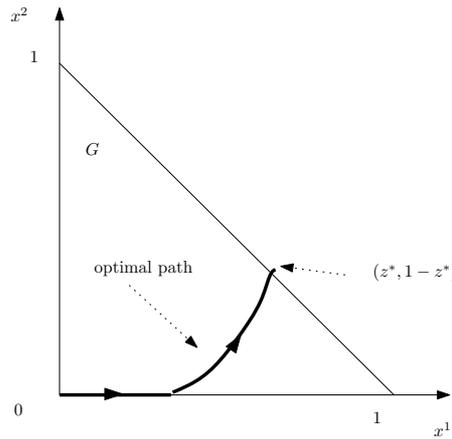}
\caption{Optimal deadlock point and path, Case 1 with $g_1(0)>0,
g_2(0)>0$, and the exit point $(z^*,1-z^*)$ in the general situation}
\label{fig:cas1_5}
\end{center}
\end{figure}

\medskip

{\bf Some specific cases for the set $\cM$.}
\begin{enumerate}
\item {\bf (Case 1).} Assume that $g_1$ and $g_2$ are strictly increasing on $[0,1]$. Then
$x_1 \in [0,1] \mapsto W(x_1,1-x_1)$ is a strictly convex function so that $\cM$ reduces to a single point.
If $g_1(0) \geq g_2(1)$, then the function is increasing and the minimum is attained at $x_1=0$, so that
$\cM=\{(0,1)\}$. If $g_1(1) \leq g_2(0)$, the function is decreasing and the minimum is attained at $x_1=1$, so that
$\cM=\{(1,0)\}$. If $g_1(0) < g_2(1)$ and $g_2(0)<g_1(1)$, then the slope is negative at 0 and positive at 1, so that
$\cM=\{(z^*,1-z^*)\},$
with $ z^* \in (0,1)$,
 the\ unique\ solution\ of\
$g_1(z^*)=g_2(1-z^*)$.
This case is illustrated by Figures \ref{fig:cas1} and \ref{fig:cas1_5}.

\item  {\bf (Case 2)}
Assume $g_1(z)= g_2(1-z)$ for all $z \in [a,b]$ ($0<a<b<1$), and
$g_1$ [resp. $g_2$] strictly increasing on $[0,a] \bigcup [b,1]$
[resp. on  $[0,1-b] \bigcup [1-a,1]$].
Then,  $ g_1(z)-g_2(1-z)$ -- as well as $[d/dz](W(z,1-z))$ -- is
increasing [resp. zero, increasing] on the interval $[0,a]$
[resp., $[a,b], [b,1]$]. Now, the set of minimizers is the interval,
$$
\cM = {\rm segment\ } \big[(a,1-a), (b,1-b)\big]
$$
as indicated in Figure \ref{fig:cas2}.
\begin{figure} [htb]
\begin{center}
\includegraphics[
width=0.4\textwidth,angle=0] {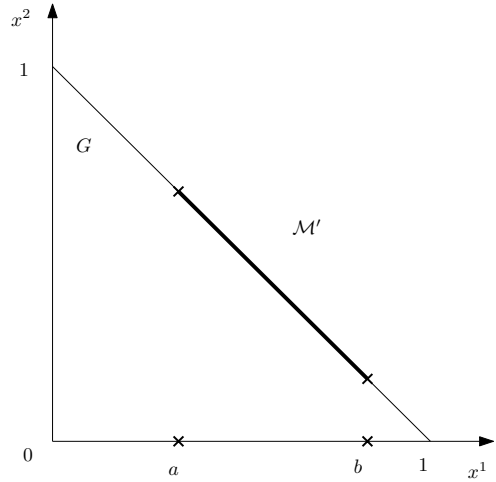} \caption{Optimal  deadlock
points, Case 2} \label{fig:cas2}
\end{center}
\end{figure}
\item  {\bf (Case 3)}
Assume $g_1(z)-g_2(1-z)$ is negative on  $[0,a)$, positive  on
$(a,c)$,  negative on $(c,b)$ and  positive  on
$(b,1]$  ($0<a<c<b<0$). Then, $W(z,1-z)$ is a double-wells, and the set of
minimizers is a pair,
$$
\cM = \{ (a,1-a), (b,1-b)\}
$$
see Figure \ref{fig:cas3}.
\begin{figure} [htb]
\begin{center}
\includegraphics[
width=0.4\textwidth,angle=0]
{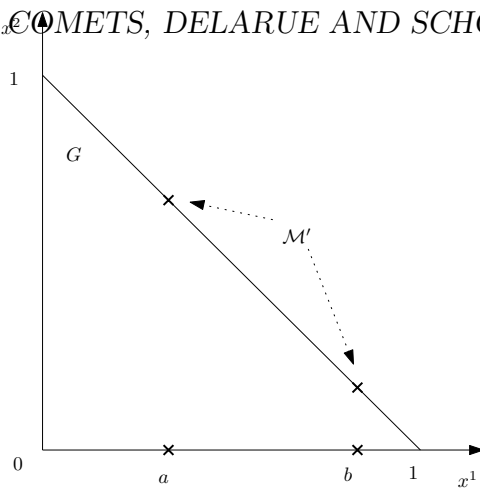}
\caption{Optimal deadlock points, Case 3}
\label{fig:cas3}
\end{center}
\end{figure}
\end{enumerate}

By Theorem \ref{th:V=W}, there is a one-to-one
correspondence between elements of $\cM$ and optimal path
(so-called instantons) to exit $G$.
Therefore, there is a unique optimal path for the deadlock in Case 1,
uncountably many in Case 2, and exactly two in Case 3.

\section{Limit Cycle} \label{sec-limtcycle}

In this section we work out an example where the system has,
in the large scale limit $m \to \8$,
a stable attractor, which is a limit cycle. Denote by $\1$ the vector $(1,1)^t$,
and consider the differential system in $\R^2$,
$$\dot x_t = h(x_t)$$
with
\begin{equation}
\label{ode-cycle}
 h(x) = \left(
\begin{array}{cc}
              0&-1\\1&0
           \end{array} \right)
\Big(x-\frac{1}{4} {\1}\Big) + \frac{1}{2} \Big[1-64\times|x-\frac{1}{4}
{\1}|^2\Big]
\Big(x-\frac{1}{4} {\1}\Big)\;,
\end{equation}
\begin{figure} [htb]
\begin{center}
\includegraphics[
width=0.4\textwidth,angle=0]
{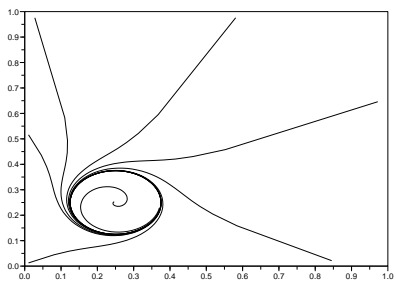}
\caption{An ordinary differential equation with a limit cycle}
\label{fig:007}
\end{center}
\end{figure}

\noindent
whose phase portrait is given in Figure \ref{fig:007}.
The circle $\mathcal C$ centered at
$1/4$ with radius $1/8$ is a stable limit set:
trajectories spiral into it  as time approaches infinity.
More precisely, it can be
checked that any point in $[0,1]^2$ is attracted by  $\mathcal C$. Moreover,
the vector field on the axis is pointing inside the first quadrant,
and, for $\ell \in [1,2]$, the vector field
on the sloping side $|x|_1=\ell$
is pointing inside the domain $G$.

Obviously, the reason for the existence of the limit cycle
is that the vector field is the superposition of
\begin{equation}
\label{eq-h1}
h_1(x)= \left(
\begin{array}{cc}
              0&-1\\1&0
           \end{array} \right)
\Big(x-\frac{1}{4} {\1}\Big)
\end{equation}
-- a rotation around the center $(1/4)  \1$ which preserves
the norm of the vector
$(x-(1/4) { \1})$--, and of
\begin{equation}
\label{eq-h2}
h_2(x)= \frac{1}{2} \big[1-64\times|x-\frac{1}{4}
{\1}|^2\big]
\Big(x-\frac{1}{4} {\1}\Big)
\end{equation}
-- whose effect is moving the system on the radius issued at the center
towards the intersection of the radius and the circle $\mathcal C$-- .
It is plain to check that the components of
$h_1, h_2$ are bounded on $[0,1]^2$ by
a constant strictly smaller than 32. Let $E=\{1,2\}$,
and
 assume that $\mu$ is the Bernoulli law $\mu(1)=\mu(2)=1/2$.
With $h^{(1)}, h^{(2)}$ the  components of $h \in \R^2$,
define the transition by
\begin{equation}
\label{eq-pcycle}
p(x,i,v) =  \left\{
\begin{array}{cc}
\frac{1}{4} \big( 1  \pm \frac{1}{32}h_i^{(1)}(x)\big) & {\rm for \ } v=\pm
e_1\\
\frac{1}{4} \big( 1  \pm \frac{1}{32}h_i^{(2)}(x)\big) & {\rm for \ } v=\pm
e_2
\end{array} \right.\quad, \qquad i=1,2
\end{equation}
Since  $\mu$ is Bernoulli, the limit ordinary differential equation
(\ref{eq:ODE}) is given here by
$$
\bar f (x) = \frac{1}{128} h(x)
$$
with $h$ from
(\ref{ode-cycle}).
 Assumptions {\bf (A.1--3)} are fulfilled, as well as the counterpart to
{\bf (A.4)} -- with the attractor $\mathcal C$ replacing the stable fixed point
 $0$. Most of the results
of Section \ref{sec-deadlock} can be generalized to this case,
with the quasi-potential computed as the minimal action over all paths
from ${\mathcal C}$ to the current point.
For instance,
(\ref{eq:tmoy1}) becomes
$$
\bbE_{x_m}[ \tau^{(m)}] = \exp \bigl[  m (\hat{V}+o(1)) \bigr]
$$
for  any sequence $x_m \rightarrow x \in G$, with
$$\hat{V} = \inf \{ \J0T(\phi); \, \phi_0\in {\mathcal C}, \,
|\phi_T|_1=\ell, \, T>0\}\;.$$
We cannot compute the exact value of the quasi-potential in this
example, but it could be estimated numerically from above.
Following \cite[Chapter 5, Theorem 4.3]{FV}, we could also provide a suitable version of Proposition \ref{prop:identification}.

\mysection{Appendix A}

\subsection{Proof of Lemma \ref{couplage}: successful coupling}
\label{sec:preuve-couplage}

The proof relies on a tricky coupling argument.
In \cite{OV},
the authors investigate the large deviations for stochastic differential
equations with a small noise: the coupling argument then follows from
standard arguments for the Brownian motion. In our own
setting, the standard stochastic analysis tools are useless and we need to construct a coupling for our purpose.

{\bf Coupling}. For an initial condition $x \in {\mathbb Z}^d$, $|x|_1 < \delta m$, the position of the walker is given by
\begin{equation*}
X_{n+1} = (2\Pi^{(m)} - \Id) \bigl(X_n + f(X_n/m,\xi_n,U_n) \bigr),
\end{equation*}
with $X_0 = x$. Here, $f(z,i,\cdot)$ denotes a function from $(0,1)$ to ${\mathcal V}$ such that $f(z,i,U)$ has
$p(z,i,\cdot)$ as distribution (typically, $f(z,i,\cdot)$ is an inverse of the cumulative distribution function
of $p(z,i,\cdot))$. For another initial condition $y \in {\mathbb Z}^d$, $|y|_1 < \delta m$, the position of the walker can be defined in a similar way.
The realizations $(U_n)_{n \geq 0}$ may be the same. Nevertheless, the position may be defined with a different sample of uniform law.
It may be also defined with the same sample but with a different function $f$.

In what follows, we are seeking for a copy $(\hat{X}_n)_{n \geq 0}$ of the walk, starting from $y$, such that $\hat{X}$ and $X$ join up in a finite time.
For this purpose, we assume $|x-y|_1 \in 2 {\mathbb N}$ (otherwise, it is impossible). We will use the same
sample of uniform law but a different function $f$. We thus write
\begin{equation*}
\hat{X}_{n+1} = (2 \Pi^{(m)} - \Id) \bigl(\hat{X}_n + \hat{f}_n(U_n) \bigr),
\end{equation*}
where $\hat{f}_n$ is some random function from $(0,1)$ into ${\mathcal V}$, depending on $X_n$,
$\hat{X}_n$ and $\xi_n$ such that the conditional law of $\hat{f}_n(U_n)$ with respect to $(X_n,\hat{X}_n,\xi_n)$ is exactly $p(\hat{X}_n,\xi_n,\cdot)$. The explicit form
of $\hat{f}_n$ has to be determined.

To simplify, we will just denote (when possible) $\hat{f}_n(U_n)$ by $\hat{f}_n$. Similarly, we will denote $f(X_n/m,\xi_n,U_n)$ by $f_n$ (or $f_n(U_n)$
when necessary).

Before providing an explicit form for $\hat{f}_n$, we investigate the $L^1$-distance $\Delta_{n} = |X_{n} - \hat{X}_{n}|_1$. Loosely speaking, we want
it to
decrease with $n$. We thus compute $\Delta_{n+1}$ in terms of $\Delta_n$. For this purpose, it is crucial to note that $\Delta_n$ is always even (because of the particular
choice for the initial conditions and for the reflection).
%
We also
recall the formula
\begin{equation*}
\forall a,b \in \R, \ |a+b| = |a| + |b| - 2 (|a|  \wedge |b|) {\mathbf 1}_{\{ab<0\}}.
\end{equation*}
If $X_n$ and $\hat{X}_n$ are not on the boundary, we deduce
\begin{equation}
\label{eq:couplage00}
|\Delta_{n+1}|_1
= |\Delta_n|_1 + |f_n-\hat{f}_n|_1
- 2 \sum_{i=1}^d
[|(\Delta_n)_i|\wedge |(f_n)_i-(\hat{f}_n)_i|] {\mathbf 1}_{\{
(\Delta_n)_i((f_n)_i-(\hat{f}_n)_i) <0\}}.
\end{equation}
If one of the two processes is on the boundary at time $n$, the difference
$|\Delta_{n+1}|_1 - |\Delta_n|_1$ has the form $|\Delta_n+ g_n - \hat{g}_n|_1 - |\Delta_n|_1$ with $g_n$ and $\hat{g}_n$
as in \eqref{eq:form_g}. We can check that it is always bounded by $|\Delta_n +  f_n - \hat{f}_n|_1 - |\Delta_n|_1$. In other words, we can forget the reflection.
To prove this assertion, it is sufficient to focus on each coordinate. If $(X_n)_i=0$ and $(\hat{X}_n)_i \geq 2$, the proof is obvious.
If $(X_n)_i=0$ and $(\hat{X}_n)_i=1$, the proof is the same except for $(\hat{X}_{n+1})_i=0$ and $(X_n)_{i+1}=1$. In this case, the processes switch.
However, the result is still true. Other cases are treated in a similar way. Hence, in any case, \eqref{eq:couplage00} is true with $\leq$ instead of $=$.

Turn back to \eqref{eq:couplage00}.
Again,  $|f_n-\hat{f}_n|_1$ is always equal to 2, except for
$f_n=\hat{f}_n$. To handle the last term, we introduce the following notations
\begin{equation*}
E^+_n ({\rm resp.} \ E^-_n, \ {\rm resp.} \ E^0_n) = \{u
\in {\mathcal V}:\,  \langle \Delta_n,u \rangle >0 \ ({\rm resp.}  <0, \
{\rm resp.} =0)\}.
\end{equation*}
If $f_n = - \hat{f}_n \in E_n^{-}$, the sum is equal to  $\min(|\langle \Delta_n,f_n \rangle|,2)$.
If $f_n = \hat{f}_n$ or $f_n=- \hat{f}_n \in E_n^+ \cup E_n^0$, the sum is zero. If $f_n \perp \hat{f}_n$, $|(f_n)_i-(\hat{f}_n)_i|$ is 0 or 1
and $(|(\Delta_n)_i| \wedge|(f_n)_i-(\hat{f}_n)_i|) {\mathbf 1}_{\{(\Delta_n)_i ((f_n)_i - (\hat{f}_n)_i)<0\}}=
|(f_n)_i-(\hat{f}_n)_i| {\mathbf 1}_{\{(\Delta_n)_i ((f_n)_i - (\hat{f}_n)_i)<0\}}$
 is also
0 or 1: it is equal to 1 if and only if $f_n \in E_n^{-}$ and $i$ is the coordinate of $f_n$ or $\hat{f}_n \in E_n^+$ and $i$ is the coordinate of $\hat{f}_n$.
Hence,
\begin{equation*}
\begin{split}
|\Delta_{n+1}|_1
&\leq |\Delta_n|_1 + 2 - 2 {\mathbf 1}_{\{f_n=\hat{f}_n\}}
- 2 {\mathbf 1}_{\{f_n \perp \hat{f}_n\}} \bigl( {\mathbf 1}_{\{f_n \in E_n^-\}} + {\mathbf 1}_{\{\hat{f}_n \in E_n^+\}} \bigr)
\\
&\hspace{15pt}
 - 2 [|\langle \Delta_n,f_n \rangle|\wedge 2 ]{\mathbf 1}_{\{f_n=-\hat{f}_n,f_n \in E_n^-\}}.
\end{split}
\end{equation*}
Noting that $\{f_n \perp \hat{f_n}\}$ is the complementary of $\{f_n=\hat{f}_n\} \cup\{f_n = -\hat{f}_n\}$, we have
\begin{equation*}
\begin{split}
|\Delta_{n+1}|_1
&\leq |\Delta_n|_1 + 2 - 2 {\mathbf 1}_{\{f_n=\hat{f}_n\}}
- 2  \bigl( {\mathbf 1}_{\{f_n \in E_n^-\}} + {\mathbf 1}_{\{\hat{f}_n \in E_n^+\}} \bigr)
\\
&\hspace{15pt} + 2 \bigl( {\mathbf 1}_{\{f_n \in E_n^-\}} + {\mathbf 1}_{\{\hat{f}_n \in E_n^+\}} \bigr)
\bigl( {\mathbf 1}_{\{f_n= \hat{f}_n\}} + {\mathbf 1}_{\{f_n = - \hat{f}_n \}} \bigr)
\\
&\hspace{15pt}
 - 2 [|\langle \Delta_n,f_n \rangle|\wedge 2 ]{\mathbf 1}_{\{f_n=-\hat{f}_n,f_n \in E_n^-\}}.
\end{split}
\end{equation*}
We have $\{f_n=\hat{f}_n\} = \{f_n=\hat{f}_n \in E_n^+ \} \cup \{f_n=\hat{f}_n \in E_n^- \}
\cup \{f_n=\hat{f}_n \in E_n^0 \}$. Moreover $\{f_n= - \hat{f}_n \in E_n^-\} =
\{\hat{f}_n= - f_n \in E_n^+\}$. Hence,
\begin{equation*}
\begin{split}
|\Delta_{n+1}|_1
&\leq |\Delta_n|_1 + 2 - 2 {\mathbf 1}_{\{f_n=\hat{f}_n \in E_n^0\}}
- 2  \bigl( {\mathbf 1}_{\{f_n \in E_n^-\}} + {\mathbf 1}_{\{\hat{f}_n \in E_n^+\}} \bigr)
\\
&\hspace{15pt} + 4 {\mathbf 1}_{\{f_n= - \hat{f}_n , f_n\in E_n^-\}}
 - 2 [|\langle \Delta_n,f_n \rangle|\wedge 2 ]{\mathbf 1}_{\{f_n=-\hat{f}_n,f_n \in E_n^-\}}.
\end{split}
\end{equation*}
Finally,
\begin{equation}
\label{expr_couplage}
\begin{split}
|\Delta_{n+1}|_1
&\leq |\Delta_n|_1 + 2 {\mathbf 1}_{\{f_n \in E_n^+\}} - 2 {\mathbf 1}_{\{\hat{f}_n \in E_n^+\}}
+ 2 {\mathbf 1}_{\{f_n \in E_n^0,f_n \not = \hat{f}_n\}}
\\
 &\hspace{15pt}
+ 2 {\mathbf 1}_{\{f_n= - \hat{f}_n , f_n\in E_n^-\}}
 {\mathbf 1}_{\{|\langle \Delta_n,f_n\rangle| =1\}}.
\end{split}
\end{equation}
We claim that, for $n < \sigma = \inf \{k \geq 0: \,
|\Delta_k|_1=0 \} \wedge \inf\{k \geq 0: \,
|\Delta_k|_1 > 2 \lfloor m \delta^{1/2} \rfloor \}$ ($\delta$ small enough), we can choose
$\hat{f}_n$ such that $\{f_n=-\hat{f}_n,f_n \in E_n^-,|\langle \Delta_n,f_n\rangle| =1\}$ is empty and such that
${\mathbb P}\{f_n \in E_n^0,f_n \not = \hat{f}_n| \, | {\mathcal F}_n\} \leq \sum_{u \in E_n^0} [
p(X_n,\xi_n,u) - p(\hat{X}_n,\xi_n,u)]^+$, with ${\mathcal F}_n = {\mathcal F}_n^{
\xi,X,\hat{X}}$.

The idea is the following. We define the random sets (i.e. they may depend on $\xi_n$, $X_n$ and $\hat{X}_n$): $A_n(u) = \{r \in (0,1): \,  f_n(r) =u \}$ and
$\hat{A}_n(u) = \{r \in  (0,1): \,  \hat{f}_n(r)=u\}$ for $u \in {\mathcal V}$. The Lebesgue measures of these sets are known: $|A_n(u)| = p(X_n,\xi_n,u)$
and $|\hat{A}_n(u)| = p(\hat{X}_n,\xi_n,u)$. In the sequel, we just write $p_n(u)$ and $\hat{p}_n(u)$ for these quantities.

\begin{figure}[htb]
\begin{center}
\psfrag{0}{$0$}
\psfrag{1}{$1$}
\psfrag{en0}{$E_n^0$}
\psfrag{-v1}{$-v_1$}
\psfrag{-v2}{$-v_2$}
\psfrag{-v3}{$-v_3$}
\psfrag{v1}{$v_1$}
\psfrag{v2}{$v_2$}
\psfrag{v3}{$v_3$}
\includegraphics[
width=0.5\textwidth,angle=0]{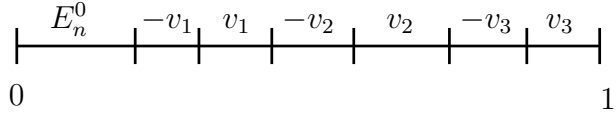}
\caption{Order for $f_n$, $M=3$.}
\label{fig:coupling}
\end{center}
\end{figure}

For each $u \in {\mathcal V}$, $A_{n}(u)$ is
an interval (because of the construction by inversion of the
 cumulative distribution
function). However, the geometry of $\hat{A}_n(u)$ is free: we will perform the coupling by choosing the form of each $\hat{A}_n(u)$ in a suitable way.
Without loss of generality, we can assume that $\cup_{u \in E_n^0} A_n(u)$ is an interval with 0 as left bound (see Figure \ref{fig:coupling}).

For $u \in E_n^0$, we can find a subinterval of $A_n(u)$ of length $p_n(u) \wedge \hat{p}_n(u)$, with the same left bound as $A_n(u)$,
and set for $r$ in this interval $\hat{f}_n(r)=u$. Hence,
${\mathbb P}\{f_n=u,\hat{f}_n \not = u|{\mathcal F}_n\} \leq (p_n(u)-\hat{p}_n(u))^+$, so that
\begin{equation}
\label{expr_couplage_E0}
{\mathbb P}\{f_n \in E_n^0,f_n \not = \hat{f}_n|{\mathcal F}_n\} \leq \sum_{u \in E_n^0} [
p(X_n,\xi_n,u) - p(\hat{X}_n,\xi_n,u)]^+.
\end{equation}
This is exactly what we were seeking for.

It remains to choose $\hat{f}_n$ such that
$\{f_n=-\hat{f}_n,f_n \in E_n^-,|\langle \Delta_n,f_n\rangle| =1\}$ is empty.
For $n < \sigma$, $\Delta_n \not =0$ and $E_n^{+}$ cannot be empty. Since
$|\Delta_n|_1$ is always even, $E^{+}_n$ cannot count one single vector
such that $|\langle \Delta_n,u \rangle|$ is odd. Hence, the set
$\{u \in {\mathcal V}: \, |\langle \Delta_n,u \rangle| \in 2{\mathbb N}+1\}$ counts
either zero element or more than two. If $M = |E_n^+| = 1$,
the set $\{u \in {\mathcal V}: \, |\langle \Delta_n,u \rangle| \in 2{\mathbb N}+1\}$ is empty and there is nothing to do.
If $M \geq 2$,
we can index $E_n^{+}$ under the form $E_n^+=\{v_1,\dots,v_M\}$ with $p_n(v_1) \geq p_n(v_2) \geq \dots \geq p_n(v_M)$.
Then, we can assume that the partition
related to $f_n$ is ordered as follows:
\begin{equation*}
\forall u \in E_n^0, \ A_n(u) \prec A_n(-v_1) \prec A_n(v_1)
\prec A_n(-v_2) \prec A_n(v_2) \dots \prec A_n(-v_M) \prec A_n(v_M),
\end{equation*}
where $B_1 \prec B_2$ means $\forall (x,y) \in B_1 \times B_2, \ x<y$,
$B_1$ and $B_2$ being two subsets of $[0,1]^2$ (see Figure \ref{fig:coupling}).

\begin{figure}[htb]
\begin{center}
\psfrag{0}{$0$}
\psfrag{1}{$1$}
\psfrag{en0}{$E_n^0$}
\psfrag{-v1}{$-v_1$}
\psfrag{-v2}{$-v_2$}
\psfrag{-v3}{$-v_3$}
\psfrag{v1}{$v_1$}
\psfrag{v2}{$v_2$}
\psfrag{v3}{$v_3$}
\includegraphics[
width=0.5\textwidth,angle=0]{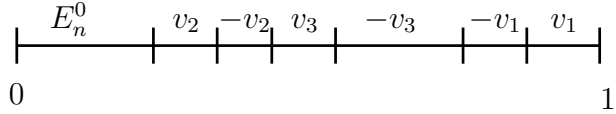}
\caption{Order for $\hat{f}_n$, $M=3$.}
\label{fig:coupling2}
\end{center}
\end{figure}

For $\hat{f}_n$, we already know that $\hat{A}_n(u)$ intersects, for $u \in
E_n^0$, $A_n(u)$ on an interval of length $p_n(u) \wedge \hat{p}_n(u)$. Then, we can complete $\hat{A}_n(u)$, if necessary, that is if
$\hat{p}_n(u)>p_n(u)$, so that $\cup_{u \in E_n^0} \hat{A}_n(u)$ is an
interval with zero as lower bound (see Figure \ref{fig:coupling2}). In particular, we have
\begin{equation*}
\forall u \in E_n^0, \ \forall v \not \in E_n^0, \ \hat{A}_n(u) \prec \hat{A}_n(v).
\end{equation*}
Then, we can complete the partition
associated to $\hat{f}_{n}$ as follows
\begin{equation*}
\hat{A}_{n}(v_2) \prec \hat{A}_n(-v_2) \prec
\hat{A}_{n}(v_3) \prec \hat{A}_{n}(-v_3)
\dots \prec \tilde{A}_{n}(v_M) \prec
\hat{A}_{n}(-v_M) \prec \hat{A}_{n}(-v_1)
\prec \hat{A}_{n}(v_1),
\end{equation*}
see Figure \ref{fig:coupling2}.

We now prove that, for $\delta$ small enough and $u \in E_n^{-}$, the sets $A_n(u)$ and $\hat{A}_n(-u)$ are
disjoint. For $2 \leq i \leq M$, the right boundary of $\hat{A}_n(v_i)$ is given by $\hat{p}_n(E_n^0) +
\hat{p}_n(v_2)+ \hat{p}_n(-v_2) + \dots + \hat{p}_n(v_i)$
and the left boundary of $A_n(-v_i)$ is given by $p_n(E_n^0)+ p_n(-v_1) + p_n(v_1) + \dots + p_n(-v_{i-1}) +
p_n(v_{i-1})$. By the Lipschitz property of $p$, the difference between $p_n(u)$ and $\hat{p}_n(u)$ is bounded
by $(C/m) |X_n -
\hat{X}_n|_1 \leq 2C \delta^{1/2}$ for every $u \in {\mathcal V}$. Since $p_n(v_1)\geq p_n(v_i)$, we have
\begin{equation*}
\begin{split}
&p_n(E_n^0) + p_n(-v_1)+ p_n(v_1) +\dots + p_n(-v_{i-1}) + p_n(v_{i-1})
\\
&\geq \hat{p}_n(E^n_0) +
\hat{p}_n(-v_2) + \hat{p}_n(v_2) + \dots + \hat{p}_n(-v_{i-1}) + \hat{p}_n(v_{i-1}) + \hat{p}_n(v_i) + p_n(-v_1) - 4 C d \delta^{1/2}
\\
&\geq \hat{p}_n(E^n_0) +
\hat{p}_n(-v_2) + \hat{p}_n(v_2) + \dots + \hat{p}_n(-v_{i-1}) + \hat{p}_n(v_{i-1}) + \hat{p}_n(v_i) + c - 4 C d \delta^{1/2},
\end{split}
\end{equation*}
with $c= \inf\{p(z,i,v); \, z \in [0,1]^d, \, i \in E, \, v \in{\mathcal V}\}>0$ (see Assumption {\bf (A.2)}).
For $\delta$ small enough, we obtain $\hat{A}_n(v_i) \cap A_n(-v_i) = \emptyset$. It remains to prove the same thing for $i=1$. The right boundary of
$A_n(-v_1)$ is given by $p_n(E_n^0) + p_n(-v_1)$ and the left boundary of $\hat{A}_n(v_1)$ is given by $\hat{p}_n(E_n^0) + \hat{p}_n(v_2) +
\hat{p}_n(-v_2) + \dots + \hat{p}_n(v_M) + \hat{p}_n(-v_M) + \hat{p}_n(-v_1) \geq p_n(E_n^0) + p_n(-v_1) + 2c - 4 C d \delta^{1/2}$.
This completes the construction of $\hat{f}_{n}$ for $M \geq 2$.
\medskip

{\bf Hitting Time.} Recall that ${\mathcal F}_n={\mathcal F}_n^{\xi,X,\hat{X}}$ for all $n \geq 0$.
 By \eqref{expr_couplage} and \eqref{expr_couplage_E0}, we have for $\delta$ small enough (say
 $\delta \leq \rho_0$ for some $\rho_0>0$) and $n<\sigma$
\begin{equation}
\label{expr_couplage_E0_2}
{\mathbb E}  \bigl[|\Delta_{n+1}|_1 | {\mathcal F}_n \bigr]
- |\Delta_n|_1
\leq  2 (p_n-\hat{p}_n)^+(E_n^0) + 2 p_n(E_n^+) - 2 \hat{p}_n(E_n^+),
\end{equation}
where $(p_n(\cdot)-\hat{p}_n(\cdot))^+(A) = \sum_{u \in A}(p_n(u)-\hat{p}_n(u))^+$
for any subset $A$ of ${\mathcal V}$ (the same holds for $p_n(A)$ and $\hat{p}_n(A)$).
By \eqref{hyp:tloi} in Theorem \ref{th:tloi}, we have
\begin{equation*}
\begin{split}
&2 (p_n - \hat{p}_n)^+(E_n^0)
+ 2 p_n(E_n^+) - 2 \hat{p}_n(E_n^+)
\\
&= 2 (p_n- \hat{p}_n)^+(E_n^0)
+ p_n(E^+_n) - \hat{p}_n(E_n^+)
\\
&\hspace{15pt}
+
\bigl(1-p_n(E^0_n)-p_n(E^-_n) \bigr)
- \bigl( 1 - \hat{p}_n(E^0_n) - \hat{p}_n(E^-_n) \bigr)
\\
&= 2 (p_n- \hat{p}_n)^+(E_n^0) - (p_n-\hat{p}_n)(E_n^0)
 + (p_n - \hat{p}_n)(E_n^+) -
(p_n - \hat{p}_n)(E_n^-)
\\
&=\sum_{u \in \Lambda, u  \perp \Delta_n} |p_n(u)-\hat{p}_n(u)|
+ \sum_{u \in \Lambda} \bigl(p_n(u)- \hat{p}_n(u)\bigr)
{\rm sgn}\bigl(\langle \Delta_n,u \rangle \bigr)
\\
&\leq -(\kappa/m) |\Delta_n|_1.
\end{split}
\end{equation*}
By
\eqref{expr_couplage_E0_2}, we can write
$
|\Delta_{n+1}|_1 = |\Delta_n|_1 + 2 \varepsilon_{n+1}
$,
with $\varepsilon_n \in \{-1,0,1\}$ and ${\mathbb
  E}(2\varepsilon_{n+1}|{\mathcal F}_n)\leq - (\kappa/m) |\Delta_n|_1$
  for $n < \sigma$.
Hence, ${\mathbb E}(|\Delta_{n+1}|_1|{\mathcal F}_n) \leq
  (1-\kappa/m) |\Delta_n|_1$ for $n < \sigma$, so that
  $((1-\kappa/m)^{-n \wedge \sigma}|\Delta_{n \wedge \sigma}|_1)_{n \geq 0}$
  is a supermartingale. We deduce that, for all $n \geq 1$,
\begin{equation*}
2 {\mathbb P}\{\sigma >n\} (1-\kappa/m)^{-n} \leq {\mathbb
  E}[(1-\kappa/m)^{-\sigma \wedge n} |\Delta_{\sigma
\wedge n}|_1] \leq 2 \delta m.
\end{equation*}
We obtain
\begin{equation}
\label{couplage_ineq1}
{\mathbb E}(\sigma) \leq  \delta m \sum_{n \geq 0}
(1-\kappa/m)^n \leq C' \delta m^2,
\end{equation}
for some constant $C'>0$.

We now investigate ${\mathbb P}\{|\Delta_{\sigma}|_1 > 2 \lfloor m
\delta^{1/2} \rfloor \}$. Since
  $((1-\kappa/m)^{-n \wedge \sigma}|\Delta_{n \wedge \sigma}|_1)_{n \geq 0}$
  is a supermartingale, we have for all $n \geq m^{3/2}$,
\begin{equation*}
(1-\kappa/m)^{-m^{3/2}} {\mathbb E} \bigl[ {\mathbf 1}_{\{\sigma \geq m^{3/2}\}} |\Delta_{\sigma \wedge n}|_1
\bigr] \leq {\mathbb E} \bigl[ (1-\kappa/m)^{-\sigma \wedge n} |\Delta_{\sigma \wedge n}|_1 \bigr] \leq 2 \delta m.
\end{equation*}
Letting $n$ tend to $+ \infty$, we deduce (changing if necessary the value of $C'$)
\begin{equation}
\label{couplage:ineq2}
(1-\kappa/m)^{-m^{3/2}} {\mathbb P}\bigl\{\sigma \geq m^{3/2},   |\Delta_{\sigma}|_1 > 2 \lfloor m\delta^{1/2} \rfloor\bigr\}
\leq C'.
\end{equation}
It remains to see what happens for $\sigma < m^{3/2}$.
We set
$q_n^+={\mathbb P}\{\varepsilon_n=1|{\mathcal F}_n\}$ and
$q_n^-={\mathbb P}\{\varepsilon_n=-1|{\mathcal F}_n\}$.
Conditionally to the past, the process $(|\Delta_k|_1)_{k \geq 0}$ doesn't move at time $n$ with probability $1-(q_n^++q_n^-)$. Conditionally to moving, it jumps with probabilities
$q_n^+/(q_n^++q_n^-)$ and $q_n^-/(q_n^++q_n^-)$. Since ${\mathbb E}(\varepsilon_{n+1}|{\mathcal F}_n)<0$, we have $q_n^+/(q_n^++q_n^-)<1/2$. Hence, the time needed by the chain
$(|\Delta_n|_1)_{n \geq 0}$ to reach $2 \lfloor m\delta^{1/2} \rfloor$ is (stochastically) larger than the time needed by the simple random walk to hit
 $2 \lfloor m\delta^{1/2} \rfloor - 2 \lfloor m \delta  \rfloor$ when starting from zero. Hence,
\begin{equation*}
 {\mathbb P}\bigl\{\sigma < m^{3/2},   |\Delta_{\sigma}|_1 > 2 \lfloor m\delta^{1/2} \rfloor\bigr\}
 \leq {\mathbb P} \{ \tau_{2 \lfloor m\delta^{1/2} \rfloor - 2 \lfloor m \delta  \rfloor} < m^{3/2} \},
\end{equation*}
 where $\tau_L$ denotes the hitting time, by the simple random walk, of a given integer $L$.
 It is well known (see e.g. \cite[Chapter 10]{Williams}) that ${\mathbb P}\{\tau_L \leq m^{3/2}\} \leq \exp(-\alpha L + m^{3/2} \ln(\cosh(\alpha)))$
 for any $\alpha >0$. Choosing
 $\alpha =m^{-3/4}$, we have ${\mathbb P}\{\tau_L \leq m^{3/2}\} \leq \exp(-L m^{-3/4} + m^{3/2} \ln(\cosh(m^{-3/4})))$. If
 $L = \eta m$, for some $\eta >0$,
 ${\mathbb P}\{\tau_{\eta m} \leq m^{3/2}\} \leq C' \exp(- \eta m^{1/4})$. Hence (changing $C'$ if necessary),
\begin{equation}
\label{couplage:ineq3}
{\mathbb P}\bigl\{\sigma < m^{3/2},   |\Delta_{\sigma}|_1 > 2 \lfloor m\delta^{1/2} \rfloor\bigr\}
\leq C' \exp(-2(\delta^{1/2}-\delta)m^{1/4}).
\end{equation}

We can complete the proof of Lemma \ref{couplage}. We have, for all $t \geq S$,
\begin{equation*}
|{\mathbb P}_x\{\tau^{(m)} > m^2 t \} - {\mathbb P}_y\{ \tau^{(m)} > m^2 t\}|
\leq 2 \P \{ \sigma >m^3 S \ {\rm or} \ |\Delta_{\sigma}|
\not = 0\}.
\end{equation*}
By \eqref{couplage_ineq1}, $\P\{\sigma >m^3S\} \leq C' \delta S^{-1} m^{-1}$. By \eqref{couplage:ineq2} and \eqref{couplage:ineq3},
$\P \{  |\Delta_{\sigma}| \not = 0\} \leq C' [ \exp(-2(\delta^{1/2}-\delta)m^{1/4}) + \exp(-\kappa m^{1/2})]$.
This completes the proof.

\subsection{Proof of Lemma \ref{lem:convex}}
\label{sec:convex}
For a given $x \in [0,1]^d$, we have to prove that the bilinear form $\nabla^2_{\alpha,\alpha} H : \lambda \in \R^d \mapsto \sum_{i,j} \lambda_i \lambda_j [\partial^2
H/\partial \alpha_i \partial \alpha_j](x,0)$ is positive definite.
We first note that the bilinear form ${\mathcal E}_f : \lambda \in \R^d \mapsto \sum_{i,j=1}^d \sum_{k \in E} \mu(k) \lambda_i \lambda_j
({\mathbb E}[(f_if_j)(x,k,U)] - {\mathbb E}[f_i(x,k,U)]{\mathbb E}[f_j(x,k,U)])$ induced by the averaged covariance matrix of the random vectors $f(x,k,U)$
($U$ following the uniform distribution on $(0,1)$) is nondegenerate. Indeed, for all $\lambda \in \R^d$, Jensen's inequality yields
\begin{equation*}
\begin{split}
{\mathcal E}_f(\lambda) &= \sum_{i=1}^d \sum_{k \in E} \mu(k) \lambda_i^2 {\mathbb E}(f_i^2(x,k,U)) - \sum_{k \in E} \mu(k) \bigl( \sum_{i=1}^d
\lambda_i {\mathbb E}(f_i(x,k,U)) \bigr)^2
\\
&=  \sum_{i=1}^d \sum_{k \in E} \lambda_i^2 \mu(k) p(x,k,\pm e_i)
 - \sum_{k \in E} \mu(k) \bigl( \sum_{i=1}^d \lambda_i [p(x,k,e_i)-p(x,k,-e_i)]\bigr)^2
\\
&\geq  \sum_{i=1}^d \sum_{k \in E} \lambda_i^2 \mu(k) p(x,k,\pm e_i)
 - \sum_{i=1}^d \sum_{k \in E} \mu(k) \lambda_i^2 \frac{[p(x,k,e_i)-p(x,k,-e_i)]^2}{p(x,k,\pm e_i)}
\\
&=  \sum_{i=1}^d \sum_{k \in E} \lambda_i^2 \mu(k) \frac{4 p(x,k,e_i)p(x,k,-e_i)}{p(x,k,\pm e_i)} >0,
\end{split}
\end{equation*}
with $p(x,k,\pm e_i)=p(x,k,e_i)+p(x,k,-e_i)$.
In what follows, we provide an explicit expression for $\nabla^2_{\alpha,\alpha} H$ and then compare it to ${\mathcal E}_f$.
We know that the leading eigenvalue of the matrix $Q(x,\alpha)$ (see \eqref{eq:perron}) is simple and equal to $\exp(H(x,\a))$. As a by-product, the coordinates of
the corresponding eigenvector $v(x,\alpha)$ (i.e. of the $\ell^1$ normalized eigenvector with positive entries) are infinitely differentiable with respect to $\alpha$. In particular, we
can differentiate twice the relationship $Q(x,\alpha) v(x,\alpha) = \exp(H(x,\alpha)) v(x,\alpha)$ with respect to $\alpha_i,\alpha_j$. We obtain
\begin{equation*}
\begin{split}
&\bigl[\frac{\partial Q}{\partial \alpha_i} v + Q \frac{\partial v}{\partial \alpha_i} \bigr](x,\alpha) = \bigl[ \bigl(\frac{\partial H}{\partial \alpha_i}
 v +  \frac{\partial v}{\partial \alpha_i} \bigr)\exp(H)\bigr](x,\alpha)
\\
&\bigl[ \frac{\partial^2 Q}{\partial \alpha_i \alpha_j} v +
\frac{\partial Q}{\partial \alpha_i} \frac{\partial v}{\partial \alpha_j}
+ \frac{\partial Q}{\partial \alpha_j}\frac{\partial v}{\partial \alpha_i}
+ Q \frac{\partial^2 v}{\partial \alpha_i \partial \alpha_j} \bigr] (x,\alpha)
\\
&\hspace{5pt} = \bigl[ \bigl(\frac{\partial^2 H}{\partial \alpha_i \partial \alpha_j}  v + \frac{\partial H}{\partial \alpha_i} \frac{\partial v}{\partial \alpha_j} +
\frac{\partial H}{\partial \alpha_i}\frac{\partial H}{\partial \alpha_j} v + \frac{\partial^2 v}{\partial \alpha_i \partial \alpha_j}
+ \frac{\partial v}{\partial \alpha_i} \frac{\partial H}{\partial \alpha_j} \bigr) \exp(H) \bigr](x,\alpha).
\end{split}
\end{equation*}
For $\alpha=0$, we know that
$Q(x,0)=P$ (so that $v(x,0)={\mathbf 1}=(1,\dots,1)^t$) and
$[\partial Q/\partial \alpha_i](x,0)
= (P_{k,k'}{\mathbb E}[f_i(x,k,U)])_{k,k' \in E}$. Hence, for every $k \in E$,
${\mathbb E}[f_i(x,k,U)]  + \sum_{k' \in E} P_{k,k'} [\partial v_{k'}/\partial \alpha_i](x,0) =
[\partial H/\partial \alpha_i](x,0)
+ [\partial v_k/\partial \alpha_i](x,0)$. Integrating with respect to the invariant measure $\mu$, we
deduce that
$[\partial H/\partial \alpha_i](x,0) = \sum_{k \in E} \mu(k) {\mathbb E}[f_i(x,k,U)] =\bar{f}_i(x)$
. Finally,
\begin{equation}
\label{eq:convex12}
\sum_{k' \in E} ({\rm Id} -P)_{k,k'} \frac{\partial v_{k'}}{\partial \alpha_i}(x,0) = {\mathbb E}[f_i(x,k,U)]
- \bar{f}_i(x).
\end{equation}
Applying the same method for the second order derivatives, we obtain for every $k \in E$:
\begin{equation}
\label{eq:convex2}
\begin{split}
&{\mathbb E}[(f_if_j)(x,k,U)]
+ \sum_{k' \in E} P_{k,k'} \frac{\partial^2 v_{k'}}{\partial \alpha_i \partial \alpha_j} (x,0)
\\
&\hspace{15pt} + \sum_{k' \in E} \bigl[ {\mathbb E}[f_i(x,k,U)] P_{k,k'}
- \bar{f}_i(x) \delta_{k,k'} \bigr] \frac{\partial v_{k'}}{\partial \alpha_j}(x,0)
\\
&\hspace{15pt}
+ \sum_{k' \in E} \bigl[ {\mathbb E}[f_j(x,k,U)] P_{k,k'} - \bar{f}_j(x) \delta_{k,k'} \bigr] \frac{\partial v_{k'}}{\partial \alpha_i}(x,0)
\\
&=\frac{\partial^2 H}{\partial \alpha_i \partial \alpha_j}(x,0) + [\bar{f}_i \bar{f}_j] (x)
+ \frac{\partial^2 v_k}{\partial \alpha_i \partial \alpha_j} (x,0).
\end{split}
\end{equation}
By \eqref{eq:convex12}, we have
\begin{equation*}
\begin{split}
&\sum_{k,k'  \in E} \mu(k) \bigl[ {\mathbb E}[f_i(x,k,U)] P_{k,k'} - \bar{f}_i(x) \delta_{k,k'} \bigr] \frac{\partial v_{k'}}{\partial \alpha_j}(x,0)
\\
&= \sum_{k,k' \in E} \mu(k) \bigl[{\mathbb E}[f_i(x,k,U)]   - \bar{f}_i(x) \bigr] P_{k,k'} \frac{\partial v_{k'}}{\partial \alpha_j}(x,0)
\\
&= \sum_{k,k',k'' \in E} \mu(k) ({\rm Id}-P)_{k,k''}  \frac{\partial v_{k''}}{\partial \alpha_i}(x,0)  P_{k,k'} \frac{\partial v_{k'}}{\partial \alpha_j}(x,0)
\\
&= \langle ({\rm Id}-P) \frac{\partial v}{\partial \alpha_i}(x,0),P \frac{\partial v}{\partial \alpha_j}(x,0) \rangle_{\mu},
\end{split}
\end{equation*}
where $\langle \cdot,\cdot \rangle_{\mu}$ denotes the scalar product on $L^2(\mu)$
and ${\rm Id}$ the identity matrix on $E$. We now integrate \eqref{eq:convex2} with respect to the invariant measure,
we deduce
\begin{equation}
\label{eq:convex3}
\begin{split}
\frac{\partial^2 H}{\partial \alpha_i \partial \alpha_j}(x,0)
&= \sum_{k \in E} \mu(k) {\mathbb E}[(f_if_j)(x,k,U) ] -  [\bar{f}_i \bar{f}_j ](x)
\\
&\hspace{5pt} + \langle ({\rm Id}-P) \frac{\partial v}{\partial \alpha_i}(x,0),P \frac{\partial v}{\partial \alpha_j}(x,0) \rangle_{\mu}
\\
&\hspace{5pt} + \langle ({\rm Id}-P) \frac{\partial v}{\partial \alpha_j}(x,0),P \frac{\partial
v}{\partial \alpha_i}(x,0) \rangle_{\mu}.
\end{split}
\end{equation}
Using \eqref{eq:convex12}, we deduce
\begin{equation}
\label{eq:convex4}
\begin{split}
&\sum_{k \in E} \mu(k) {\mathbb E}[f_i(x,k,U)] {\mathbb E}[f_j(x,k,U)]
 - [\bar{f}_i \bar{f}_j](x)
\\
&\hspace{15pt} = \langle ({\rm Id}-P) \frac{\partial v}{\partial \alpha_i}(x,0),({\rm Id}-P)
\frac{\partial v}{\partial \alpha_j}(x,0) \rangle_{\mu}.
\end{split}
\end{equation}
Plugging \eqref{eq:convex4} into \eqref{eq:convex3}, we obtain
\begin{equation*}
\begin{split}
\frac{\partial^2 H}{\partial \alpha_i \partial \alpha_j}(x,0)
&=\sum_{k \in E} \mu(k) {\mathbb E}[(f_if_j)(x,k,U)] - \sum_{k \in E} \mu(k)
 {\mathbb E}[f_i(x,k,U)] {\mathbb E}[f_j(x,k,U)]
\\
&+ \langle ({\rm Id}-P) \frac{\partial v}{\partial \alpha_i}(x,0),({\rm Id}-P) \frac{\partial v}{\partial \alpha_j}(x,0) \rangle_{\mu}
\\
&+ \langle ({\rm Id}-P) \frac{\partial v}{\partial \alpha_i}(x,0),P \frac{\partial v}{\partial \alpha_j}(x,0) \rangle_{\mu}
+ \langle ({\rm Id}-P) \frac{\partial v}{\partial \alpha_j}(x,0),P \frac{\partial v}{\partial \alpha_i}(x,0) \rangle_{\mu}
\end{split}
\end{equation*}
For all $(\lambda_i)_{1 \leq i \leq d} \in \R^d$ and $k\in E$, we set $u^{\lambda}(k) =
\sum_{i=1}^d \lambda_i [\partial v_k/\partial \alpha_i](x,0)$. Then,
\begin{equation*}
\begin{split}
\nabla^2_{\alpha,\alpha} H(\lambda)
&= {\mathcal E}_f(\lambda)
+ \langle ({\rm Id}-P)u^{\lambda} ,(I-P) u^{\lambda} \rangle_{\mu}
+ 2 \langle ({\rm Id}- P) u^{\lambda},P u^{\lambda} \rangle_{\mu}.
\\
&= {\mathcal E}_f(\lambda) +   \langle u^{\lambda} ,u^{\lambda} \rangle_{\mu} -
 \langle P u^{\lambda},P u^{\lambda} \rangle_{\mu}
 \\
&= {\mathcal E}_f(\lambda) + \int_{E} d\mu(k) {\mathbb E}^k \bigl[ (u^{\lambda}(\xi_1))^2 \bigr] - \int_E d\mu(k)
\bigl[ {\mathbb E}^k u^{\lambda}(\xi_1) \bigr]^2 \geq
{\mathcal E}_f(\lambda).
\end{split}
\end{equation*}
This completes the proof.

{\small


\begin{thebibliography}{99}

\bibitem{Atar-Dupuis}
Atar, R., Dupuis, P.: Large deviations and queueing networks: methods for rate function identification.
{\it Stochastic Process. Appl.} {\bf  84}  255--296, 1999.

\bibitem{AZR}
Azencott, R., Ruget, G.:
M\'elanges d'\'equations diff\'erentielles et grands \'ecarts \`a la loi des grands nombres.
{\it Z. Wahrscheinlichkeitstheorie und Verw. Gebiete}  {\bf 38}  1--54, 1977.

\bibitem{AZ} Azencott, R.: Grandes d{\'e}viations et applications. In:
 Eighth Saint Flour Probability Summer School---1978,  pp. 1--176,
 Lecture Notes in Math., 774, Springer, Berlin, 1980.

\bibitem{baldi88}
Baldi, P.:
Large deviations and stochastic homogenization.
{\it Ann. Mat. Pura Appl.} {\bf  151}  161--177, 1988.

\bibitem{bardi:capuzzo}
Bardi, M., Capuzzo-Dolcetta, I.
Optimal control and viscosity solutions of Hamilton-Jacobi-Bellman equations.
Birkhäuser Boston, Inc., Boston, MA, 1997.


\bibitem{barles}
Barles, G.
Solutions de viscosit\'e des \'equations de Hamilton-Jacobi. Springer-Verlag, Paris, 1994.

\bibitem{lions:90}
Capuzzo-Dolcetta, I., Lions, P.-L. 
Hamilton-Jacobi equations with state constraints.  Trans. Amer. Math. Soc.,  318, 643--683, 1990.

\bibitem{CDS1}
Comets, F., Delarue, F., Schott, R.:
 Distributed Algorithms in an Ergodic  Markovian Environment.
{\it Random Structures and Algorithms} {\bf 30},
131-167, 2007.

\bibitem{DZ} Dembo, A., Zeitouni, O.:
 Large deviations techniques and applications. Second edition. Applications of Mathematics 38. Springer-Verlag, New York, 1998.

\bibitem{dupuis 87}
Dupuis, P.: Large deviations analysis of reflected diffusions and constrained stochastic approximation algorithms in convex sets.
{\it Stochastics}  {\bf 21}, 63--96, 1987.

\bibitem{dupuis:88}
Dupuis, P.:
Large deviations analysis of some recursive algorithms with state dependent noise.
{\it Ann. Probab.} {\bf 16}, 1509--1536, 1988.

\bibitem{Dupuis-Ellis}
Dupuis, P. and Ellis, R.S.: The large deviations principle
for a general class of queueing systems I.   {\it Trans. Amer. Math. Soc.}
{\bf 347} 2689--2751, 1995.

\bibitem{Dupuis-Ramanan}
Dupuis, P., Ramanan, K.:
A time-reversed representation for the tail probabilities of stationary
reflected Brownian motion.  {\it Stochastic Process. Appl.} {\bf  98} 253--287,
2002.
 
\bibitem{FK}
Feng, J., Kurtz T.:
Large deviations for stochastic processes, 2005.
{\texttt http://www.math.wisc.edu/~kurtz/feng/ldp.htm}

\bibitem{flajolet}
Flajolet, P.:
The evolution of two stacks in bounded space and random walks in a triangle.
Proceedings of FCT'86, LNCS 233, 325--340, Springer Verlag, 1986.

\bibitem{FV}
Freidlin, M., Wentzell, A.D.: Random perturbations of dynamical systems.
Grundlehren der Mathematischen Wissenschaften,
260. Springer-Verlag, New York, 1984.

\bibitem{GS1} 
Guillotin-Plantard, N., Schott, R.: Distributed algorithms with dynamic random transitions. 
{\it Random Structures and Algorithms}  {\bf 21}  371--396, 2002.  

\bibitem{GS2}
Guillotin-Plantard, N., Schott, R.: Dynamic random walks. Theory and applications. Elsevier B.V., Amsterdam, 2006.

\bibitem{GV}
Gulinsky, O., Veretennikov, A.: Large deviations for discrete-time processes
with averaging. VSP, Utrecht, 1993.

\bibitem{ignatiouk}
Ignatiouk-Robert, I.: Large deviations for processes with discontinuous
statistics.  {\it Ann. Probab.}  {\bf 33}  1479--1508, 2005.

\bibitem{ignatiouk01}
Ignatiouk-Robert, I.:
Sample path large deviations and convergence parameters. 
{\it Ann. Appl. Probab.} {\bf 11} 1292--1329, 2001. 

\bibitem{knuth} 
Knuth D.E.,
The art of computer programming, Vol. 1, Addison-Wesley, 1973.  

\bibitem{lions sznitman 1984}
Lions, P.-L., Sznitman, A.-S.:
Stochastic differential equations with reflecting boundary conditions.
{\it Comm. Pure Appl. Math.} {\bf 37} 511--537, 1984.

\bibitem{lions:85}
Lions, P.-L. Neumann type boundary conditions for Hamilton-Jacobi equations.
Duke Math. J.,  52, 793--820, 1985.

\bibitem{lions:88}
Lions, P.-L. Optimal control of reflected diffusion processes: an example of state constraints. In: 
Stochastic differential systems (Bad Honnef, 1985), 269--276, Lecture Notes in Control and Inform. Sci., 78, Springer, Berlin, 1986. 


\bibitem{louchard1} 
Louchard, G.: Some distributed algorithms revisited. 
{\it Commun. Statist. Stochastic Models}  {\bf 4}  563--586, 1995.  

\bibitem{louchard2}
Louchard, G., Schott, R.: Probabilistic analysis of some distributed algorithms.
{\it Random Structures and Algorithms}  {\bf 2}  151--186, 1991. 

\bibitem{louchard3} 
Louchard, G., Schott, R., Tolley, M., Zimmermann, P.: Random walks, heat equations and distributed algorithms algorithms. 
{\it Computat. Appl. Math.}  {\bf 53}  243--274, 1994.  

\bibitem{Maier}
Maier, R.:
Colliding stacks: a large deviations analysis.
{\it Random Structures Algorithms}  {\bf 2}  379--420, 1991.

\bibitem{OV}
Olivieri, E., Vares, M.E.: Large deviations and metastability.
Encyclopedia of Mathematics and its Applications, 100.
Cambridge University Press, Cambridge, 2005.

\bibitem{Williams} Williams, D.: Probability with Martingales.
Cambridge University Press, Cambridge, 1991.

\bibitem{yao} 
Yao, A.: 
An analysis of a memory allocation scheme for implementing stacks. 
{\it SIAM J. Comput}  {\bf 10}  398--403, 1981. 
 
\end{thebibliography}
\end{document}